%% file: main.tex
\documentclass[11pt]{article}
\usepackage[utf8]{inputenc}
\usepackage{geometry}
\usepackage{tikzit}
\input{codecomps.tikzstyles}

\usepackage{verbatim}
\usepackage[dvipsnames]{xcolor}

\usepackage{titlesec}

\titleformat{\section}
  {\normalfont\scshape\filcenter}{\thesection}{1em}{}
\titleformat{\subsection}
 {\normalfont\scshape\filcenter}{\thesubsection}{1em}{}
\titlespacing*{\section}
{0pt}{5pt plus 1pt minus 1pt}{5pt plus 1pt minus 1pt}
\titlespacing*{\subsection}
{0pt}{5pt plus 1pt minus 1pt}{5pt plus 1pt minus 1pt}

\input{preamble}

\title{\vspace{-1cm} Co-Decompositions and FO Model Checking}
\author{Benjamin Merlin Bumpus \and Guilherme C.~M.~Gomes \and Laura Larios-Jones \and Will J.~Turner}
\date{}

\begin{document}

\maketitle

\thispagestyle{empty}

\setcounter{page}{0}

\vspace{-0.75cm}
\begin{center}
\begin{minipage}{0.9\textwidth}
\small
\textbf{Abstract.} 
Tree decompositions describe large graphs as gluings of smaller component parts. Here we study their categorical duals, \textit{tree co-decompositions}, which construct large graphs via fibered products of small constituent parts. Co-decompositions are compressed representations of graphs or relational structures: the decomposed object can be exponentially larger than the space required to specify the co-decomposition. 

Our main contribution is an algorithmic meta-theorem for first order logic: given \textit{only} a co-decomposition of a $\sigma$-structure $G$ with $b$ bags, each of size at most $w$ and an FO sentence $\phi$, there is an algorithm which, \textit{without ever constructing or inspecting $G$}, operates on the co-decomposition and determines whether $G \models \phi$ in time $f_\sigma(w,|\varphi|)\cdot b$. We pair this algorithm with structural results on tree co-width, including its relationship to tree-width and a separation from twin-width. Furthermore, we characterize bounded path co-width using a graph grammar. In contrast to our logical meta-theorem, we prove that \textsc{$(s,t)$-Path} becomes NP-hard on co-decomposition inputs. These results demonstrate how categorical constructions can enable efficient computation on exponentially compressed structures. 
\end{minipage}
\end{center}
\vspace{2cm}

\pagenumbering{gobble} 

\newpage 

\tableofcontents 

\pagebreak

\pagenumbering{arabic}

\input{sections/intro}

\input{sections/codecomps}

\input{sections/path-hardness}

\input{sections/existsFO}

\input{sections/fo}

\input{sections/discussion}

\bibliographystyle{plain}
\bibliography{ref}

\appendix

\input{sections/relational_structures}
\input{sections/appendix_binary}
\input{sections/appendix_tree_vs_co}
\input{sections/appendix_3_build}
\input{sections/appendix_grammar}

\end{document}

%% file: codecomps.tikzstyles
\tikzstyle{black dot}=[fill=black, draw=none, shape=circle, scale=0.3]

\tikzstyle{c0}=[-, line width=0.05em]
\tikzstyle{c0 arrow}=[->, line width=0.05em]
\tikzstyle{c0 bold}=[-, line width=0.2em]
\tikzstyle{c1}=[-, draw={rgb,255: red,94; green,164; blue,255}, line width=0.05em]
\tikzstyle{c1 arrow}=[->, draw={rgb,255: red,94; green,164; blue,255}, line width=0.05em]
\tikzstyle{c2}=[-, draw={rgb,255: red,255; green,162; blue,23}, line width=0.05em]
\tikzstyle{c2 dashed}=[-, draw={rgb,255: red,255; green,162; blue,23}, dashed]
\tikzstyle{c2 arrow}=[->, draw={rgb,255: red,255; green,162; blue,23}, line width=0.05em]
\tikzstyle{d1}=[-, dashed, line width=0.05em]
\tikzstyle{d1 arrow}=[->, line width=0.05em, dashed]

%% file: preamble.tex
\usepackage{pgfgantt} 

\usepackage{pifont}
\usepackage{xcolor}
\definecolor{parchment}{RGB}{214, 204, 169}

\usepackage{extpfeil}
\usepackage{mdframed}
\usepackage{microtype}
\usepackage{verbatim}
\usepackage{graphicx}
\usepackage{rotating}
\usepackage{todonotes}
\usepackage{tkz-graph}
\usepackage{amsthm}
\usepackage{amsmath} 

\usepackage{svg}

\usepackage{mathtools}

\usepackage{thm-restate}

\usepackage{url}
\definecolor{linkColor}{RGB}{156,78,13}
\usepackage{hyperref}
\hypersetup{
    colorlinks = true,
    allcolors = linkColor 
    }

\usepackage[nameinlink,capitalize]{cleveref}
\usepackage{pslatex}

\usepackage{datetime}
\newdateformat{versiondate}{%
\THEMONTH\THEDAY}

\usepackage{algorithm}
\usepackage{algpseudocodex}

\usepackage{verbatim}

\usepackage{framed}

\usepackage{pgf,tikz,environ}
\usetikzlibrary{arrows}
\usetikzlibrary{cd}

\usepackage{xpatch} 

\usepackage{adjustbox}

\newtheoremstyle{zoltanstyle}
  {0.2em} 
  {0.2em} 
  {} 
  {} 
  {\bfseries} 
  {.} 
  {.5em} 
  {} 
\renewcommand{\qed}{\hfill\textbf{Qed.}}

\newcommand{\pname}[1]{{\sc #1}}
\theoremstyle{zoltanstyle}
\swapnumbers
\xpatchcmd\swappedhead{~}{.~}{}{}
\newtheorem{body}{}
\crefname{body}{Paragraph}{Paragraphs}
\numberwithin{body}{section}

\newtheorem{conjecture}[body]{Conjecture}
\newtheorem{corollary}[body]{Corollary}
\newtheorem{definition}[body]{}
\crefname{definition}{Paragraph}{Paragraphs}
\newtheorem{example}[body]{Example}
\newtheorem{addendum}[body]{Addendum}

\newtheorem{lemma}[body]{Lemma}

\newtheorem{proposition}[body]{Proposition}
\newtheorem{remark}[body]{Remark}
\newtheorem{observation}[body]{Observation}
\newtheorem{theorem}[body]{Theorem}
\newtheorem{intuition}[body]{Intuition}

\newtheorem{question}[body]{Question}

\colorlet{exampleColor}{linkColor!65}

\surroundwithmdframed{example}
\surroundwithmdframed{lemma}
\surroundwithmdframed{theorem}
\surroundwithmdframed{proposition}

\surroundwithmdframed{addendum}

\usepackage{pgf,tikz,environ}
\usepackage{quiver}
\usepackage{proof}
\usetikzlibrary{arrows}
\usetikzlibrary{arrows.meta,positioning,fit,calc}
\usetikzlibrary{cd}
\usetikzlibrary{positioning}
\usetikzlibrary{fadings}
\usetikzlibrary{decorations.pathreplacing}
\usetikzlibrary{decorations.pathmorphing}
\tikzset{snake it/.style={decorate, decoration=snake}}

\newcommand{\xmark}{\text{\ding{55}}}
\newcommand{\mycmark}{\text{\ding{51}}}

\definecolor{parchment}{RGB}{214, 204, 169}

\makeatletter
\renewenvironment{proof}[1][\proofname]{\par
  \pushQED{\qed}%
  \normalfont
  \topsep0pt \partopsep0pt 
  \trivlist
  \item[\hskip\labelsep
        \itshape
    #1\@addpunct{.}]\ignorespaces
}{%
  \popQED\endtrivlist\@endpefalse
  \addvspace{4pt} 
}
\makeatother

\expandafter\let\expandafter\oldproof\csname\string\proof\endcsname
\let\oldendproof\endproof

\renewenvironment{proof}[1][\proofname]{%
  \oldproof[\normalfont \bfseries #1.]%
}{\oldendproof}

\usepackage[inline]{enumitem}
\usepackage{mathptmx}

\usepackage[T1]{fontenc}
\usepackage[utf8]{inputenc}

\usetikzlibrary {positioning}
\usetikzlibrary{fadings}
\usetikzlibrary{decorations.pathreplacing}
\usetikzlibrary{decorations.pathmorphing}
\tikzset{snake it/.style={decorate, decoration=snake}}
\definecolor{gold}{RGB}{255,215,0}
\definecolor{softBlack}{RGB}{45, 47, 49}
\definecolor{creamWhite}{RGB}{245,244,241}
\definecolor{softGray}{RGB}{220, 216, 214}
\definecolor{brick}{RGB}{232, 48, 48}

\definecolor{goodred}{HTML}{CC6677}
\definecolor{goodblue}{HTML}{332288}
\definecolor{goodyellow}{HTML}{DDCC77}
\definecolor{goodgreen}{HTML}{117733}
\definecolor{goodcyan}{HTML}{88CCEE}
\definecolor{goodwine}{HTML}{882255}
\definecolor{goodteal}{HTML}{44AA99}
\definecolor{goodolive}{HTML}{999933}
\definecolor{goodpurple}{HTML}{AA4499}

\newcommand{\typesetoperator}[1]{\mathsf{#1}}

\usepackage{xspace}
\definecolor{gold}{RGB}{255,215,0}
\definecolor{softBlack}{RGB}{45, 47, 49}
\definecolor{creamWhite}{RGB}{245,244,241}
\definecolor{softGray}{RGB}{220, 216, 214}
\definecolor{brick}{RGB}{232, 48, 48}

\usepackage{enumitem}

\newcommand{\define}[1]{\textbf{#1}}

\newcommand{\cat}[1]{\mathsf{#1}}

\DeclareMathOperator{\op}{^{\typesetoperator{op}}}

\newcommand{\typesetproblem}[1]{\textsc{#1}}

\DeclareMathOperator{\np}{\typesetproblem{NP}}

\DeclareMathOperator{\Da}{\mathcal{D}}  
\DeclareMathOperator{\Ea}{\mathcal{E}}  
\DeclareMathOperator{\Fa}{\mathcal{F}}  
  
\DeclareMathOperator{\Ha}{\mathcal{H}}

\DeclareMathOperator{\Pa}{\mathcal{P}}  
\DeclareMathOperator{\Qa}{\mathcal{Q}}

\DeclareMathOperator{\Ta}{\mathcal{T}}

\DeclareMathOperator{\Xa}{\mathcal{X}}  
\DeclareMathOperator{\Ya}{\mathcal{Y}}

\newcommand{\bigO}[1]{{\mathcal{O}\!\left(#1\right)}}

\newcommand{\pts}[1]{\underline{#1}}

%% file: sections/intro.tex
\section{Introduction}


\body{
    Whether mathematical concepts are discovered or invented, the path by which we arrive at them is not always the best route to understanding. We came to \textit{tree co-decompositions} by dualizing a categorical generalization of tree-width~\cite{StructuredDecompositions}, where the topological viewpoint is replaced by a more algebraic one.
    Yet the notion of a co-decomposition extends beyond graphs and tree-shapes: co-decompositions of vector spaces are precisely \textit{cellular sheaves}, which  play a growing role in applied topology~\cite{Curry2014,HansenGhrist2019,Robinson2017}. In this introduction, however, we choose a more natural way, at least in hindsight, and which hints at potential applications of our algorithms: the use of coordinates to organize distinct concerns in scientific modeling~\cite{LibkindEtAl2022Epidemic,BaezEtAl2023StockFlowFramework}.
}

\body{
    Suppose Letícia, an epidemiologist, is trying to understand the spread of a new disease; let's call it X. She knows that various factors influence its spread: infection dynamics, host population movement patterns, population demographics, other diseases, etc. She would prefer to model each of these aspects independently, making it is easier to reason about, communicate, and explain her findings. She sets up a disease spread model, a model to account for ageing, a model to record movement from one city to another, and a model for a separate disease, Y, which is known to interact with X. Of course, what Letícia is really after is a combined model which is appropriately stratified, in the sense that she can combine states across the various models: for instance, she might want a part corresponding to those individuals who live in São Paulo, are elderly, are exposed to X, and have recently recovered from~Y. This sort of situation arises in practice, and several mathematical frameworks have been developed for this purpose \cite{LibkindEtAl2022Epidemic,BaezEtAl2023StockFlow,BaezEtAl2023StockFlowFramework}. Some of these use \textit{fibered products}, also known as \textit{pullbacks}, to negotiate the relationships between the constituent models so that they may be assembled into a much larger model of the process Letícia actually cares about~\cite{LibkindEtAl2022Epidemic,BaezEtAl2023StockFlowFramework}. A co-decomposition is precisely such a ``recipe'' for assembling large models as fibered products of smaller pieces.
}\label{para: leticia}

\body{
    As the number of variables of the combined model increases, inference and information retrieval becomes exponentially more difficult. In practice, Letícia might cope with this intractability by describing the individual models as small combinatorial structures, which interact through what we call \emph{adhesions}. For example, she may use Petri nets to organize modelling at a high level, since the corresponding systems of differential equations can be derived automatically from the nets once the appropriate rates have been specified \cite{LibkindEtAl2022Epidemic}. Suppose abstractly that Letícia has to combine two models $G$ and $H$. She does this by first selecting an adhesion $K$, together with homomorphisms $g\colon G\to K$ and $h\colon H\to K$. She then constructs the larger model, $G\times_K H$, by taking pairs of elements of $G$ and $H$ mapped to the same point in $K$. We state this formally for graphs, but there is an analogous definition for arbitrary relational structures.
}

\begin{figure}
    \centering
    \includegraphics[width=0.99\linewidth, page=10]{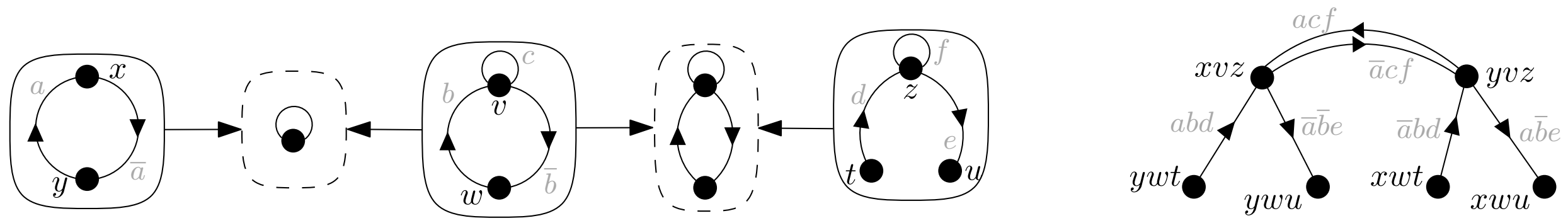}
    \caption{Left: A tree co-decomposition of graphs whose auxiliary graph is the path with two edges. The two bags on the right are drawn such that the vertices are horizontally aligned with the vertices of the adhesion to which they are mapped. That is, $v$ and $z$ are mapped to the top vertex in the adhesion. Right: The limit of the co-decomposition. For ease of reading, edges are labeled in gray and vertex labels are black.}
    \label{fig:limit}
\end{figure}

    \begin{definition}\label{def:pullback}
        For graphs $G,H,K$ and homomorphisms $g \colon G \to K$, $h \colon H \to K$, the \define{pullback} (or fibered product) of $G$ and $H$ along $K$ is the graph $G \times_K H$ with vertex set $V:=\{(x,y)\in V(G)\times V(H) \mid g(x)=h(y)\}$ and edge set $E:=\{((x,y),(x',y'))\in V\times V \mid (x,x')\in E(G)\wedge (y,y')\in E(H)\}$.
    \end{definition}
\body{
    The key insight here is that the relations (edges) present in $G\times_K H$ are completely determined, coordinate-wise, by the relations in $G$ and $H$, subject to the compatibility imposed by $K$. By repeating this operation for each of the smaller models she has constructed, Letícia can systematically assemble the large model she is after. If the relationships between Letícia's small models describe a tree, then the combined model is what we call a \emph{tree co-decomposition}, an example of which can be seen in the left-hand side of \cref{fig:limit}. 
}\label{para: pullback Leticia}
\newpage
\begin{definition}\label{def:tree_codecomp}
    A \define{tree co-decompsition} $\Ta$ is a triple $(T, \{G_a\},\{g_t^e\})$ such that: \vspace{-.5em}
    \begin{enumerate}
        \item $T$ is a tree;\vspace{-.5em}
        \item the family of graphs $\{G_a\}:=\{G_e\}_{e \in E(T)}\cup\{G_t\}_{t \in V(T)}$ is indexed over \textit{both} the edges $e\in E(T)$ and nodes $t\in V(T)$ of $T$; and \vspace{-.5em}
        \item each $\{g_t^e\}_{t\in V(T),e\in E(t)}$ is a family of an epimorphisms $G_t\to G_e$ indexed over the incident edges $e\in E(t)$ of each node $t\in V(T)$.\vspace{-.5em}
    \end{enumerate} 
    We call $T$ the \define{auxiliary tree}. The graphs $\{G_t\}$ indexed by nodes are the \define{bags of $\Ta$} and the graphs $\{G_e\}$ indexed by edges are the \define{adhesions of $\Ta$}; the \define{width} of $\Ta$ is its maximum bag size.
    We say that $G$ is the \define{limit} of $\Ta$ if $G$ is isomorphic to the graph obtained by recursively taking the pullbacks along the edges of $T$, i.e., $G \simeq \lim \Ta$.
    The \define{tree co-width} of $G$, denoted by $\mathsf{cow}(G)$, is the minimum among the widths of all tree co-decompositions whose limit is $G$.
    When the tree $T$ is a path, we call $\Ta$ a \define{path co-decomposition} (these admit an alternative characterization in terms of a grammar; c.f.~Appendix~\ref{appendix: grammar}).
\end{definition}

\body{
    As our choice of name and the tree shape mentioned above suggest, the central object of this paper is deeply connected to tree-width. Tree decompositions first appeared in the work of Bertel{\`e} and Brioschi~\cite{BerteleBrioschi1972}, were independently studied by Halin~\cite{Halin1976}, and were later reintroduced by Robertson and Seymour~\cite{ROBERTSON198449}, who used them to great effect to characterize graphs with a forbidden large clique minor in the Graph Minors project. Since then, they have become indispensable tools in combinatorics and   algorithmics~\cite{BODLAENDER06,CyganLPPS14,LOKSHTANOVMS2018,Courcelle_Engelfriet_2012,langer_practical_2014}. A \emph{tree decomposition} of a graph~$G$ can be thought of as a collection of subgraphs of $G$ (the \textit{bags}), indexed by the nodes of an auxiliary tree, which together cover $G$. The overlaps between these bags are controlled by the tree in the sense that, if a vertex belongs to two bags, it must also belong to every bag along the unique path connecting them. The important point is that $G$ can be reconstructed by gluing neighboring bags along their shared subgraphs (the \textit{adhesions}). In this paper, we reverse this picture: with tree \textbf{co}-decompositions, instead of embedding adhesions into bags and gluing these together, we map bags onto adhesions and reconstruct the global graph from compatible choices of coordinates; we call these \textit{matching families}.
}

\body{
    The reversal described above, in which bags map to adhesions rather than embedding adhesions into bags, is an instance of \textit{categorical duality}.
    Indeed, in recent years, tree decompositions have received several category-theoretic treatments and generalizations, including through graph grammars and cospan decompositions~\cite{BauderonCourcelle1987,BlumeEtAl2013}, monoidal width~\cite{DiLavoreSobocinski2023}, spined categories~\cite{BumpusKocsis2023,BumpusThesis2021}, and, most relevantly to this paper, the framework of \textit{structured decompositions}~\cite{StructuredDecompositions}. In part, the present paper originates from asking what happens when we dualize the latter construction\footnote{Structured decompositions~\cite{StructuredDecompositions} are common generalizations of: (1) tree decompositions, (2) Bass-Serre `graphs of groups'~\cite{BASS19933, serre1977arbres, serre2002trees, serre1970groupes}, and (3) hybrid dynamical systems~\cite{ames2006categorical}. A co-decomposition in a category $\cat C$ is a structured decomposition in $\cat C^{\mathrm{op}}$. Here we study tree co-decompositions of graphs and relational structures.}: replacing inclusions of adhesions into bags with epimorphisms from bags onto adhesions, and replacing pushouts with pullbacks. Although tree co-decompositions are categorical duals of tree decompositions, the associated width parameters behave rather differently.
}    
    \begin{proposition}[Appendix~\ref{appendix:tree-co-vs-tree}]\label{prop: different from tree-width}
        For every graph $G$, $\mathbf{cow}(G) \leq \mathbf{tw}(G) + 2$. Moreover, the class of bounded tree-width graphs is strictly contained in the class of bounded tree co-width graphs. 
    \end{proposition}

\body{
Another important graph parameter related to tree co-width is \emph{twin-width}, introduced by Bonnet, Kim, Thomassé and Watrigant~\cite{TWINWIDTH1} as a generalization of a permutation invariant studied by Guillemot and Marx~\cite{GuillemotMarx}, and which has received considerable attention in the literature~\cite{BonnetGKTW2024, gajarsky_et_al:LIPIcs.ICALP.2022.123, bonnet_et_al:LIPIcs.STACS.2023.15, heinrich_et_al:LIPIcs.WG.2026.26, ganian_et_al:LIPIcs.STACS.2026.42, bonnet_et_al:LIPIcs.SoCG.2025.23,pettersson2023bounds} due to its structural and algorithmic properties; in particular, we highlight the fact that first-order model checking is FPT-time solvable when parameterized by twin-width and formula length~\cite{TWINWIDTH1}.
In fact, since direct products of graphs are special cases of tree co-decompositions, we obtain the following separation from twin-width.

\begin{theorem}\label{prop: different from twin-width}
    There are graphs of arbitrarily large twin-width whose \textit{path} co-width is at most $3$.\footnote{Let $G_d$ be the $d$-fold direct product of $K_3$ with itself. As each factor maps onto the terminal graph, $\mathbf{cow}(G_d)\leq 3$. However, by \cite[Proposition~3.6]{pettersson2023bounds}, we have $\operatorname{tww}(G_d)=2^d$ for $d\geq 2$. Note that tree co-width is upper-bounded by path co-width.}
\end{theorem}
}

\body{
    Another way of understanding tree co-decompositions is as a \textit{succinct} representation of relational structures. Indeed, a co-decomposition with $b$ bags, each of size at most $w$, can give rise to a structure with as many as $w^b$ points. Given the enormous gap between the size of the representation and that of the object being constructed, the co-decompositions \textit{themselves} will be the objects of our study. In particular, all of our algorithmic results assume that we are given a co-decomposition of a graph $G$, but that \textit{we are neither given $G$ nor allowed to construct it}. Thus, all of our algorithmic results are only \textit{implicitly} about~$G$. There is a sense in which this is the most natural setting to consider; after all, Letícia's large, composite models are always given as co-decompositions: she wanted to be able to see all the constituent parts, point to them in scientific discussions and tinker with them individually, but infer global information.
}

\body{
    Since the limit $G$ of a co-decomposition $\Ta$ with $b$ bags of width at most $w$ can be so large, an algorithm that efficiently computes the properties of $G$ directly from $\Ta$ could be exponentially faster than one that takes $G$ directly. As appealing as this is, we might also expect this to mean that some simple problems are made intractable by this succinct representation. As with the work of Galperin and Widgerson~\cite{GALPERIN1983183} on another representation of size $o(|V(G)|)$, this is indeed the case: although \pname{$(s,t)$-Path} is efficiently solvable when given $G$ explicitly, the problem becomes $\np$-hard when we are given only a tree co-decomposition of $G$ and two matching families (i.e., the coordinate-wise descriptions of two vertices of $G$).
}

\begin{theorem}\label{thm: s-t path hardness}
    Given only a tree co-decomposition and two matching families, \pname{$(s,t)$-Path} is $\np$-hard.
\end{theorem}

\body{
    \cref{thm: s-t path hardness} portends a rough road ahead for algorithm designers and, specifically, returning to our story, to Letícia. Recall that she wanted to study a large model only via its building blocks and their relationships, i.e., via the the co-decomposition. Although each component makes sense individually, she is unsure whether the resulting global model behaves as intended. The issue is that the entire model is too large to reason about directly, and running simulations may be prohibitively expensive. So she has an idea: perhaps she can determine whether certain structural patterns, such as particular kinds of feedback loops, occur in the global model just by looking at the co-decomposition. Fortunately, in contrast to \cref{thm: s-t path hardness}, it turns out there are many properties of the global structure that can be decided efficiently by working \textit{only} with the co-decomposition. As our main algorithmic contribution, we prove that we can decide if the limit of a tree co-decomposition is a model for a first-order formula in FPT-time when parameterized by the formula length and the width of the co-decomposition. The running time of our algorithm is rather galactic (c.f. \cref{para: explicit bounds}). However, we get a much faster -- runnable -- algorithm for $\exists$-FO (c.f. \cref{thm:exists-fo}). 
    }

\begin{theorem}
    \label{thm: FO model checking on codecompositions}
    Fix a finite relational signature $\sigma$. There is an algorithm that, given an FO sentence $\phi$ and a tree co-decomposition $\Ta$ of $\sigma$-structures of width at most $w$ with $b$ bags, decides whether $\lim\Ta\models\phi$ in time $f_\sigma(w,|\phi|)\cdot b$, for some computable $f_\sigma$.
\end{theorem}

\body{
    From the algorithms designer's point of view, \cref{thm: FO model checking on codecompositions} immediately implies membership in FPT for several problems, such as \pname{$k$-Clique}, \pname{$k$-Matching}, and \pname{$k$-Dominating Set} and several others when co-width and solution size are the parameters. In forthcoming work, we will show that the optimization variants of several such particular problems also lie in FPT when parameterized by co-width only.
}

\body{
    The proof \cref{thm: FO model checking on codecompositions} is obtained by means of a delicate dynamic program. However, writing such an algorithm on a co-decomposition $\Ta$ feels rather different from doing so on a tree decomposition. For one thing, as mentioned earlier, we can explore $\Ta$, but are never allowed to see the decomposed limit. For another, each bag of $\Ta$ corresponds to a coordinate axis, so there is a sense in which, at any given bag, we are seeing a projection of the \textit{entire} global object. This gives dynamic programming on co-decompositions a somewhat Bayesian feel: as we traverse the tree, we discover new coordinates, as though shifting the position of a lamp casting shadows of the global object. Each new shadow may reveal previously hidden structure, or show that what we had conjectured from earlier shadows could not have been true globally. These shadows are well behaved in the following sense:
    \begin{enumerate*}
        \item a relation that fails in any one of the coordinates (bags) cannot hold globally; and
        \item a relation that holds in every coordinate we have inspected so far remains, at best, a conjecture until we have traversed the entirety of $\Ta$.
    \end{enumerate*}
    This motivates the central abstractions of our algorithm, which we call \textit{hypotheses}: small relational structures that probe the truth of atomic statements in the global structure. We collect hypotheses at each bag and endeavour to disprove them whenever we encounter a new ``shadow'' or coordinate axis. As we explore more bags, we combine these hypotheses and discard those that are no longer \textit{supported} (in a technical sense) by the coordinates we have encountered.
}\label{para: bayesian}

\body{
    The intuitive outline we have just given suffices to show that model checking for $\exists$-FO lies in FPT when parameterized by co-width and formula length. The case of full FO, allowing both universal and existential quantification, is significantly more challenging and most of the paper is devoted to it. Our approach is to turn to the familiar game-theoretic interpretation of quantification in which one looks for winning strategies for the existential player, or equivalently, Skolem-like instructions specifying how existential choices depend on preceding universal choices. Our algorithm tracks the projections of candidate strategies onto the individual bags, together with hypotheses recording which atomic relations may hold according to a strategy. We then combine this information using pullbacks, discarding unsupported possibilities as we proceed. 
}

\body{
    From a bird's-eye view, the proof of \cref{thm: FO model checking on codecompositions} retains the same structure\footnote{The relationship between the two algorithms is summarized via a ``table of translations'' (\cref{tab:existential-to-full-fo}).} as the much easier case of $\exists$-FO. Although the proof of the case of full FO is substantially more category-theoretic, much of the paper can be read without categorical pre-requisites.  The structural properties of co-decompositions (Section~\ref{sec: tree co-decompositions}, and Appendices~\ref{appendix:tree-co-vs-tree}, \ref{appendix: grammar}, \ref{appendix:binary}, \ref{appendix: tree shape important}) and or our hardness results (Section~\ref{sec: hardness}) are developed purely graph-theoretically. Similarly, our algorithm for $\exists$-FO model checking (Section~\ref{sec: exists FO}) only requires a sprinkle of model-theoretic terminology and familiarity with dynamic programming. The proof of \cref{thm: FO model checking on codecompositions} (Section~\ref{sec: full FO}) does need more category theory and indeed we do not know how to prove it without this perspective. 
}

\subsection{Related Work}

\body{
    \cref{thm: FO model checking on codecompositions} is an example of an \emph{algorithmic metatheorem}: a result that proves the (in)tractability of large families of problems by using a unifying language to describe them~\cite{grohe2009methods}; we point to~\cite{siebertz2024advancesalgorithmicmetatheorems} for a recent survey on the subject.
    These metatheorems often rely on proving the efficient decidability of some fragment of second-order (SO) logic for a broad graph class; for example, Courcelle's celebrated theorem~\cite{COURCELLE199012}, shows that small \emph{monadic} SO$_2$ formulae can be efficiently verified on graphs of tree-width at most a fixed constant $k$.
    Formula satisfability and SO are by no means the only targets of these metatheorems, which encompass several tasks with varying objectives and constraints~\cite{DELIMA2026103819,COURCELLE20092675,ARNBORG1991308,jelinek2025monadicsecondorderlogicpermutations,MIMWIDTH,METAKERNELIZATION}.
}

\body{
    We are interested in the landscape around FO model checking which, despite having limited expressivity~\cite{Hella1995,IMMERMAN198276}, models a many problems and has a rich literature~\cite{FrickGrohe2001,DvorakKralThomas2013,GroheKreutzerSiebertz2017,DreierMaehlmannSiebertz2023,DreierEtAl2024MonadicallyStable,DreierMaehlmannTorunczyk2024,FrickGrohe2002,GajarskyEtAl2024Elementary}. In~\cref{eg:half_graphs}, we build all (symmetric) half-graphs as limits of path co-decompositions of width at most four, proving that the classes of graphs of bounded tree co-width are neither structurally-nowhere-dense nor monadically stable.
}

\body{
    Other related topics are Feferman-Vaught composition~\cite{FefermanVaught1959} and its algorithmic uses, studied by Makowsky~\cite{Makowsky2004}. Rather than auxiliary formulas or logical types, our dynamic program tracks projections of candidate global assignments and strategies, together with hypotheses recording their relational information\footnote{See Section~\ref{sec: discussion} for more on the relationship between classical Feferman-Vaught-like results and ours.}. Witness-oriented precedents include Kneis, Langer and Rossmanith's game-theoretic proof of Courcelle's theorem~\cite{KneisLangerRossmanith2011} and Jakl, Marsden and Shah's categorical account of logical composition via game comonads~\cite{JaklMarsdenShah2025}. We also build on the framework of structured decompositions developed by Bumpus et al.~\cite{StructuredDecompositions} and categorical algorithms of Althaus et al. for deciding sheaves~\cite{AlthausEtAl2023Sheaves} and limit emptiness~\cite{AlthausEtAl2026InLim}. Our construction uses pullbacks both to assemble structures and to reason about them compositionally in categories of witnesses.
}

\body{
    Another closely related line of work concerns vector-space-valued \emph{cellular sheaves} on graphs and simplicial complexes, studied by Curry~\cite{Curry2014}, Hansen and Ghrist~\cite{HansenGhrist2019}, Hansen and Gebhart~\cite{HansenGebhart2020SheafNeuralNetworks}, and Robinson~\cite{Robinson2017}. Their global sections are precisely the limits of the underlying diagrams of linear maps; on trees, those with surjective restriction maps are tree co-decompositions valued in vector spaces. Their applications in topology, spectral graph theory, sensor integration and sheaf neural networks further motivate algorithms on diagrams and their limits.
}\vspace{-.5em}

\subsection{Outlook}
\body{
    We view this paper as further evidence for a broader programme of \textit{categorical algorithmics}: the use of category theory to develop algorithmic ideas and, conversely, of algorithmic ideas to develop new category theory. This means studying categorical constructions both as objects of computation and as tools for algorithm design. After all, one of the difficulties of algorithmics is the tremendous amount of data and bookkeeping involved, and category theory has a long tradition of managing such complexity in geometry, topology and algebra. Why not in algorithmics? Universal constructions also arise in scientific practice (think of Letícia) where they uncover tacit, underlying patterns and new questions. Even our present contribution, which is relatively narrow in its scope, reveals many new research questions. We compile a dozen of these in Section~\ref{sec: discussion} and organize them by the fields they most naturally inhabit. Most interesting to us, however, is that the same constructions can appear on both sides of an algorithmic problem: Letícia uses pullbacks to assemble her complicated model from simpler ones, while we use pullbacks to assemble information about that model. We suspect there more is to be learned from this correspondence. 
}

%% file: sections/codecomps.tex
\section{Preliminaries: Pullbacks, Matching Families and Limits}\label{sec: tree co-decompositions}

\body{
    Throughout, all graphs are directed and allowed to contain loops, but not parallel edges. We describe a graph $G$ as a vertex set $V(G)$ along with an edge relation $E(G)\subseteq V(G)\times V(G)$. A \define{graph homomorphism} $\phi \colon G \to H$ is a function $\phi\colon V(G)\to V(H)$ so that if $(x,y)\in E(G)$, then $(\phi(x),\phi(y))\in E(H)$. We sometimes write the composition $g\circ f$ as $gf$. We think of undirected graphs $G$ as \textit{symmetric}, meaning that for each arc $(x,y) \in E(G)$ the arc $(y,x)$ is also in $G$. A \define{graph epimorphism} is a graph homomorphism which is surjective (as a function on vertex sets). See Hell and Ne{\v{s}}et{\v{r}}il~\cite{HellNesetril2004} for further background on homomorphisms. 
}

\body{
    A \define{category} $\cat{C}$ is, roughly, a directed, reflexive graph equipped with a binary operation on edges called `composition'. Formally, it consists of a family of objects, morphisms $f:X\to Y$ between objects, an identity morphism $\cat{id}_X \colon X\to X$ for each object, and an associative composition law $g \circ f$  whenever the codomain of $f$ equals the domain of $g$, such that $f\circ \cat{id}_X=f$ and $\cat{id}_Y\circ f=f$. A \define{functor} $F \colon \cat{C} \to \cat{D}$ `is' a graph homomorphism: formally, it is a map sending objects (resp. arrows) in $\cat{C}$ to objects (resp. arrows) in $\cat{D}$ which preserves identities ($F(\cat{id}_X) = \cat{id}_{F(X)}$) and composition ($F(f \circ g) = F(f) \circ F(g)$). For background on category theory, appearing mostly in Section~\ref{sec: full FO}, we recommend Riehl~\cite{Riehl2016} and Awodey~\cite{Awodey2010}, and refer to Jacobs~\cite{JacobsCLTT} for the more advanced material on fibrations and categorical logic used later in the paper.
}

\body{
    Throughout the paper, when dealing with structures, we work with a fixed signature $\sigma$. Given a structure~$G$, we write $\pts{G}$ for the \define{underlying point set} of $G$. For example, if $\sigma$ is the signature of graphs encoded as a binary edge relation on a vertex set, then $\pts{G}$ is just $V(G)$; however, if $\sigma$ encodes graphs as an incidence relation between edges and vertices, then $\pts{G}$ is the set $E(G)\sqcup V(G)$. See Appendix~\ref{appendix:RS} for details.
}\label{para: prelim}

\vspace{-.5em}
\begin{minipage}[t]{0.70\textwidth}
\body{
Our constructions are based on \textit{pullbacks} (also known as \textit{fibered products}). Recall that, in any category~$\cat{C}$, given a \define{cospan}, that is a pair of morphisms with common codomain $X \xrightarrow{f} Z \xleftarrow{g} Y$, its \define{pullback} is an object $X\times_Z Y$ equipped with projection morphisms $X \xleftarrow{\pi_X} X\times_ZY \xrightarrow{\pi_Y} Y$ satisfying $f\pi_X=g\pi_Y$, and universal with this property: for any other object $W$ equipped with morphisms $u:W\to X$ and $v:W\to Y$ satisfying $fu=gv$, there exists a unique morphism $W\to X\times_ZY$ commuting with the projections.
}
\end{minipage}
\begin{minipage}[t]{0.20\textwidth}
\[\hspace{1em}\begin{tikzcd}
	& W & \\
	X & {X\times_ZY} & Y \\
	& Z
	\arrow["u"', from=1-2, to=2-1]
	\arrow["{\exists!}"{description}, from=1-2, to=2-2]
	\arrow["v", from=1-2, to=2-3]
	\arrow["f"', from=2-1, to=3-2]
	\arrow["{\pi_X}"', from=2-2, to=2-1]
	\arrow["{\pi_Y}", from=2-2, to=2-3]
	\arrow["g", from=2-3, to=3-2]
\end{tikzcd}\]
\end{minipage}
\vspace{.4em}

\begin{figure}
    \begin{minipage}{0.58\textwidth}
        \includegraphics[width=\textwidth, page=8]{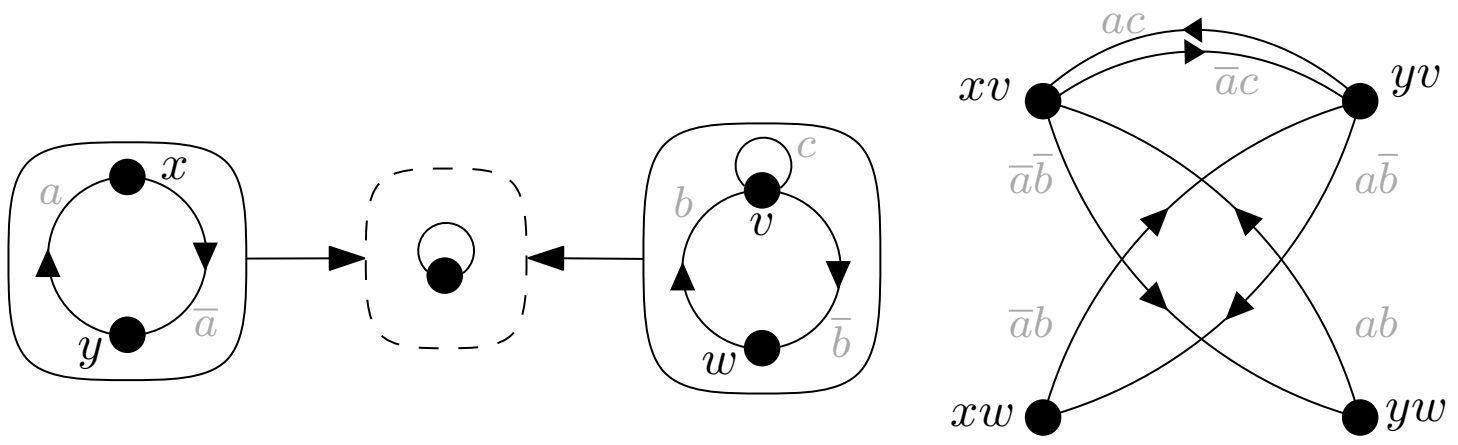}
    \end{minipage}%
    \hspace{\fill}
    \begin{minipage}{0.38\textwidth}
        \caption{Left: A cospan of graph epimorphisms, where the graph in the dashed box is the adhesion. Right: The cospan's pullback. Edge labels are in gray, vertex labels in black and no arrows on loops.}
    \label{fig:pullback}
    \end{minipage}%
\end{figure}

\body{
In the category of sets, pullbacks have the concrete description $X\times_Z Y := \{(x,y)\in X\times Y : f(x)=g(y)\}$ with projections sending pairs to their first and second coordinates, respectively. The  {pullback} of $\sigma$-structures $G \xrightarrow{g} K \xleftarrow{h} H$ is the $\sigma$-structure $G\times_K H$ whose underlying set of points is $\pts{G\times_K H}:=\{(x,y)\in\pts{G}\times\pts{H}\mid g(x)=h(y)\},$ and whose relations are defined, for each $i\in I$, by $((x_1,y_1),\ldots,(x_{a_i},y_{a_i}))\in R_i^P
    \iff
    (x_1,\ldots,x_{a_i})\in R_i^G
    \ \text{and}\
    (y_1,\ldots,y_{a_i})\in R_i^H.$
Whether a given category $\cat{C}$ admits pullbacks is a property that must be checked. The category $\cat{Grph}$ of graphs and graph homomorphisms admits pullbacks, and we already unpacked this definition in Paragraph~\ref{para: pullback Leticia}; see Figure~\ref{fig:pullback} for an example.
}

\begin{example}\label{ex:direct_product}
    The direct product $G\times H$ of a graph $G$ and a graph $H$ is equal to the pullback $G\times_o H$, where $o$ is the graph consisting of a single vertex with a loop.
\end{example}

\begin{definition}\label{def: matching families}
Let $\Ta = (T, \{G_a\},\{g_t^e\})$ be a tree co-decomposition. A \define{matching family} (cf.~\cite{Rosiak2022}) for $\Ta$ is a tuple of vertices $(x_t)_{t\in V(T)}$ such that for all edges $e = (s,t)$ of $T$, the images of $x_s$ and $x_t$ agree in the adhesion $G_e$; explicitly, $g_s^e(x_s) = g_t^e(x_t)$. The \define{limit} $\lim \Ta$ of $\Ta$ is the graph whose vertices are the set of matching families of $\Ta$ and where there is an edge from $(x_t)_{t\in V(T)}$ to $(y_t)_{t\in V(T)}$ whenever $(x_t,y_t)$ is an edge in $G_t$ for every $t\in V(T)$. We say that $\Ta$ \define{decomposes} a graph $G$ if $G \cong \lim \Ta$. The limit of a tree co-decomposition whose auxiliary tree is a single edge is exactly the pullback of the two bags over the adhesion.
\end{definition}

\begin{figure}
    \begin{minipage}{0.58\textwidth}
    \includegraphics[width=\textwidth, page=6]{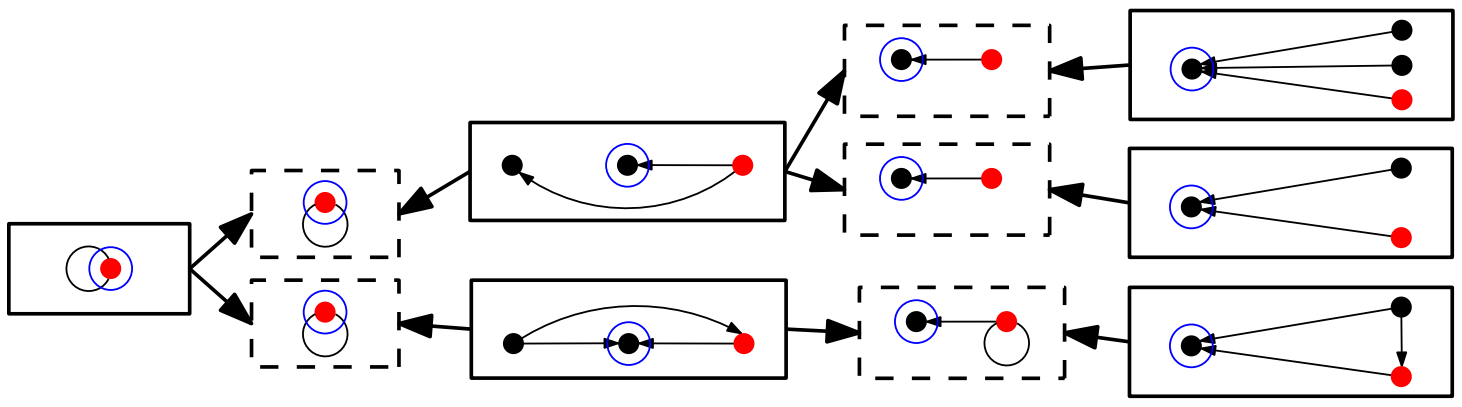}
    \end{minipage}%
    \hspace{\fill}
    \begin{minipage}{0.38\textwidth}
    \caption{Two matching families: the red vertex and blue/double circled vertex. Since they differ in some bags, these are two different vertices in the limit. Since there is one in every bag, there is an edge from the red to blue vertex in the limit.}\label{fig:extension-picture}
    \end{minipage}%
\end{figure}

\body{
Given a co-decomposition $\Ta = (T, \{G_a\},\{g_t^e\})$ and a graph $G$ so that $G\cong\lim \Ta$, we can think of $\Ta$ as providing a coordinate system to the vertices of $G$, indexed by the nodes of $T$. Indeed, the graph $G$ can be viewed explicitly as a subgraph of the categorical product of all of the bags of $\Ta$.
}

\begin{observation}
    Let $\Ta=(T, \{G_a\},\{g_t^e\})$ be a tree co-decomposition. Then $\lim\Ta\subseteq\prod_{t\in V(T)}G_t$, where $\prod$ is the iterated direct product.
\end{observation}

\body{
The construction of \cref{def: matching families} satisfies the universal property characterizing the categorical limit (c.f.~Riehl~\cite{Riehl2016} for background), as we show below \footnote{Co-decompositions are categorical duals of \textit{structured decompositions}~\cite{BumpusKocsis2023} (or, in some sense equivalently, of cellular sheaves~\cite{Curry2014, HansenGhrist2019, Robinson2017}) and hence~\cite[Section 2.8]{BumpusKocsis2023} the dual of familiar tree decompositions. To make this explicit, denote coproducts as $A + B$ and let $T$ be a tree. The \define{barycentric subdivision} of $T$ is a finite category $\smallint T$ having $V(T) + E(T)$ as objects and having exactly one morphism from $e \in E(T)$ to $x \in V(T)$ whenever $e$ is incident with $x$. From this, it follows definitionally that \textit{there is a bijective correspondence between {tree co-decompositions} and diagrams of the form $d \colon (\smallint T)^{\op} \to \cat{Grph}$ with the property that, for any morphism $f \colon e \to x$ in $\smallint T$, the morphism $d(f)$ is epic.}
}.
}

\begin{definition}
A \define{cone} over a tree co-decomposition $\Ta=(T, \{G_a\},\{g_t^e\})$ consists of a graph $H$ together with a family of homomorphisms $\phi_a \colon H \to G_a$ for $a\in V(T)\cup E(T)$ such that for every node $t\in V(T)$ and $t$-incident edge $e$, we have that $\phi_e=g_t^e\circ\phi_t$. This, and matching families, yields the following imeditately. 
\end{definition}

\begin{proposition}\label{prop:universal}
Let $\Ta=(T, \{G_a\},\{g_t^e\})$ be a tree co-decomposition. Then $\lim \Ta$ is universal among cones over $\Ta$. That is, for every cone $(H,\{\phi_a\})$ there exists a unique graph homomorphism $u \colon H \to \lim \Ta$ commuting with the projection maps.
\end{proposition}

\body{
Since the indexing graph is a tree, the limit of $\Ta$ can also be computed recursively by taking iterated pullbacks~\cite{Riehl2016}. This observation will be particularly useful in our algorithmic constructions. For completeness, we collect several structural results concerning co-decompositions in the appendices: we establish the relationship between tree co-width and tree-width (c.f.\ Appendix~\ref{appendix:tree-co-vs-tree}), we give a grammar characterization of path co-width (c.f.\ Appendix~\ref{appendix: grammar}) and we also show that every tree co-decomposition can be transformed into a \define{binary tree co-decomposition} without increasing its width (c.f. Appendix~\ref{appendix:binary}). Finally, in Appendix~\ref{appendix: tree shape important} we show that the restriction to tree-shaped diagrams is essential to obtaining a nontrivial width parameter: if arbitrary indexing shapes are allowed, every graph admits a co-decomposition of width at most three.
}

%% file: sections/path-hardness.tex
\section{NP-hardness for \pname{$(s,t)$-Path}}\label{sec: hardness}
\body{
    Throughout this paper, a co-decomposition is treated as the input object. That is, when we are given a tree co-decomposition $\mathcal T$, we are not also given its limit $\lim \mathcal T$. The limit is, of course, determined by $\mathcal T$, but it may be much larger than the data specifying $\mathcal T$: if the co-decomposition has $b$ bags of width at most $k$, then the vertex set of the limit may have size on the order of $k^b$. Thus, any algorithm which constructs the limit is not polynomial in the size of the co-decomposition, even if it is polynomial in the size of the limit.\label{para:bags-width-size}
}

\body{
    This distinction has computational consequences. For a given graph $G$, it is well-known that deciding whether there is an $s$-to-$t$ path is polynomial-time computable. However, when $G$ is given only \textit{implicitly}, as the limit of a co-decomposition, the vertices $s$ and $t$ are themselves matching families, and an $s$-to-$t$ path is a global sequence of compatible local moves across all bags. The purpose of this section is to show that, when the input consists only of a tree co-decomposition together with two vertex matching families, deciding whether these vertices are connected by a path in the limit is \textsf{NP}-hard. In contrast, later in the paper we show that subgraph isomorphism, and more generally first-order model checking, is fixed-parameter tractable when parameterized by tree co-width and linear in the number of bags. Since $(s,t)$-path is not FO-expressible, there is a sense in which this section demonstrates the algorithmic tightness of our algorithmic meta-theorems. 
}

\body{
  Towards our hardness result, we reduce from \pname{Exact Cover By 3-Sets} (\pname{X3C}) \cite{garey2002computers}, where, given a ground set $U$ and a family of triples $\mathcal{S}\subseteq {U\choose 3}$ covering $U$, one must decide if
  there exists some $\mathcal{S}' \subseteq \mathcal{S}$ that partitions $U$. For
  simplicity, let $(U, \mathcal{S})$ be our input instance with $U = \{u_1, \dots,
  u_n\}$, and $\mathcal{S} = \{S_1, \dots, S_m\}$.
}

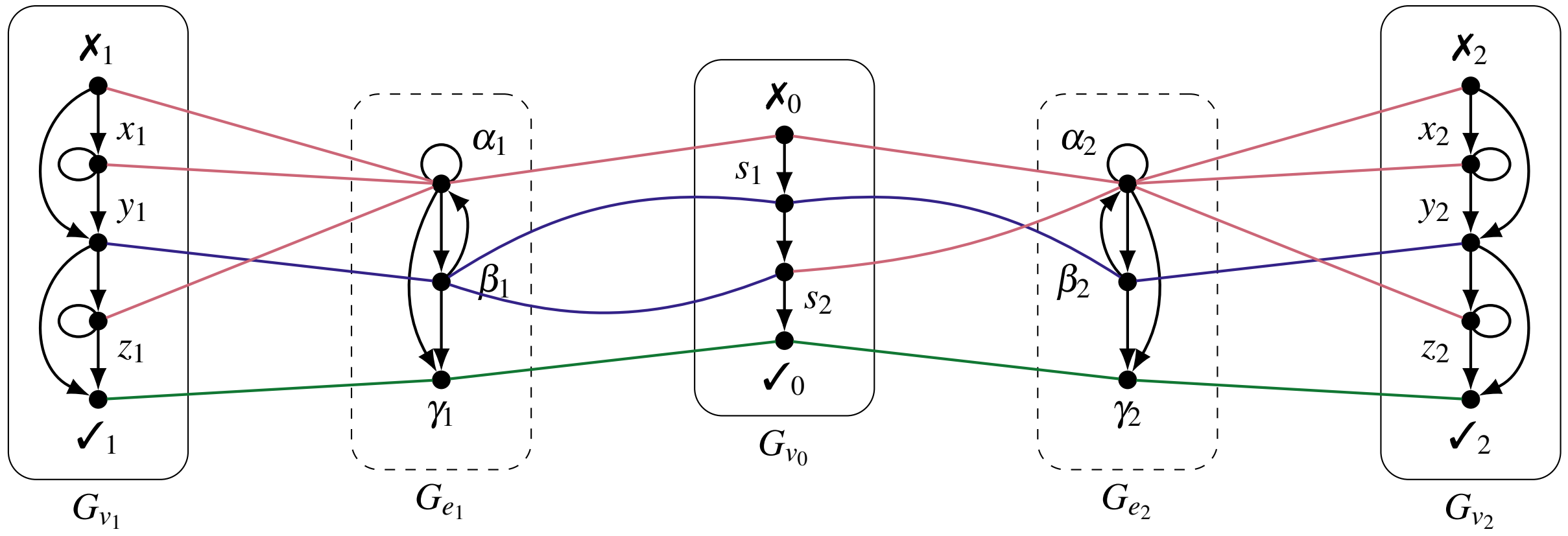
\begin{figure}
  \centering
  \newcommand{\cyscale}{0.7}
  \begin{tikzpicture}[>=Latex]
    \tikzset{
      VertexStyle/.style={
        shape=circle,
        draw=black,
        fill=black,
        minimum size=5pt,
        inner sep=0pt
      },
      graphbox/.style={
        draw,
        rounded corners=8pt,
        inner sep=8pt
      },
      mapbox/.style={
        draw,
        dashed,
        rounded corners=6pt,
        inner sep=5pt,
        align=center
      },
      perfect loop/.style 2 args={to path={
          arc[start angle=#1-180, end angle=#1+180, radius=#2] \tikztonodes
      }},
      functional graph/v label pos/.store in=\fgVLabelPos,   functional graph/v label pos=45,
      functional graph/loop dir/.store in=\fgLoopDir,       functional graph/loop dir=180,
      functional graph/id/.store in=\fgId,                 functional graph/id=1,
      functional graph/bend dir/.store in=\fgBendDir,       functional graph/bend dir=right,
      functional graph/.pic={
        \begin{scope}[local bounding box=-bbox,yscale=0.8]

          \Vertex[x=0, y=0,  Math, Lpos=90,           L={\protect\xmark_{\fgId}}]{mark\fgId}
          \Vertex[x=0, y=-1, Math, Lpos=\fgVLabelPos, L={x_{\fgId}}]{x\fgId}
          \Vertex[x=0, y=-2, Math, Lpos=\fgVLabelPos, L={y_{\fgId}}]{y\fgId}
          \Vertex[x=0, y=-3, Math, Lpos=-\fgVLabelPos, L={z_{\fgId}}]{z\fgId}
          \Vertex[x=0, y=-4, Math, Lpos=-90,          L={\protect\mycmark_{\fgId}}]{c\fgId}

          \Edge[style={->}](mark\fgId)(x\fgId)
          \Edge[style={->, bend \fgBendDir=60}](mark\fgId)(y\fgId)
          \Edge[style={perfect loop={\fgLoopDir}{0.2cm}}](x\fgId)(x\fgId)
          \Edge[style={->}](x\fgId)(y\fgId)
          \Edge[style={->}](y\fgId)(z\fgId)
          \Edge[style={->, bend \fgBendDir=60}](y\fgId)(c\fgId)
          \Edge[style={perfect loop={\fgLoopDir}{0.2cm}}](z\fgId)(z\fgId)
          \Edge[style={->}](z\fgId)(c\fgId)

          \node (extra10) at ($(mark\fgId)+(0.5,0.5)$) {};
          \node (extra11) at ($(c\fgId)-(0.5,0.5)$) {};

          \node[
            graphbox,
            fit=(mark\fgId)(x\fgId)(y\fgId)(z\fgId)(c\fgId)(extra10)(extra11),
            label=south:$G_{v_{\fgId}}$
          ] (Gv\fgId) {};

        \end{scope}
      },
      adhesion graph/alpha label pos/.store in=\agAlphaLabelPos,
      adhesion graph/alpha label pos=135,
      adhesion graph/beta label pos/.store in=\agBetaLabelPos,
      adhesion graph/beta label pos=0,
      adhesion graph/loop dir/.store in=\agLoopDir,       adhesion graph/loop dir=90,
      adhesion graph/id/.store in=\agId,                 adhesion graph/id=1,
      adhesion graph/bend dir/.store in=\agBendDir,       adhesion graph/bend dir=right,
      adhesion graph/.pic={
        \begin{scope}[local bounding box=-bbox]

          \Vertex[x=0, y=0,  Math, Ldist=4pt,Lpos=\agAlphaLabelPos, L={\alpha_{\agId}}]{alpha\agId}
          \Vertex[x=0, y=-1, Math, Ldist=4pt,Lpos=\agBetaLabelPos, L={\beta_{\agId}}]{beta\agId}
          \Vertex[x=0, y=-2, Math, Lpos=-90, L={\gamma_{\agId}}]{gamma\agId}

          \Edge[style={->}](alpha\agId)(beta\agId)
          \Edge[style={->}](beta\agId)(gamma\agId)
          \Edge[style={->, bend \agBendDir=30}](alpha\agId)(gamma\agId)

          \Edge[style={perfect loop={\agLoopDir}{0.2cm}}](alpha\agId)(alpha\agId)

          \Edge[style={->, bend \agBendDir=45}](beta\agId)(alpha\agId)

          \node (extra10) at ($(alpha\agId)+(0.5,0.5)$) {};
          \node (extra11) at ($(gamma\agId)-(0.5,0.5)$) {};

          \node[
            graphbox, dashed,
            fit=(alpha\agId)(beta\agId)(gamma\agId)(extra10)(extra11),
            label=south:$G_{e_{\agId}}$
          ] (Ge\agId) {};

        \end{scope}
      }
    }
    \SetVertexLabelOut
    
    \begin{scope}[yshift=-0.5cm,yscale=\cyscale]
      \Vertex[x=0,y=0,Math,Lpos=90,L={\protect\xmark_0}]{mark0}
      \Vertex[x=0,y=-1,Math,Lpos=135,L={s_1}]{s1}
      \Vertex[x=0,y=-2,Math,Lpos=-45,L={s_2}]{s2}
      \Vertex[x=0,y=-3,Math,Lpos=-90,L={\protect\mycmark_0}]{c0}

      \Edge[style={->}](mark0)(s1)
      \Edge[style={->}](s1)(s2)
      \Edge[style={->}](s2)(c0)
      \node (extra00) at ($(mark0)+(0.5,0.5)$) {};
      \node (extra01) at ($(c0)-(0.5,0.5)$) {};

      \node[
        graphbox,
        fit=(mark0)(s1)(s2)(c0)(extra00)(extra01),
        label=south:$G_{v_0}$
      ] (Gv0) {};
    \end{scope}

    \begin{scope}[xshift=-7cm,yscale=\cyscale]
      \pic at (0,0) [functional graph/id=1] {functional graph};
    \end{scope}

    \begin{scope}[xshift=-3.5cm,yshift=-1cm,yscale=\cyscale]
      \pic at (0,0) [adhesion graph/id=1, adhesion graph/alpha label pos=45] {adhesion graph};
    \end{scope}

    \begin{scope}[xshift=3.5cm,yshift=-1cm,yscale=\cyscale]
      \pic at (0,0) [adhesion graph/id=2, adhesion graph/alpha label pos=135, adhesion graph/beta label pos=180,adhesion graph/bend dir=left] {adhesion graph};
    \end{scope}

    \begin{scope}[xshift=7cm,yscale=\cyscale]
      \pic at (0,0) [functional graph/id=2, functional graph/bend dir=left, functional
      graph/loop dir=0, functional graph/v label pos=135] {functional graph};
    \end{scope}

    \Edge[style={color=goodred}](mark0)(alpha1)
    \Edge[style={color=goodred}](mark0)(alpha2)
    
    \Edge[style={color=goodblue,bend right=20}](s1)(beta1)
    \Edge[style={color=goodblue,bend left=20}](s1)(beta2)
    
    \Edge[style={color=goodblue,bend left=20}](s2)(beta1)
    \Edge[style={color=goodred,bend right=10}](s2)(alpha2)
    
    \Edges[style={color=goodgreen}](c1,gamma1,c0,gamma2,c2)

    \foreach \i in {1,2} {
        \Edge[style={color=goodred}](mark\i)(alpha\i)
        \Edge[style={color=goodred}](x\i)(alpha\i)

        \Edge[style={color=goodblue}](y\i)(beta\i)

        \Edge[style={color=goodred}](z\i)(alpha\i)
    }

  \end{tikzpicture}
  \vspace{-1em}
  \caption{Part of the constructed co-decomposition for a \pname{X3C} instance with $\{1,2\} \subset S_1$ and $\{1,3\} \subset S_2$, with some arcs in $G_{v_0}$ omitted for clarity. The coloured edges represent the adhesion mappings: red edges are incident to vertices that, when swapped into the matching family, indicate that the corresponding element is not covered by the set vertex; blue edges are incident to vertices that correspond to element-covering swaps; the green edges are incident to the vertices of the bags that, together, form the target matching family.}
  \label{fig:path_hardness}
\end{figure}

\body{
  Construct a tree co-decomposition $\mathcal{T} = (T, \{G_a\}, \{g_t^{e}\})$ from
  $(U, \mathcal{S})$ as follows (c.f.  \cref{fig:path_hardness}): (i) $T$ is a star with $V(T) = \{{v_0}, {v_1}, \dots,
  {v_n}\}$ and $E(T) = \{{e_1}, \dots, {e_n}\}$, with ${v_0}$ being the
  central node, and ${e_i}$ the edge between ${v_0}$ and ${v_i}$; (ii) $G_{v_0}$
  is a tournament on the vertex set $\{\xmark_0, s_1, \dots, s_m,
  \mycmark_0\}$  obtained as the transitive closure of $\xmark_0\to s_1\to s_2\to\ldots\to s_n\to\mycmark_0$; (iii) the other bags $G_{v_i}$ are isomorphic to one another, and have $V(G_{v_i}) = \{\xmark_i,
  x_i, y_i, z_i, \mycmark_i\}$ and $E(G_{v_i}) = \{(\xmark_i, x_i), (\xmark_i, y_i),
  (x_i, x_i), (y_i, z_i), (y_i, \mycmark_i), (z_i, z_i), (z_i, \mycmark_i)\}$; (iv) each
  adhesion $G_{e_i}$ is the graph with vertex set $\{\alpha_i, \beta_i,
  \gamma_i\}$ and the arcs $\{(\alpha_i,\alpha_i),(\alpha_i,\beta_i),(\alpha_i,\gamma_i),(\beta_i,\alpha_i)\}$; (v) for each $u_i \in U$ we define
  $g_{v_i}^{{e_i}}$ by $(g_{v_i}^{{e_i}})^{-1}(\alpha_i) = \{\xmark_i, x_i, z_i\}$,
  $(g_{v_i}^{{e_i}})^{-1}(\beta_i) = \{y_i\}$, and $(g_{v_i}^{{e_i}})^{-1}(\gamma_i) =
  \{\mycmark_i\}$; (vi) for each $S_j \in \mathcal{S}$ and each $u_i \in U$, we set:\vspace{-1.5em}

  \[
    g_{v_0}^{e_i}(s_j)
    =
    \begin{cases}
      \alpha_i & \text{if $u_i \notin S_j$},\\
      \beta_i & \text{if $u_i \in S_j$};
    \end{cases}\vspace{-.3em}
  \]

  (vii) set $g_{v_0}^{e_i}(\xmark_0) = \alpha_i$ and $g_{v_0}^{e_i}(\mycmark_0) =
  \{\gamma_i\}$; (viii) finally, we write $\vec{\xmark} := (\xmark_0, \dots, \xmark_n)$ and $\vec{\mycmark} := (\mycmark_0,
  \dots, \mycmark_n)$, which will be the start and destination for our path.
}

\body{
  Intuitively, we see this as a reconfiguration problem in the co-decomposition,
  where selecting a set of $\mathcal{S}$ equates with moving one token in
  each bag of $\mathcal{T}$: upon picking $S_j$, the token on $G_{v_0}$ goes from its
  previous vertex to $s_j$, which forces a token movement in $G_{v_i}$ to either $x_i$
  (if $u_i \notin S_j$ and $u_i$ has not been covered), $z_i$ (if $u_i$ has already been
  covered by another $S_r$), or $y_i$ (if $u_i$ is being covered by $S_j$). The adhesions
  and homomorphisms are constructed so as to impede the re-covering of any $u_i$. The following completes the proof of Theorem~\ref{thm: s-t path hardness}: given only a co-decomposition and two of its matching families, \pname{$(s,t)$-Path} is $\np$-hard.
}

\begin{proposition}\label{prop:path_equivalence}\vspace{-.05cm}
  There is a subset of $\mathcal{S}$ partitioning $U$ if and only if there is a
  path from $\vec{\xmark}$ to $\vec{\mycmark}$ in~$L$.
\end{proposition}

\begin{proof}
  Let $\mathcal{S}'$ be a solution to $(U, \mathcal{S})$ and $L = \lim \Ta$. Write $\mathcal{S}'=\{S_{j_1},S_{j_2},\ldots,S_{j_t}\}$ with $j_1<j_2<\ldots<j_t$. For a
  given matching family $\vec{\ell}\in V(L)$ and index $i \in \{0, \dots, n\}$, we write $\vec{\ell}[i]$ for the coordinate of $\vec{\ell}$ in bag $G_{v_i}$. We will define a sequence of matching families $\vec{\ell}_1,\ldots,\vec{\ell}_t$ so that $\xmark\vec{\ell}_1\ldots\vec{\ell}_t\mycmark$ is our required path in $L$.
  For the first step, write  $S_{j_1}=\{u_a, u_b, u_c\}$; we define $\vec{\ell}_1$ so that $\vec{\ell_1}[0] = s_j$, for $i \in
  \{a,b,c\}$ we have $\vec{\ell_1}[i] = y_i$, and for every $i \notin \{a,b,c\}$, we have
  $\vec{\ell_i} = x_i$. Observe that every $\xmark_i$ has an arc to $\ell_1[i]$.
  After the first step, we proceed inductively: $\vec{\ell}_{k+1}$ is
  obtained from $\vec{\ell}_{k}$ by: (i) replacing $\vec{\ell}_{k}[0] = s_{j_k}$ with
  $s_{j_{k+1}}$ (since $j_k < j_{k+1}$, we have $s_{j_k}s_{j_{k+1}} \in E(G_{v_0})$); (ii) for each
  $u_i \in S_{j_{k+1}}$, replace $\vec{\ell}_{k}[i] = x_i$ with $y_i$; and (iii) for each
  $u_i \notin S_{j_{k+1}}$:\vspace{-1.5em}

  \[
    \vec{\ell}_{k+1}[i] =
    \begin{cases}
      z_i & \text{if $u_i$ was covered by $S_{j_k}$},\\
      \vec{\ell}_{k}[i] & \text{otherwise.}
    \end{cases}\vspace{-.3em}
  \]

  To verify $(\vec{\ell}_{k}, \vec{\ell}_{k+1}) \in E(L)$, note that, during step (ii),
  $u_i$ has not been covered before, and so it holds that $\vec{\ell}_{k}[i] =
  \vec{\ell}_{k-1}[i] = \dots = \vec{\ell}_{1}[i] = x_i$, which has an arc to $y_i$;
  during step (iii), if $u_i$ was covered by $S_{j_k}$, then $\vec{\ell}_{k}[i] = y_i$, which
  has an arc to $z_i$, otherwise $\vec{\ell}_{k}[i] \notin \{y_i, \mycmark_i\}$ and
  $\vec{\ell}_{k}$ has a loop.
  Finally, after processing $S_{j_t}$, every $\vec{\ell}_{t}[i]$ is
  either $s_{j_t}$ if $i = 0$, $y_i$ or $z_i$, all of which have arcs to their
  corresponding to $\mycmark_i$, and so $\vec{\ell}_{t}$ has an arc to $\vec{\mycmark}$
  in $L$. Consequently, $(\vec{\xmark}, \vec{\ell}_{1}, \dots, \vec{\ell}_{t},
  \vec{\mycmark})$ is a path in $L=\lim\Ta$. For the converse, let $P = (\vec{\xmark}, \vec{\ell}_{1}, \dots, \vec{\ell}_{t},
  \vec{\mycmark})$ be a sequence of matching families of $\mathcal{T}$ corresponding to a
  path in $L$. For each $k\in [t]$, let $j_k$ be the index so that $\vec{\ell}_k[0]=s_{j_k}$. By construction, each $s_{j_k}$ is distinct as
  $G_{v_0}$ is acyclic. We claim that $\mathcal{S}' := \{S_{j_k} \colon k\in[t] \}$ is a partition of $U$.
  First, suppose for a contradiction that there is some $u_i$ that is not covered by $\mathcal{S}'$. By construction, for each $k\in[t]$ we have $g_{v_0}^{e_i}(s_{j_k}) = \alpha_i$, so it holds that
  $\vec{\ell}_{1}[i] = x_i$ as the only arc from $\xmark_i$ to
  $(g_{v_i}^{e_i})^{-1}(\alpha_i)$ is to $x_i$; moreover, this implies that
  $\vec{\ell}_{t}[i] = \dots = \vec{\ell}_{1}[i] = x_i$, but $x_i$ has no arc to $\mycmark_i$,
  so $(\vec{\ell}_{t}, \vec{\mycmark}) \notin E(L)$, a contradiction. To see that $u_i$ is covered only once, let $k$ be the first index for which $u_i \in
  S_{j_k}$; by a similar argument to the previous one, $\vec{\ell}_{k}[i] = y_i$ and for
  every $r \geq k$ it must be that $\vec{\ell}_{r}[i] \notin \{x_i, y_i\}$ as no other
  vertex of $G_{v_i}$ reachable from $y_i$ has an arc to any of $\{x_i, y_i\}$.
  Consequently, $\vec{\ell}_{r}[i] = z_i$, no other vertex of $\vec{\ell}_{r}[0] \notin
  (g_{v_0}^{e_i})^{-1}(\beta_i)$, and so $S_{j_r}$ cannot contain $u_i$. This completes our proof.
\end{proof}

%% file: sections/existsFO.tex
\section{Warm-up: Model Checking for $\exists$-FO}\label{sec: exists FO}

\body{
    This whole section is devoted to proving the following Theorem~\ref{thm:exists-fo}. We will work with \textit{grounding} of logical formulas, that is the act of replacing quantifiers with logical connectives ranging over the entire domain. For example we will turn a formula $\exists x \exists y \phi(x,y)$ into the formula $\bigvee_{x' \in G} \bigvee_{y' \in G} \phi(x', y')$. 
}

\begin{theorem}\label{thm:exists-fo}
Let $\Ta$ be a rooted binary tree co-decomposition of $\sigma$-structures with $b$ bags, each containing at most $w$ points, and let $\phi$ be an existential first-order sentence with $\ell$ variables. Letting $q=\sum_{R\in\sigma}\ell^{\operatorname{ar}(R)}$, then whether $\lim\Ta\models\phi$ can be decided, without constructing $\lim\Ta$, in time\vspace{-.5em}
\[
    \bigO{b(\ell+q)w^{2\ell}2^{2q} + w^\ell 2^q |\phi|}.\vspace{-.8em}
\]
In particular, for graphs, this gives a running time of $b\cdot 2^{\bigO{\ell^2+\ell\log(w+1)}} \operatorname{poly}(|\phi|).$ 
\end{theorem}

\begin{definition}
    Fix an ambient $\ell$-variable formula $\phi$. An \define{assignment} on a structure $G$ is a function $[\ell]\rightarrow \pts{G}$.
\end{definition}

\body{
    Since grounding the formula over the entire limit might take exponential time, we will use dynamic programming to reason about the global grounding of the formula in terms of these local assignments. To this end, we need to record not only the assignment, but also whether we already local evidence that the overall formula is false under the assignment. A canonical way of doing this is by inducing all possible relations on the points chosen by the assignment; we call these \textit{hypotheses}.
}

\begin{definition}
    Given an assignment $\beta \colon [\ell] \to \pts{G}$, denote by $\beta[G] \colon N \to G$ the unique\footnote{The structure $N$ is unique and so is the map.} homomorphism of structures -- which we call the \define{hypothesis induced on $G$ by $\beta$} -- with the following two properties: (1) the underlying function $\pts{\beta[G]} \colon \pts{N} \to \pts{G}$ is exactly the assignment $\beta$, (2) for each relation $R$ in the signature and each tuple $\vec{x}\in[\ell]^{\operatorname{ar}(R)}$, we have $R^N(\vec{x})$ if and only if $R^G(\beta(\vec{x}))$. 
\end{definition}

\body{
    As noted in Paragraph~\ref{para: bayesian}, dynamic programming over co-decompositions has a certain Bayesian feel: until one reaches the root of the co-decomposition tree, one can disprove atomic statements about the entire limit, but one can only conjecture their proofs. Hence, on their own, induced hypotheses are not enough: as we traverse the co-decomposition, we might want to remember that some relation that holds locally actually didn't hold further down the previously-explored tree. This motivates the following weaker notion.
}

\begin{figure}
    \begin{minipage}{0.68\textwidth}
    \includegraphics[width=0.7\linewidth, page=15]{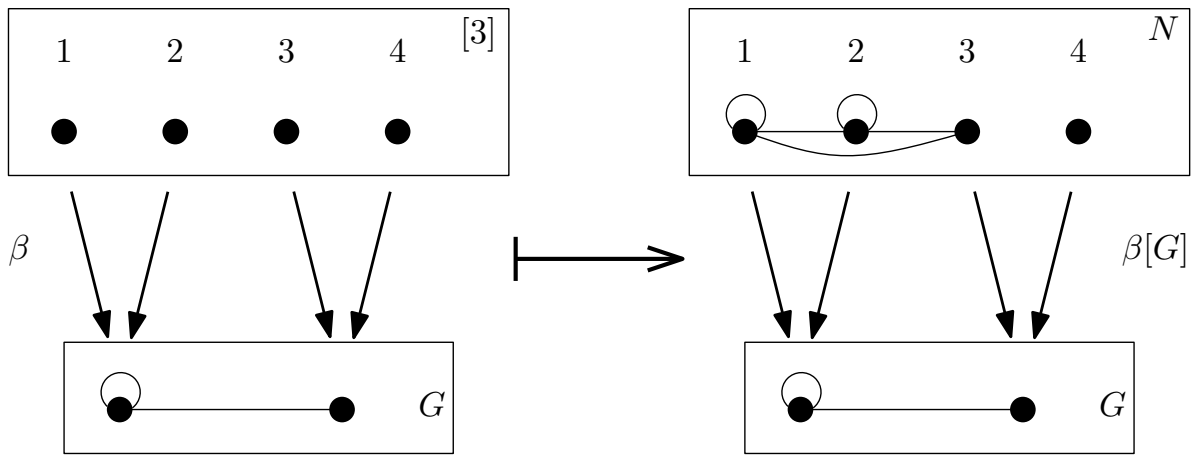}
    \vspace{1em}
    \end{minipage}%
    \hspace{\fill}
    \begin{minipage}{0.28\textwidth}
    \caption{(Left) An assignment $\beta$ sits above a graph $G$. (Right) This induces an induced hypothesis $\beta[G]$ sitting over $G$.}
    \label{fig:induced_env}
    \end{minipage}%
\end{figure}

\begin{definition}\label{def: hypothesis}
    A \define{hypothesis} in a structure $G$ is a homomorphism $N\xrightarrow{\nu } G$ where $\underline{N} = [\ell]$. We call $N$ the \define{head} of $\nu$ and denote it by $\text{head}(\nu)$. The \define{underlying assignment for $\nu$} is the function $\underline{\nu} \colon [\ell] \rightarrow \pts{G}$ on the sets of points. We think of a hypothesis as picking out $\ell$ points from $G$ and recording some of their relations to one another (but not necessarily all of the relational data). An induced hypothesis is a genuine induced substructure in the case where the map $\beta$ is injective. See Figure~\ref{fig:induced_env} for an example. 
\end{definition}

\body{
    We define morphisms of hypotheses below. These are important for two reasons. The first is that thinking in terms of morphisms will allow us to show that induced hypotheses satisfy a certain universal property (c.f. Proposition~\ref{prop: Universal property of induced hypotheses}) which we will use later on in proofs. The second is that morphisms of hypotheses will allow us to speak of pullbacks of hypotheses, which is the crucial ingredient towards establishing our dynamic programming algorithms both for $\exists$-FO and full FO model checking (Theorems~\ref{thm:exists-fo} and~\ref{thm: FO model checking on codecompositions}).
}


\vspace{-0.8em}
\begin{minipage}[t]{0.50\textwidth}
\begin{definition}\label{def: morphism of hypotheses}
    Given hypotheses $N\xrightarrow{\nu} G$ and $M\xrightarrow{\mu} H$ in structures $G$ and $H$, respectively, a \define{morphism of hypotheses} is a homomorphism $G\xrightarrow{f} H$ such that there exists a map $N\xrightarrow{q} M$ whose underlying function is $\mathsf{id}_{[\ell]}$ and which makes the diagram on the right commute.
\end{definition}
\end{minipage}
\begin{minipage}[t]{0.40\textwidth}{\small
\[\hspace{1em}\begin{tikzcd}
	{[\ell]} && {[\ell]} \\
	N && M \\
	G && H
	\arrow["{\mathsf{id}}"{description}, from=1-1, to=1-3]
	\arrow[hook, from=1-1, to=2-1]
	\arrow[hook, from=1-3, to=2-3]
	\arrow["\exists", dashed, from=2-1, to=2-3]
	\arrow["\nu"{description}, from=2-1, to=3-1]
	\arrow["\mu"{description}, from=2-3, to=3-3]
	\arrow["f"{description}, from=3-1, to=3-3]
\end{tikzcd}\]
}
\end{minipage}
\vspace{-0.3em}

\begin{lemma}\label{lemma: limits of hypotheses}
    Fix $\ell\in\mathbb{N}$. The \define{category $\mathbf{Hyp}^\sigma_{\ell}$ of $[\ell]$-hypotheses of $\sigma$-structures} is defined as follows. Its objects are hypotheses (that is, homomorphisms $\nu:N\rightarrow G$ of $\sigma$-structures such that $\pts{N}=[\ell]$) and its morphisms are given as in \cref{def: morphism of hypotheses}. Identities and composition are defined componentwise. Furthermore, the category of hypotheses and hypothesis morphisms admits pullbacks.
\end{lemma}\vspace{-.1cm}
\begin{proof}
    Verifying that $\mathbf{Hyp}^\sigma_{\ell}$ is a category is routine; we focus on proving that it admits pullbacks. The pullback of a span $(\nu:N\to G)\xrightarrow{g}(\eta:L\to K)\xleftarrow{h}(\mu:M\to H)$ of hypotheses over a span $G\xrightarrow{g}K\xleftarrow{h}H$ of structures is the map $\nu\times_\eta\mu:N\times_L M\rightarrow G\times_K H$. Viewing $\text{head}(\nu\times_\eta\mu)=N\times_L M$ again as a structure with point set $[\ell]$, we get that the pullback is again a hypothesis. In particular, a hypothesis in $G\times_K H$.
\end{proof}

\body{
    Recall that the limit of a tree co-decomposition is just an iterated pullback. Our next result shows the same is true of induced hypotheses: given an assignment to the global structure, the global\footnote{Here, and throughout, we mean ``global'' in the sense of being an assignment into the limit of the co-decompostiion.} hypothesis it induces can be recovered by taking the limit of the hypotheses induced by its projections to the bags.
}

\begin{lemma}\label{lem:lim_env_is_env_lim}
    Let $\Ta=(T,\{G_a\},\{g_t^e\})$ be a tree co-decomposition and let $\beta:[\ell]\rightarrow\pts{\lim\Ta}$ be a global assignment. Write $\beta_a$ for the projection of $\beta$ to $G_a$. Then the tuple $\Ta'=(T,\{\beta_a[G_a]\},\{g_t^e\})$ is a co-decomposition of hypotheses. Furthermore, the limit of $\Ta'$ is the hypothesis induced on the limit of $\Ta$ by $\beta$. That is, $\lim\Ta'=\beta[\lim\Ta]$. 
\end{lemma}\vspace{-.2cm}

\body{
    The universal property of pullbacks, on its own, is almost enough to prove Lemma~\ref{lem:lim_env_is_env_lim}. Indeed we require only one further universal property enojoyed by induced hypotheses. A similar universal property (c.f. Lemma~\ref{lemma: canonical comparison from projected induced to induced projected}) will be leveraged later when dealing with the case of full FO model checking. 
}

\body{
    The point of the preceding lemma is that a hypothesis on a limit is just a matching family of hypotheses in a suitably-defined co-decomposition. To make use of this, we need to project hypotheses: given a homomorphism $G\xrightarrow{\pi} H$ and a hypothesis $N\xrightarrow{\nu} G$, we call the composite $ N\xrightarrow{\pi\circ\nu} H$ the \define{projection of $\nu$ along $\pi$}. Usually $\pi$ will be the canonical projection from a (partial) limit to a bag: the resulting projection allows us to store, inside that bag, the relational information accumulated over the entire partial limit. 
}

\begin{observation}
    Notice that the projection of a hypothesis $\nu \colon N\rightarrow G$ along a map $\pi \colon G\rightarrow H$ leaves the head $N$, and hence the relations it records, unchanged. Furthermore: (1) the size of $N$ depends on the number of variables in the formula; and (2) the resulting hypothesis need not be induced (capturing the intuitive idea that, while proceeding up the tree, some relation visible locally may already have been disproved elsewhere).
\end{observation}

\begin{definition}\label{para: def: supported hypoethis}
    Now let $\Ta=(T,\{G_a\},\{g_t^e\})$ be a rooted co-decomposition (the tree $T$ has a root). For each node $i\in V(T)$, write $G[i]$ for the \define{partial limit} of $\Ta$ at $i$. That is, the limit of the tree co-decomposition obtained from $\Ta$ by restricting to the subtree rooted at $i$.
    Now fix $i$ and let $\pi_i\colon G[i]\rightarrow G_i$ be the canonical projection from the partial limit $G[i]$ to its root bag $G_i$. We say that a hypothesis $\nu_i\colon N_i\rightarrow G_i$ is \define{supported at $i$} if $\nu_i=\pi_i\circ\beta[G[i]]$ for some assignment $\beta\colon [\ell]\rightarrow\pts{G[i]}$. Given an assignment $\beta_i\colon [\ell]\rightarrow\pts{G_i}$, the \define{supported hypotheses for $\beta_i$} are those whose underlying assignment is $\beta_i$. In our analogy to Bayesian thinking: a hypothesis is supported if it has not been disproved by anything in the subtree rooted at $i$. From a dynamic programming point of view, all we need to check is the supported hypotheses at the root. \vspace{-0.1cm}
\end{definition}

\begin{lemma}\label{lemma: support at root is enough}
    Let $\phi$ be an $\ell$-variable $\exists$-FO formula, and $\Ta=(T,\{G_a\},\{g_t^e\})$ be a co-decomposition rooted at $r$. Then $\lim\Ta\models\phi$ if and only if $\operatorname{head}(\nu)\models\phi$ for some hypothesis $\nu$ supported at $r$. 
\end{lemma}
\begin{proof}
    The formula $\phi$ is satisfied in $\lim\Ta$ if and only if there is a satisfying assignment $\beta\colon [\ell]\rightarrow\pts{\lim\Ta}$. By definition, the induced hypothesis $\beta[\lim\Ta]$ records exactly the relations satisfied by the points selected by $\beta$, so $\beta$ is satisfying if and only if $\operatorname{head}(\beta[\lim\Ta])\models\phi$. Since $G[r]=\lim\Ta$, the hypotheses supported at $r$ are precisely the projections of these induced hypotheses to $G_r$. Projection leaves the head unchanged.
\end{proof}

\body{
    Lemma~\ref{lemma: support at root is enough} suggests a dynamic program. We proceed up the co-decomposition tree, starting from the leaves. To each bag $G_t$ we associate the set of supported hypotheses, and to each adhesion $G_e$ the set of hypotheses induced by assignments to $G_e$. These sets form a co-decomposition of sets, where the connecting maps send each hypothesis to the hypothesis induced by its underlying assignment on the corresponding adhesion. Matching families in this auxiliary co-decomposition are families of hypotheses whose underlying assignments agree on the adhesions. We proceed as follows. At the leaves, we simply enumerate the induced hypotheses. At internal nodes, we form matching families from the supported hypotheses of the children and induced hypotheses on the current bag. Then we combine them by taking \textit{pullbacks of hypotheses}. Finally, we project the results back to the current bag. By Lemma~\ref{lem:lim_env_is_env_lim}, we recover exactly the supported hypotheses. 
}\label{para: algorithm roadmap}

\body{
    As to the bound on the running time, note that hypotheses' heads always have precisely $\ell$ points, regardless of the size of the limit. Hence, at each node we only ever deal with at most $w^\ell 2^q$ hypotheses, where $q=\sum_{R\in\sigma}\ell^{\operatorname{ar}(R)}$. Since the algorithm considers at most pairs of such hypotheses at each node, the claimed running time follows. Finally, by Lemma~\ref{lemma: support at root is enough}, it suffices to inspect the heads of the supported hypotheses at the root to decide whether $\lim\Ta\models\phi$. The formal recurrences and their proofs are given in Section~\ref{sec: exists FO algorithms nitty-gritty}.
}

\body{
    We reiterate that the structure outlined in Paragraph~\ref{para: algorithm roadmap} is a roadmap to the case of full FO. Although the objects will be much more complicated (we will have \textit{strategies}, rather than assignments and \textit{hypotheticals}, rather than hypotheses), the algorithm is basically the same: (1) construct a co-decomposition of auxiliary data, (2) find matching families and then (3) combine these using pullbacks defined on this auxiliary data. 
}

\newpage
\subsection{Proof of Lemma~\ref{lem:lim_env_is_env_lim}}\label{sec: proof of lem:lim_env_is_env_lim}

\body{
    The proof of Lemma~\ref{lem:lim_env_is_env_lim} is straightforward, once one establishes that induced hypotheses satisfy a certain universal property (Proposition~\ref{prop: Universal property of induced hypotheses}). 
}
\begin{proposition}\label{prop: Universal property of induced hypotheses} 
    Given an assignment $\beta \colon [\ell] \to \pts{G}$, an induced hypothesis $\beta[G]$ satisfies: for any other hypothesis of the form $\gamma \colon N\to G$  with $\underline{\gamma} = \beta$, there is a unique map of hypotheses $\gamma \to \beta[G]$. 
\end{proposition}
\begin{proof}
     The unique map on the ground sets $q \colon N \to \text{head}(\beta[G])$ is fixed; we need to show it is a morphism of relational structures. Take a relation symbol $R$ and a tuple $\vec{x}$ in $N$ of the appropriate arity. By definition of morphisms of relational structures, we have that, if $R(\vec{x})$, then $R(\gamma(\vec{x}))$. But then, by definition of induced hypothesis, we have that $R(q(\vec{x}))$. Since $R$ was arbitrary $R$, we get $q$ is a morphism of structures. 
\end{proof}

\body{
    The following result states that morphisms of structures yield morphisms of induced hypotheses. We will use this result in Section~\ref{sec: exists FO} to establish Lemma~\ref{lem:lim_env_is_env_lim}, the main result in this section. In particular, we use it to prove that the tuple $\Ta'$ built from induced hypotheses is indeed a co-decomposition of hypotheses.
}

\begin{lemma}\label{lem:helper_1}
    Let $G$ and $H$ be structures and $f:G\rightarrow H$ be a homomorphism. Let $\beta:[\ell]\rightarrow \pts{G}$ be an assignment. Then $f$ is a morphism of hypotheses from $\beta[G]$ to $(f\circ\beta)[G]$.
\end{lemma}
\begin{minipage}[t]{0.45\textwidth}
    \begin{proof}
    Consider the diagram on the right. By Proposition~\ref{prop: Universal property of induced hypotheses}, and since $(f  \circ \beta)[H]$ is an induced hypothesis, there is a unique map $\text{head}(\beta[G]) \to \text{head}((f \circ \beta)[H])$ making the diagram to the right commute. 
    \end{proof}
\end{minipage}
\begin{minipage}[t]{0.45\textwidth}
\vspace{-2em}
    \[\hspace{2em}\begin{tikzcd}
	& {\text{head}(\beta[G])} & {\text{head}((f \circ \beta)[H])} \\
	\\
	{[\ell]} & G & H
	\arrow["{\beta[G]}", from=1-2, to=3-2]
	\arrow["{(f\circ \beta)[H]}", from=1-3, to=3-3]
	\arrow[from=3-1, to=1-2]
	\arrow["\beta"', from=3-1, to=3-2]
	\arrow["f"', from=3-2, to=3-3]
    \end{tikzcd}\]
    \end{minipage}

\begin{proof}[Proof of Lemma~\ref{lem:lim_env_is_env_lim}]
    For $\Ta'$ to even be well-defined, we must first check that the morphisms $g_t^e$ are morphisms of hypotheses. For every edge $e\in E(T)$ and incident node $t$, Lemma~\ref{lem:helper_1} shows that $g_t^e$ is a morphism $\beta_t[G_t]\to\beta_e[G_e]$. Thus $\Ta'$ is a co-decomposition of hypotheses.

    It remains to identify the limit of $\Ta'$. Since the limit of a tree-shaped diagram can be computed by iterated pullbacks, it suffices to prove the result for a single pullback. To that end, consider a cospan $G\xrightarrow{g}K\xleftarrow{h}H$, an assignment $\beta:[\ell]\to\pts{G\times_KH}$, and its induced hypothesis $\beta[G\times_KH]$. Write $h^\# :G\times_KH\to G$ and $g^\#:G\times_KH\to H$ for the pullback projections. By Lemma~\ref{lem:helper_1}, the following diagram commutes.
    \textbf{Note} in the rest of this proof, we use overline notation to denote the head of a hypothesis.
    
\[\begin{tikzcd}[sep = small]
	{\overline{\beta[G \times_{K} H]}} &&&&& \\
	&& {\overline{(g^{\#} \circ \beta)[H]} \times_{\overline{(h \circ g^{\#} \circ \beta)[K]}} \overline{(h^{\#} \circ \beta)[G]}} && {\overline{(h^{\#} \circ \beta)[G]}} \\
	&&& {\overline{(g^{\#} \circ \beta)[H]}} && {\overline{(h \circ g^{\#} \circ \beta)[K]}} \\
	&& {G \times_{K} H} && G \\
	&&& H && K
	\arrow[curve={height=-18pt}, from=1-1, to=2-5]
	\arrow[curve={height=30pt}, from=1-1, to=3-4]
	\arrow[curve={height=24pt}, from=1-1, to=4-3]
	\arrow[from=2-3, to=2-5]
	\arrow[from=2-3, to=3-4]
	\arrow[from=2-5, to=3-6]
	\arrow[from=2-5, to=4-5]
	\arrow[from=3-4, to=3-6]
	\arrow[from=3-4, to=5-4]
	\arrow[from=3-6, to=5-6]
	\arrow["{h^{\#}}"{description}, from=4-3, to=4-5]
	\arrow["{g^{\#}}"{description}, from=4-3, to=5-4]
	\arrow["g"{description}, from=4-5, to=5-6]
	\arrow["h"{description}, from=5-4, to=5-6]
\end{tikzcd}\]
    The two morphisms from $\overline{\beta[G\times_KH]}$ to the induced hypotheses on $G$ and $H$ have the same composite with the induced hypothesis on $K$. They therefore form a cone over the pullback diagram of heads. Thus, by the universal property of this pullback, there is a unique homomorphism
    \[
    c \colon \overline{\beta[G\times_KH]} \to
    \overline{(h^\#\circ\beta)[G]} \times_{\overline{(g\circ h^\#\circ\beta)[K]}} \overline{(g^\#\circ\beta)[H]}
    \]
    making the upper triangles commute. This is shown dotted in the following diagram.
\[\begin{tikzcd}[sep = small]
	{\overline{\beta[G \times_{K} H]}} &&&&& \\
	&& {\overline{(g^{\#} \circ \beta)[H]} \times_{\overline{(h \circ g^{\#} \circ \beta)[K]}} \overline{(h^{\#} \circ \beta)[G]}} && {\overline{(h^{\#} \circ \beta)[G]}} \\
	&&& {\overline{(g^{\#} \circ \beta)[H]}} && {\overline{(h \circ g^{\#} \circ \beta)[K]}} \\
	&& {G \times_{K} H} && G \\
	&&& H && K
	\arrow[curve={height=-6pt}, dashed, from=1-1, to=2-3]
	\arrow[curve={height=-18pt}, from=1-1, to=2-5]
	\arrow[curve={height=30pt}, from=1-1, to=3-4]
	\arrow[curve={height=24pt}, from=1-1, to=4-3]
	\arrow[from=2-3, to=2-5]
	\arrow[from=2-3, to=3-4]
	\arrow[from=2-5, to=3-6]
	\arrow[from=2-5, to=4-5]
	\arrow[from=3-4, to=3-6]
	\arrow[from=3-4, to=5-4]
	\arrow[from=3-6, to=5-6]
	\arrow["{h^{\#}}"{description}, from=4-3, to=4-5]
	\arrow["{g^{\#}}"{description}, from=4-3, to=5-4]
	\arrow["g"{description}, from=4-5, to=5-6]
	\arrow["h"{description}, from=5-4, to=5-6]
\end{tikzcd}\]
    By Lemma~\ref{lemma: limits of hypotheses}, the pullback of the two induced hypotheses over the induced hypothesis on $K$ is itself a hypothesis $\overline{(h^\#\circ\beta)[G]} \times_{\overline{(g\circ h^\#\circ\beta)[K]}} \overline{(g^\#\circ\beta)[H]}
    \to G\times_KH$; this is the red, dashed arrow below.
\[\begin{tikzcd}[row sep=small]
	{\overline{(G \times_{K} H)[\beta]}} &&&&& \\
	&& {\overline{H[g^{\#} \circ \beta]} \times_{\overline{K[h \circ g^{\#} \circ \beta]}} \overline{G[h^{\#} \circ \beta]}} && {\overline{G[h^{\#} \circ \beta]}} \\
	&&& {\overline{H[g^{\#} \circ \beta]}} && {K[h \circ g^{\#} \circ \beta]} \\
	&& {G \times_{K} H} && G \\
	&&& H && K
	\arrow[curve={height=-6pt}, dashed, from=1-1, to=2-3]
	\arrow[curve={height=-18pt}, from=1-1, to=2-5]
	\arrow[curve={height=30pt}, from=1-1, to=3-4]
	\arrow[curve={height=24pt}, from=1-1, to=4-3]
	\arrow[from=2-3, to=2-5]
	\arrow[from=2-3, to=3-4]
	\arrow[draw={rgb,255:red,214;green,92;blue,92}, dashed, from=2-3, to=4-3]
	\arrow[from=2-5, to=3-6]
	\arrow[from=2-5, to=4-5]
	\arrow[from=3-4, to=3-6]
	\arrow[from=3-4, to=5-4]
	\arrow[from=3-6, to=5-6]
	\arrow["{h^{\#}}"{description}, from=4-3, to=4-5]
	\arrow["{g^{\#}}"{description}, from=4-3, to=5-4]
	\arrow["g"{description}, from=4-5, to=5-6]
	\arrow["h"{description}, from=5-4, to=5-6]
\end{tikzcd}\]
    By Proposition~\ref{prop: Universal property of induced hypotheses}, the existence of the red dashed arrow implies the existence of the arrow $q$ below (depicted in blue below).
\[\begin{tikzcd}[sep = small]
	{\overline{\beta[G \times_{K} H]}} &&&&& \\
	&& {\overline{(g^{\#} \circ \beta)[H]} \times_{\overline{(h \circ g^{\#} \circ \beta)[K]}} \overline{(h^{\#} \circ \beta)[G]}} && {\overline{(h^{\#} \circ \beta)[G]}} \\
	&&& {\overline{(g^{\#} \circ \beta)[H]}} && {\overline{(h \circ g^{\#} \circ \beta)[K]}} \\
	&& {G \times_{K} H} && G \\
	&&& H && K
	\arrow["c"{description}, curve={height=-6pt}, dashed, from=1-1, to=2-3]
	\arrow[curve={height=-18pt}, from=1-1, to=2-5]
	\arrow[curve={height=30pt}, from=1-1, to=3-4]
	\arrow[curve={height=24pt}, from=1-1, to=4-3]
	\arrow["q"{description, pos=0.8}, draw={rgb,255:red,74;green,80;blue,232}, curve={height=-6pt}, dotted, from=2-3, to=1-1]
	\arrow[from=2-3, to=2-5]
	\arrow[from=2-3, to=3-4]
	\arrow[draw={rgb,255:red,214;green,92;blue,92}, dashed, from=2-3, to=4-3]
	\arrow[from=2-5, to=3-6]
	\arrow[from=2-5, to=4-5]
	\arrow[from=3-4, to=3-6]
	\arrow[from=3-4, to=5-4]
	\arrow[from=3-6, to=5-6]
	\arrow["{h^{\#}}"{description}, from=4-3, to=4-5]
	\arrow["{g^{\#}}"{description}, from=4-3, to=5-4]
	\arrow["g"{description}, from=4-5, to=5-6]
	\arrow["h"{description}, from=5-4, to=5-6]
\end{tikzcd}\]
    It remains to check that $c$ and $q$ are inverse. Under the canonical identification of both point sets with $[\ell]$, the underlying functions of $c$ and $q$ are both $\mathsf{id}_{[\ell]}$. Hence the underlying functions of $q\circ c$ and $c\circ q$ are identity functions. Since homomorphisms of relational structures are determined by their underlying functions, this implies that $\beta[G\times_KH] \cong (h^\#\circ\beta)[G] \times_{(g\circ h^\#\circ\beta)[K]} (g^\#\circ\beta)[H].$
\end{proof}

\subsection{Computing Supported Hypotheses by Dynamic Programming}\label{sec: exists FO algorithms nitty-gritty}

\body{
    It remains to show how to efficiently calculate the supported hypotheses at a given bag. To that end, given a rooted tree co-decomposition $\Ta$ and any node $t$ therein, denote by $H(t)$ the set of hypotheses on $G_t$ and by $\Ha_{supp}(t) := \{\nu \colon N \to G_t \mid \nu \text{ is a supported hypothesis at } t\}$ the set of supported hypotheses at $t$. 
}\label{para: set of supported hypotheses for dyn prog}

\begin{lemma}
    Let $\Ta$ be a rooted tree co-decomposition and let $i$ denote a leaf node without children. For an assignment $\beta_i$ to $G_i$, there is a single supported hypothesis for $\beta_i$, namely the induced hypothesis $\beta_i[G_i]$. Consequently $H(i) = \Ha_{supp}(i)$.
\end{lemma}
\begin{proof}
    The partial limit $G[i]$ at $i$ is identical to $G_i$.
\end{proof}

\begin{lemma}\label{lem:running_env}
Let $\Ta$ be a rooted tree co-decomposition, let $t$ be a node with unique child $c$, and write $e=ct$. Let $p_t\colon G_c\times_{G_e}G_t\rightarrow G_t$ be the pullback projection. Then $\Ha_{\mathrm{supp}}(t)$ can be computed from $\Ha_{\mathrm{supp}}(c)$ as follows. For $\rho\in\Ha_{\mathrm{supp}}(c)$, write $\beta_c$ for its underlying assignment. An assignment $\beta_t$ to $G_t$ is \define{compatible with $\rho$} if $g_c^e\circ\beta_c=g_t^e\circ\beta_t$. For such a compatible pair, write $\beta_e$ for this common assignment to $G_e$. Then\vspace{-.5em}
\[
\Ha_{\mathrm{supp}}(t)
=
\left\{
p_t\circ
\left(
\rho\times_{\beta_e[G_e]}\beta_t[G_t]
\right)
\;\middle|\;
\rho\in\Ha_{\mathrm{supp}}(c)
\text{ and }
\beta_t\text{ is compatible with }\rho
\right\}.
\]
\end{lemma}

\begin{proof}
To help keep track of variables, we will use hats on assignments and hypotheses on the partial limit.
Let $\rho\in\Ha_{\mathrm{supp}}(c)$ and let $\beta_t$ be compatible with $\rho$. Since $\rho$ is supported, there is an assignment $\widehat{\beta}_c\colon[\ell]\rightarrow\pts{G[c]}$ such that the induced hypothesis $\widehat{\rho}:=\widehat{\beta}_c[G[c]]$ projects to $\rho$. The assignments $\widehat{\beta}_c$ and $\beta_t$ agree on $G_e$, so they determine an assignment $\beta$ to $G[t]=G[c]\times_{G_e}G_t$. By Lemma~\ref{lem:lim_env_is_env_lim}, the hypothesis induced by $\beta$ on $G[t]$ is $\widehat{\rho}\times_{\beta_e[G_e]}\beta_t[G_t]$.

Consider also the pullback $\rho\times_{\beta_e[G_e]}\beta_t[G_t]$. Since projection leaves the head unchanged, $\widehat{\rho}$ and $\rho$ have the same head, and the two pullbacks also have the same head. For readability in the diagram, write $R$ for the common head of $\widehat{\rho}$ and $\rho$, write $I=\operatorname{head}(\beta_t[G_t])$ and $K=\operatorname{head}(\beta_e[G_e])$, and write $S$ for the common head of the two pullbacks. The canonical projection ${\color{red}q}\colon G[t]\rightarrow G_c\times_{G_e}G_t$ then makes the following diagram commute. \vspace{-0.5em}
\[
\begin{tikzcd}[row sep=small]
	&&&& S \\
	R &&&& {G[t]=G[c]\times_{G_e}G_t} \\
	{G[c]} && S \\
	R && {G_c\times_{G_e}G_t} && I \\
	{G_c} && K && {G_t} \\
	&& {G_e}
	\arrow["{\widehat{\rho}\times_{\beta_e[G_e]}\beta_t[G_t]}", from=1-5, to=2-5]
	\arrow["{\mathrm{id}_S}"', dashed, from=1-5, to=3-3]
	\arrow["{\widehat{\rho}}", from=2-1, to=3-1]
	\arrow["{\mathrm{id}_R}"', curve={height=18pt}, dashed, from=2-1, to=4-1]
	\arrow[from=2-5, to=3-1]
	\arrow["q", dashed, from=2-5, to=4-3, color=red]
	\arrow[curve={height=-30pt}, from=2-5, to=5-5]
	\arrow[curve={height=12pt}, from=3-1, to=5-1]
	\arrow["{\rho\times_{\beta_e[G_e]}\beta_t[G_t]}"', from=3-3, to=4-3]
	\arrow["\rho", from=4-1, to=5-1]
	\arrow[from=4-3, to=5-1]
	\arrow[from=4-3, to=5-5]
	\arrow["{\beta_t[G_t]}"', from=4-5, to=5-5]
	\arrow["{g_c^e}"', from=5-1, to=6-3]
	\arrow["{\beta_e[G_e]}"', from=5-3, to=6-3]
	\arrow["{g_t^e}", from=5-5, to=6-3]
\end{tikzcd}
\]

Thus $\rho\times_{\beta_e[G_e]}\beta_t[G_t]$ is the projection of the induced hypothesis $\widehat{\rho}\times_{\beta_e[G_e]}\beta_t[G_t]$ along $q$. It follows that $p_t\circ(\rho\times_{\beta_e[G_e]}\beta_t[G_t])$ is its projection to $G_t$, and hence belongs to $\Ha_{\mathrm{supp}}(t)$.

Conversely, let $\nu\in\Ha_{\mathrm{supp}}(t)$. By supportedness, $\nu$ is the projection to $G_t$ of $\beta[G[t]]$ for some assignment $\beta$ to $G[t]$. Let $\widehat{\beta}_c$ and $\beta_t$ be the projections of $\beta$ to $G[c]$ and $G_t$, respectively. Let $\widehat{\rho}=\widehat{\beta}_c[G[c]]$, and $\rho$ be the projection of $\widehat{\rho}$ to $G_c$. Then $\rho\in\Ha_{\mathrm{supp}}(c)$ and $\beta_t$ is compatible with $\rho$. By Lemma~\ref{lem:lim_env_is_env_lim}, the hypothesis $\beta[G[t]]$ is $\widehat{\rho}\times_{\beta_e[G_e]}\beta_t[G_t]$. The preceding diagram therefore gives $\nu=p_t\circ(\rho\times_{\beta_e[G_e]}\beta_t[G_t])$, as required.
\end{proof}

In the following, each $\nu_j\rightarrow\beta_t[G_t]$ is the canonical morphism to the induced hypothesis for its underlying assignment, and we use the canonical identification $G_t\times_{G_t}G_t\cong G_t$.
\begin{lemma}\label{lem:two_children}
Let $\Ta$ be a rooted tree co-decomposition and let $t$ be a node with exactly two children $c_1$ and $c_2$. For each $j\in\{1,2\}$, apply Lemma~\ref{lem:running_env} as if $c_j$ were the unique child of $t$, and denote the resulting set of hypotheses on $G_t$ by $\Ha_{c_j\to t}$. Then\vspace{-.5em}
\[
\Ha_{\mathrm{supp}}(t)
=
\left\{
\nu_1\times_{\beta_t[G_t]}\nu_2
\;\middle|\;
\nu_j\in\Ha_{c_j\to t}\text{ for }j\in\{1,2\},
\text{ and }\nu_1,\nu_2
\text{ have underlying assignment }\beta_t
\right\}.\vspace{-.5em}
\]
That is, $\Ha_{\mathrm{supp}}(t)$ is obtained by taking pairwise pullbacks of hypotheses in $\Ha_{c_1\to t}$ and $\Ha_{c_2\to t}$ with the same underlying assignment.
\end{lemma}
\begin{proof}
For $j\in\{1,2\}$, let $\Ta_j$ be the co-decomposition obtained by restricting $\Ta$ to $t$ and the subtree rooted at $c_j$, and write $G_j[t]$ for its limit. Thus $G[t]\cong G_1[t]\times_{G_t}G_2[t]$.

First, let $\nu_j\in\Ha_{c_j\to t}$ for $j\in\{1,2\}$, and suppose that both have underlying assignment $\beta_t$. Since $\nu_j$ is supported with respect to $\Ta_j$, there is an assignment $\widehat{\beta}_j\colon[\ell]\rightarrow\pts{G_j[t]}$ such that the induced hypothesis $\widehat{\nu}_j:=\widehat{\beta}_j[G_j[t]]$ projects to $\nu_j$. The assignments $\widehat{\beta}_1$ and $\widehat{\beta}_2$ agree on $G_t$, since both project to $\beta_t$, and therefore determine an assignment $\beta$ to $G[t]$.

By Lemma~\ref{lem:lim_env_is_env_lim}, the hypothesis induced by $\beta$ is $\widehat{\nu}_1\times_{\beta_t[G_t]}\widehat{\nu}_2$. By the universal property of pullbacks, its projection to $G_t$ is $\nu_1\times_{\beta_t[G_t]}\nu_2$. Hence this hypothesis belongs to $\Ha_{\mathrm{supp}}(t)$.

Conversely, let $\nu\in\Ha_{\mathrm{supp}}(t)$, witnessed by an assignment $\beta\colon[\ell]\rightarrow\pts{G[t]}$. For each $j\in\{1,2\}$, project $\beta$ to an assignment $\widehat{\beta}_j$ on $G_j[t]$, let $\widehat{\nu}_j=\widehat{\beta}_j[G_j[t]]$, and let $\nu_j$ be its projection to $G_t$. Then $\nu_j\in\Ha_{c_j\to t}$, and both $\nu_1$ and $\nu_2$ have the same underlying assignment $\beta_t$, namely the projection of $\beta$ to $G_t$. Applying Lemma~\ref{lem:lim_env_is_env_lim} and projecting to $G_t$ gives $\nu=\nu_1\times_{\beta_t[G_t]}\nu_2$, as required.
\end{proof}

\begin{proof}[Proof of Theorem~\ref{thm:exists-fo}]
Note that a hypothesis with underlying assignment $\beta \colon [\ell]\to\pts{G_t}$ is determined by which of the $q : =\sum_{R\in\sigma}\ell^{\operatorname{ar}(R)}$ possible atomic relations on $[\ell]$ hold in its head. Thus, since each bag has at most $w$ points, there are at most $w^\ell 2^q$ hypotheses on any given bag.

With this in mind, process the nodes of $T$ from the leaves towards the root. At a leaf $t$, the supported hypotheses are precisely the induced hypotheses $\beta_t[G_t]$. If $t$ has one child, compute the supported hypotheses $\Ha_{\mathrm{supp}}(t)$ from those of its child using Lemma~\ref{lem:running_env}. If $t$ has two children, first perform the single-child computation separately for each child and then combine the resulting sets using Lemma~\ref{lem:two_children}. Since $T$ is binary, these are all the possible cases. At a node with one child, we consider at most $w^{2\ell}2^q$ compatible pairs of hypotheses and assignments. At a node with two children, we additionally consider at most $w^{2\ell}2^{2q}$ pairs of hypotheses. Checking compatibility and computing pullbacks of hypotheses takes $O(\ell+q)$ time, since their heads have only $\ell$ points. Thus, the total time required to compute all supported-hypothesis tables is $O\left(b(\ell+q)w^{2\ell}2^{2q}\right)$. Finally, at the root $r$, we evaluate the quantifier-free part of $\phi$ on the distinguished variables of each head of a supported hypothesis. Since there are at most $w^\ell 2^q$ such hypotheses, this requires an additional $O(w^\ell 2^q|\phi|)$ time. By Lemma~\ref{lemma: support at root is enough}, this decides whether $\lim\Ta\models\phi$, giving the claimed running time $O\left( b(\ell+q)w^{2\ell}2^{2q} +w^\ell 2^q|\phi| \right).$ 
\end{proof}

%% file: sections/fo.tex
\section{The General Case: Model Checking for FO}\label{sec: full FO}

\body{
    Notice that our algorithm for $\exists$-FO used only four ingredients: (1) an induced-object construction, (2) projections to the bags (which behaves functorially), (3) pullbacks for combining compatible local states, and (4) a parameter-dependent bound on the number of states over each bag. In this section, we will carefully generalize each of these ingredients. Assignments will be replaced by \textit{strategies}, hypotheses by \textit{hypotheticals}, and induced hypotheses by \textit{induced hypotheticals}. The corresponding projection operations and pullback constructions will be developed alongside them. This is summarized in Table~\ref{tab:existential-to-full-fo}.
}

\begin{table}[htbp]
\centering
\caption{Dictionary between the dynamic programs for $\exists$-FO and full FO.}
\label{tab:existential-to-full-fo}

\begingroup
\setlength{\arrayrulewidth}{1.2pt}
\setlength{\tabcolsep}{1em}
\renewcommand{\arraystretch}{1.35}

\begin{tabular}{|c|c|}
\hline
\textbf{$\exists$-FO} & \textbf{Full FO}\\
\hline
assignment $\beta$ & strategy $S$\\
\hline
hypothesis $\nu$ & hypothetical $\mathbb H$\\
\hline
induced hypothesis $\beta[G]$ & induced hypothetical $S[G]$\\
\hline
underlying assignment of $\nu$ & underlying strategy $U(\mathbb H)$\\
\hline
projection of hypotheses & projection of hypotheticals\\
\hline
compatible assignments & compatible underlying strategies\\
\hline
pullbacks in $\mathbf{Hyp}_{\ell}^{\sigma}$ &
pullbacks in $\cat{HTh}_{m}^{\sigma}$\\
\hline
supported hypotheses & supported hypotheticals\\
\hline
tables of hypotheses & tables of reduced hypotheticals\\
\hline
$\operatorname{head}(\nu)\models\phi$ &
all hypotheses of $0$-hypotheticals satisfy $\psi$\\
\hline
\end{tabular}

\endgroup
\end{table}

\body{
    Once these generalizations have been established, the algorithmic result for full FO will fall out: we will be able to lift the algorithm of Theorem~\ref{thm:exists-fo} almost verbatim, using the dictionary in Table~\ref{tab:existential-to-full-fo}, to obtain an algorithm for full FO model checking. This is the following theorem (c.f.~Paragraph~\ref{para: explicit bounds} for the monstrous, but explicit bound depending only on $w$, $|\phi|$ and $b$).
}

\begin{theorem}[Theorem~\ref{thm: FO model checking on codecompositions}, restated]
    Let $\sigma$ be a fixed finite relational signature. There is an algorithm that, given an FO formula $\phi$ and a tree co-decomposition $\Ta$ of width at most $w$ with $b$ bags, decides whether $\lim\Ta\models\phi$ in time $f_\sigma(w,|\phi|)\cdot b$ for some computable function $f_\sigma$. 
\end{theorem}

\subsection{From Assignments to Strategies: A First Pass}\label{sec: from assignments to strategies}

\body{
    When moving from $\exists$-FO to full first-order logic, assignments no longer suffice. In the existential fragment, we need only choose points that satisfy the formula. However, in full first-order logic, existential choices may depend on the values chosen for earlier universally quantified variables. To deal with this dependence, we introduce \emph{strategies}. We point out that, although we favour our application-specific notation, what we define below can also be described as Skolemization.  
}

\begin{definition}
    A formula $\phi$ is in \define{$\forall\exists$-alternating prenex normal form} if there is a quantifier-free formula $\psi$ such that $\phi$ is written as: 
    $\phi = \forall x_1\exists y_1\cdots\forall x_\ell\exists y_\ell\, \psi(\vec{x},\vec{y}).$
\end{definition}

\body{
    For an $\ell$-variable formula $\phi$, we think of the quantifier prefix of $\phi$ as a game consisting of $\ell$ rounds and two players: universal and existential. In round $i$, the universal player chooses a value for $x_i$, after which the existential player chooses a response for $y_i$. A strategy with $m$ rounds remaining must specify an existential response, together with a strategy for the remaining $m-1$ rounds, for every possible universal choice.}
\body{
    It will be necessary later to project strategies along homomorphisms of structures. In this case, it might be that the moves of the players ``stack up'' on one another. To see what we mean by this, consider an epimorphism $f \colon H \twoheadrightarrow G$ and suppose, for the sake of example, that we are dealing with a formula with only two quantifiers in the form $\forall x \exists y \psi(x,y)$. A strategy on $H$ consists of a function $s \colon \pts{H} \to \pts{H}$ which represents the fact that for each point $x \in \pts{H}$ played by the universal player, the existential player responds with the point~$s(x)$. Now suppose we want to use this strategy to generate a strategy on $G$; what would this look like? Well, from the perspective of a point $w \in \pts{G}$, the strategy $s$ on $H$ is \textit{nondeterministic}: there are potentially many points $\{x_1, x_2 \dots\} \in f^{-1}(w)$ which project to $w$ under $f$, so how should we build a strategy on $G$ starting from a strategy $s$ on $H$? The solution is to think of $w$ as having many different \textit{personalities}: $w$ could behave like $x_1$ or it could behave like $x_2$, or $x_3$, and so on. Given any choice of such a personality then it is easy to know what a strategy on $G$ should be: let's say that $w$ `behaves' like $x_1$. Then, if the universal player picks $w$, the existential player picks $f(s(x_1))$. Packaging this intuition in a different language, one should think of personalities as \textit{sections} of the map $f$ (i.e. maps $ \sigma \colon G \hookrightarrow H$ such that $f \circ \sigma = \mathsf{id}_G$). Any such section determines a strategy for $G$ from the strategy on $H$, as in the following diagram. 
    \[\begin{tikzcd}
	H && H \\
	G && G
	\arrow["s", from=1-1, to=1-3]
	\arrow["f", two heads, from=1-1, to=2-1]
	\arrow["f", two heads, from=1-3, to=2-3]
	\arrow["\sigma", shift left=3, color={rgb,255:red,92;green,92;blue,214}, curve={height=-6pt}, hook', from=2-1, to=1-1]
	\arrow["{f \circ s \circ \sigma}", color={rgb,255:red,92;green,92;blue,214}, dotted, from=2-1, to=2-3]
\end{tikzcd}\]
In what follows, we will refer to this behavior as a \emph{nondeterministic response}.
}\label{para: personalities}

\body{
    We will need to fix some notation before defining the set of strategies on a structure $G$: denote \define{disjoint unions} (i.e. coproducts) by $\Sigma$, \define{cartesian products} (i.e. categorical products) by $\Pi$ and the \define{powerset operator} on $X$ either by $\Pa(X)$ or by $2^X$. Given a partial assignment $\beta \colon [k] \to \pts{G}$ for some $k < \ell$ and a point $x \in \pts{G}$, write $\beta \cdot x$ to denote \define{continuation of $\beta$ by $x$}, meaning the assignment obtained by {extending} the partial assignment $\beta$ so as to send the next point (the point with index $k+1$) to $x$: \[(\beta\cdot x) \colon i \mapsto \begin{cases} \beta(i) &\text{if } i \in [k]\\ x &\text{otherwise}. \end{cases}\]

}

\begin{definition}
    Given a partial assignment $\beta \colon [2(\ell - m)] \to \pts{G}$, then we recursively define $\operatorname{Strat}(G, \beta)$, the set of \define{$m$-strategies} relative to $\beta$ as follows:
    \[
\operatorname{Strat}_m(G,\beta)
:=
\begin{cases}
    \{\beta\} & \text{if } m=0, \\[1ex]
    \underbrace{
    \prod_{x\in\pts{G}}
        \underbrace{
        \Pa\left(
            \sum_{y\in\pts{G}}
            \underbrace{
                \operatorname{Strat}_{m-1}(G,\beta\cdot x\cdot y)
            }_{\substack{\text{continuations of the}\\
                         \text{assignment }\beta\cdot x\cdot y}}
        \right)
        }_{\text{nondeterministic response to }x}
    }_{\text{strategy from the assignment }\beta}
    & \text{otherwise.}
\end{cases}
\]
    Another way of thinking of the case where $m > 0$ is to see $\operatorname{Strat}_m(G, \beta)$ as the set of functions taking points $x$ in $G$ to non-empty subsets of $\sum_{y \in \pts{G}} \operatorname{Strat}_{m-1}(G, \beta \cdot x \cdot y)$. An \define{$m$-strategy relative to $\beta$} is an element of~$\operatorname{Strat}_m(G, \beta)$.
\end{definition}    

\begin{example}
    Given an assignment $\beta \colon [\ell] \to \pts{G}$, there is a single $0$-strategy relative to $\beta$, namely $\beta$ itself.
\end{example}

\subsection{From Hypotheses to Hypotheticals: A First Pass}

\body{
    We generalized assignments to strategies in Section~\ref{sec: from assignments to strategies}. Now we will generalize hypotheses. Recall that a hypothesis is roughly a way of decorating an assignment $\beta$ with local relational data that is ``seen'' by $\beta$ (c.f. \cref{def: hypothesis}). Now we wish to perfrom the same decoration on strategies. However, this will generate considerably more data. Recall from Paragraph~\ref{para: personalities} that we needed to keep track of the various potential personalities of each point: this forced us to define $\operatorname{Strat}_m(G, \beta)$ in terms of \textit{nondeterminstic responses}. Obviously, we will have to decorate each such response with a hypothesis; but there is more: to perform dynamic programming, we need to associate to each assignment generated by a strategy to a \textit{set} of hypotheses representing the local `guesses' at what the `correct' hypothesis ought to be. This is all deliberately vague, so as to convey intuition, but note that, despite being much easier, we already did something similar in the $\exists$FO case (c.f. Lemma~\ref{lemma: limits of hypotheses} and Paragraph~\ref{para: set of supported hypotheses for dyn prog}). Stating all of this formally, yields the definition of \textit{hypotheticals}.
}
\begin{definition}\label{def: hypothetical}
    Let $\operatorname{Hyp}(G,\beta)$ denote the set of hypotheses on $G$ with underlying assignment $\beta$.
    Given a partial assignment $\beta \colon [2(\ell - m)] \to \pts{G}$, then we recursively define the set of \define{potential hypotheticals} $\operatorname{HTh}(G, \beta)$ as follows:
    \begin{equation}\label{eqn: hypothetical}
        \operatorname{HTh}_m(G, \beta) := \begin{cases}
            2^{\operatorname{Hyp}(G, \beta)} &\text{if } m = 0, \\
            \prod_{x \in \pts{G}} \Pa(\Pa(\sum_{y \in \pts{G}} \operatorname{HTh}_{m-1}(G, \beta \cdot x \cdot y))) &\text{otherwise.}
        \end{cases}
    \end{equation}
\end{definition}

\body{
    Equation~\ref{eqn: hypothetical} can be a little hard to parse at first sight, so we restate it with the following annotations:
    \begin{equation}\label{eqn: hypothetical annotated}
    \operatorname{HTh}_m(G,\beta)
    :=
    \begin{cases}
        \underbrace{
            2^{\operatorname{Hyp}(G,\beta)}
        }_{\text{sets of hypotheses extending }\beta}
        & \text{if } m=0,\vspace{1.2em}
        \\ 
        \underbrace{
        \prod_{x\in\pts{G}}
            \underbrace{
            2^{
                \underbrace{
                2^{
                    \sum_{y\in\pts{G}}
                    \underbrace{
                        \operatorname{HTh}_{m-1}
                        (G,\beta\cdot x\cdot y)
                    }_{\substack{\text{continuations from the}\\
                                 \text{assignment }\beta\cdot x\cdot y}}
                }
                }_{\text{non-determinisitic response for personality of }x}
            }
            }_{\text{responses for the different personalities of }x}
        }_{\text{potential hypothetical from }\beta}
        & \text{otherwise.}
    \end{cases}
\end{equation}
}

\begin{example}
    We can view a $1$-hypothetical as a function $f:\pts{G}\to2^{2^{\pts{G}}}$ along with, for each pair $(x,y\in P\in f(x))$, a set of hypotheses with $(x,y)$ as their underlying assignment.
\end{example}

\body{
    A hypothetical is said to be \define{induced} if all of the hypotheses it contains are induced.
}
\begin{lemma}
    Let $G$ be a structure and $\phi = \forall x_1 \exists y_1 \ldots\forall x_\ell\exists y_\ell \psi(\vec{x},\vec{y})$ a formula in $\forall\exists$-alternating prenex normal form. Then $G\models\phi$ if and only if there is an induced $\ell$-hypothetical all whose hypotheses satisfy~$\psi$.
\end{lemma}
\begin{proof}
    Let $\mathbb{H}$ be an induced $\ell$-hypothetical all whose hypotheses satisfy~$\psi$. For each $\ell$-assignment $\beta$ in (the underlying strategy) of $\mathbb{H}$, there is exactly one hypothesis that corresponds to it, namely $\beta[G]$. By definition, $\beta[G]$ records exactly the relations between the points in $G$ specified by $\beta$.
\end{proof}

\subsection{Categories of Strategies and Hypotheticals}

\body{In the previous sections, strategies and hypotheticals were introduced set-theoretically. This was useful for getting a handle on the objects we were defining. However, it quickly becomes a headache when we begin to manipulate them recursively. In this section, we recast both constructions categorically. To do so, we will need to assume some further familiarity with categories; specifically, notions such as \textit{adjunctions} (in particular, \textit{unit}, \textit{counit}), \textit{limit-preserving} and \textit{categories of functors}. All of these can be found, for example, in Riehl~\cite{Riehl2016} or Awodey~\cite{Awodey2010}. Without being too explicit, we make some use of the machinery of (op)fibrations, which can be in Jacobs~\cite{JacobsCLTT}). We do this for four reasons which we will outline below.} 

\subsubsection{A Categorical Apology}

\body{The first advantage of promoting from sets of strategies and hypotheticals to categories thereof is that we can express these recursively-defined constructions using standard categorical operations such as products, sums, and categories of families. We do this so that we do not have to reinvent the wheel: the behaviour of these operations with respect to limits and functors is already well understood and can be found in standard texts. As a case in point, throughout this section we will often draw upon Jacobs' textbook~\cite{JacobsCLTT}. In practical terms, besides allowing us to draw upon known results, the categorical framing means that almost all the constructions in this section no longer have to be defined by hand through several layers of nested data, but instead arise canonically from the corresponding categorical operations.}

\body{The second advantage concerns the relationship between underlying and induced constructions. In our warmup, where we showed how to perform $\exists$-FO model checking in FPT-time, we pointed out two obvious facts: every hypothesis has an underlying assignment and every assignment gives rise to a canonical induced hypothesis. Furthermore, we used relationships between these two constructions, together with a universal property of induced hypotheses (c.f.~Prop.~\ref{prop: Universal property of induced hypotheses}), to simplify our arguments. Here, we wish to extend this approach: it will turn out to be very useful to establish a correspondence between underlying strategies and induced hypotheticals, and this correspondence will be expressed by an \textit{adjunction} (see Awodey~\cite{Awodey2010} or Riehl~\cite{Riehl2016} for the definition) between the corresponding categories. As usual, the benefit of having an adjunction is that it provides canonical comparison maps precisely where the set-theoretic presentation would otherwise require repeated bookkeeping.}

\body{The third advantage is closely related. In the previous sections, essentially all of our constructions could be made explicitly. In the present setting the relevant constructions are organized only up to isomorphism\footnote{For the category theorist, the keyword to keep in mind is: \textit{pseudofunctor}. This, and our reference to Jacobs and fibered constructions should already alert the category theorist to the kinds of manipulations we will be employing.}. If we were to organize this data combinatorially, we would be have to keep track of all these isomorphisms by hand. In contrast, the categorical formulation, and in particular the adjunctions introduced below, gives us a systematic way of handling these comparison isomorphisms rather than constructing and checking them separately at every stage of the recursion.}

\body{The fourth and final advantage (perhaps the most important one) is that the categorical presentation is very useful when it comes to proving that the relevant categories admit pullbacks. This is the key property that allows us to present our algorithm for full FO-model checking in essentially the same form as we did in the case of $\exists$-FO: if the objects are defined in the right way, then both the proof and the resulting algorithm become canonical.}

\begin{intuition}
    It might be useful to the reader to keep the following in mind throughout the rest of this section: although one can define categories of strategies and hypotheticals -- we think of these as \textit{fibers} -- over each fixed structure $G$ -- we think of this as the \textit{base}, these categories will need to be \textit{compared with each other}; all at the same time; as if they lived in a larger, ambient category\footnote{A useful analogy to keep in mind is that of modules over rings: over any ring $R$, we can define a category $\cat{Mod}(R)$ of modules over $R$; however we wish to allow compare modules over \textit{different} rings and let them comunicate via appropriate morphisms.}. In keeping with standard category-theoretic terminology, we call these the \textit{total category of strategies} and the \textit{total category of hypotheticals}. Although we won't use it explicitly, in the back of minds and on our scratch paper we know that the canonical way of ``assembling these local categories into one, total category'' is via what is known as a \textit{Grothendieck construction}. The practical upshot of this perspective is that our proofs will often proceed by calling upon standard results from fibered category theory~\cite{JacobsCLTT, StreicherFibredCategories} and combining these with the recursive structures of strategies and hypotheticals.
\end{intuition}

\subsubsection{Families and Reduced Families}

\body{
    The definitions of strategies and hypotheticals rely on the following operations on sets: products, coproducts (disjoint unions) and powersets. When passing from sets to categories, products and coproducts are still well-defined, but some choices need to be made about how to represent powersets. We opt for \textit{families}, which we define below and which can basically be thought of as indexed sequences, since they behave nicely when combined with our other constructions and since all the properties we will require of them are already well established facts in fibered category theory (c.f.~Jacobs~\cite{JacobsCLTT}). 
}

\begin{definition}\label{def: C-families}
    Let $\cat{C}$ be a category. For a set $S$, write $S^*$ for the
    discrete category whose set of objects is~$S$ (that is to say: $S^*$ is the category whose objects are the points of $S$ and whose only morphisms are the identities). The category $\cat{Fam}(\cat{C})$ is the full\footnote{A full subcategory of a category includes all the objects of the subcategory and all morphisms between those objects.} subcategory of the category $\cat{Diag}(\cat{C})$ consisting of diagrams of the form $F\colon S^*\to\cat{C}$. We call such a diagram an \define{$S$-indexed family in $\cat{C}$} and usually denote it by $(F_s)_{s\in S}.$ A \define{morphism of indexed families} $(F_s)_{s\in S} \to (G_t)_{t\in T}$ consists of a pair $(u^* \colon S^* \to T^*, \eta\colon F\to G\circ u^*)$ where $u^*$ is a functor and $\eta$ is a natural transformation as in the following diagram. 
    \[
    \begin{tikzcd}[column sep=large,row sep=large]
        S^*
            \arrow[rr,"F"]
            \arrow[dr,"u^*"']
        &&
        \cat{C}
        \\
        &
        T^*
            \arrow[ur,"G"']
        &
    \end{tikzcd}
    \qquad
    \eta\colon F\to G\circ u^*.
    \]
    Composition and identities are inherited from the category of diagrams, namely identities are of the form $(id_S, (id_{S(s)})_{s\in S})$ and composition is defined by $(v,\theta)\circ(u,\eta) := \left( v\circ u,\, \bigl(\theta_{u(s)}\circ\eta_s\bigr)_{s\in S} \right).$
\end{definition}

\body{
    Spelling out \cref{def: C-families}: a morphism of $\cat{C}$-families $(F_s)_{s\in S} \to (G_t)_{t\in T}$ consists of a \textit{function of sets} $u\colon S\to T$ together with an $S$-indexed family of morphisms in $\cat{C}$ of the form $(\eta_s\colon F_s\to G_{u(s)})_{s\in S}.$ Furthermore, it is useful to be be explicit and observe that we allow $S=\varnothing$, so that $\cat{Fam}(\cat{C})$ contains the empty family. 
}

\begin{definition}\label{def: reduced families}
    For a a small category $\cat{C}$, any subset $A\subseteq\operatorname{Ob}(\cat{C})$ induces the canonical diagram
    \begin{align*}
        \iota_A\colon A^* &\to\cat{C} \\
        a &\mapsto a.
    \end{align*}
    We call any such diagram a \define{reduced family in $\cat{C}$}. We write $\cat{Fam}_{\mathrm{red}}(\cat{C})$ for the full subcategory of $\cat{Fam}(\cat{C})$ consisting of \define{reduced families}. 
\end{definition}

\body{
    With an eye towards our intended algorithmic goals, we need to account for a subtlety: we do not want to work with arbitrary families, since there are too many of these. Instead we want to work with \textit{reduced families}, since $\operatorname{Ob} \bigl(\cat{Fam}_{\mathrm{red}}(\cat{C})\bigr) \cong 2^{\operatorname{Ob}(\cat{C})}$ where $\cat{Fam}_{\mathrm{red}}(\cat{C})$. Having said that, however, it is far easier to use arbitrary families to \textit{prove} our theorems and then only worry about reduced families when we eventually get into the weeds of the algorithmics (specifically to prove our running-time bounds). 
}

\body{
    Categories of families, being categories of diagrams, are particularly well-behaved.
}

\begin{proposition}\label{prop:fam-pullbacks}
    If $\cat{C}$ is a category with pullbacks, then $\cat{Fam}(\cat{C})$ has pullbacks. 
\end{proposition}
\begin{proof}
    Although this is standard and can be found in any introductory text in category theory, we include the proof to build intuition. We will show that, when viewed as a full subcategory of $\cat{Diag}(\cat{C})$, $\cat{Fam}(\cat{C})$ is closed under pullbacks. This implies that the inclusion $\cat{Fam}(\cat{C}) \hookrightarrow \cat{Diag}(\cat{C})$ creates pullbacks and hence the result follows.
    
    To that end, consider a cospan $(S^*,F) \xrightarrow{(u,\varphi)} (R^*,L) \xleftarrow{(v,\psi)} (T^*,G)$ in $\cat{Fam}(\cat{C})$. As Peschke and Tholen show~\cite[Proposition~2.5 and Remark~2.8]{PeschkeTholenDiagrams}, pullbacks in $\cat{Diag}(\cat{C})$ are computed by first taking the pullback of the indexing categories and then taking the corresponding pointwise pullbacks in $\cat{C}$. Thus the desired pullback diagram is given by
    \begin{align*}
        F \times_L G \colon S^* \times_{R^*} T^* &\to \cat{C} \\
        (s,t) &\mapsto F(s)\times_{L(u(s))}G(t),
    \end{align*}
    where the pullback $F(s)\times_{L(u(s))}G(t)$ is taken over the components $\varphi_s\colon F(s)\to L(u(s))$ and $\psi_t\colon G(t)\to L(v(t))$. The fact that the discrete category functor $(-)^* \colon \cat{Set} \to \cat{Cat}$ is right adjoint~\cite[Remark 8.3.5]{riehl2014categoricalhomotopytheory}, and hence preserves limits~\cite{Riehl2016}, implies that $S^*\times_{R^*}T^* \cong (S\times_R T)^*$ and hence $F \times_L G$ is discrete (i.e. a family) as desired. 
\end{proof}

\body{We will use the following standard fact repeatedly in our proofs. It can be found, for instance, in Rosicky~\cite[Section~2]{AdamekRosickyFreeCompletions}. It will allow us to turn a summation into a product, which turns out to be the most convenient form to work with. 
}

\begin{proposition}\label{prop: fam of sum is product of fam}
    Let $(\cat{C}_i)_{i\in I}$ be a family of categories. There is a canonical equivalence
    \[
        \cat{Fam}\left(\sum_{i\in I}\cat{C}_i\right)
        \simeq
        \prod_{i\in I}\cat{Fam}(\cat{C}_i).
    \]
    
\end{proposition}

\subsubsection{A Second Pass at Strategies and Hypotheticals: Promotion to Categories}

\body{
    We can now recast the definitions of strategies and hypotheticals in terms of categories of families. We proceed in two steps. First, we define the corresponding \emph{diagrammatic} categories. For the same reasons as pointed out above, these categories will contain too many objects: our categories of families retain indexing data, whereas our original set-theoretic definitions do not distinguish between different enumerations of the same family. We will therefore later restrict attention to reduced families and show that, on objects, the resulting definitions recover the set-theoretic notions introduced above.
}

\begin{definition}\label{def: diagrammatic strategies}
    Given a partial assignment $\beta\colon[2(\ell-m)]\to\pts{G},$ define the \define{category of diagrammatic $m$-strategies relative to $\beta$ on $G$} recursively as follows:
    \begin{equation}\label{eqn: diagrammatic strategies}
    \cat{Strat}_m(G,\beta) :=
    \begin{cases}
        \underbrace{
            \mathbf{1}
        }_{\text{the $0$-strategy }\beta}
        &\text{if }m=0,\\[5.5ex]
        \underbrace{
        \prod_{x\in\pts{G}}
            \underbrace{
            \cat{Fam}\left(
                \sum_{y\in\pts{G}}
                \underbrace{
                    \cat{Strat}_{m-1}(G,\beta\cdot x\cdot y)
                }_{\substack{\text{continuations from the}\\
                             \text{assignment }\beta\cdot x\cdot y}}
            \right)
            }_{\text{family of responses to }x}
        }_{\text{diagrammatic strategy from }\beta}
        &\text{otherwise.}
    \end{cases}
    \end{equation}
\end{definition}

\body{
    The definition of diagrammatic hypotheticals is obtained by replacing the base case by families of hypotheses and, recursively, allowing a family of possible ``personalities'' for each move $x$. 
}

\begin{definition}\label{def: diagrammatic hypotheticals}
    Given a partial assignment $\beta\colon[2(\ell-m)]\to\pts{G},$ define the \define{category of diagrammatic $m$-hypotheticals relative to $\beta$ on $G$} recursively as follows:
    \begin{equation}\label{eqn: diagrammatic hypotheticals}
    \cat{HTh}_m(G,\beta) :=
    \begin{cases}
        \underbrace{
            \cat{Fam}\bigl(\cat{Hyp}(G,\beta)\bigr)
        }_{\text{families of hypotheses with underlying assignment }\beta}
        &\text{if }m=0,\\[5.5ex]
        \underbrace{
        \prod_{x\in\pts{G}}
            \underbrace{
            \cat{Fam}\left(
                \underbrace{
                \cat{Fam}\left(
                    \sum_{y\in\pts{G}}
                    \underbrace{
                        \cat{HTh}_{m-1}(G,\beta\cdot x\cdot y)
                    }_{\substack{\text{continuations from the}\\
                                 \text{ assignment }\beta\cdot x\cdot y}}
                \right)
                }_{\text{response of one personality of }x}
            \right)
            }_{\text{family of personalities of }x}
        }_{\text{diagrammatic hypothetical from }\beta}
        &\text{otherwise.}
    \end{cases}
    \end{equation}
\end{definition}

\subsubsection{Comparison With the Set-Theoretic Definitions}

\body{
    As a sanity check, we would like to be sure that, if we only use reduced families, then the sets of objects of our categories of strategies and hypotheticals are isomorphic to our earlier set-theoretic definitions. To reïterate, note that only difference is how we take the families: do we keep track of all possible enumerations of indexes (families), or do we just retain one canonnical family per subset (reduced families)?
}

\begin{definition}\label{def: reduced diagrammatic strategies and reduced diagrammatic hypotheticals}
    Given a partial assignment $\beta\colon[2(\ell-m)]\to\pts{G}$, define the category of \define{reduced diagrammatic $m$-strategies relative to $\beta$} recursively as follows:
    \begin{equation}\label{eqn: reduced diagrammatic strategies}
        \cat{Strat}^{\mathrm{red}}_m(G,\beta)
        := \begin{cases}
            \mathbf{1} &\text{if } m = 0,\\
            \displaystyle
            \prod_{x\in\pts{G}}
            \cat{Fam}_{\mathrm{red}}
            \left(
                \sum_{y\in\pts{G}}
                \cat{Strat}^{\mathrm{red}}_{m-1}
                (G,\beta\cdot x\cdot y)
            \right)
            &\text{otherwise.}
        \end{cases}
    \end{equation}
    Similarly, we define the category of \define{reduced diagrammatic $m$-hypotheticals relative to $\beta$} as:
    \begin{equation}\label{eqn: reduced diagrammatic hypotheticals}
        \cat{HTh}^{\mathrm{red}}_m(G,\beta)
        := \begin{cases}
            \cat{Fam}_{\mathrm{red}}\bigl(\cat{Hyp}(G,\beta)\bigr)
            &\text{if } m = 0,\\
            \displaystyle
            \prod_{x\in\pts{G}}
            \cat{Fam}_{\mathrm{red}}
            \left(
                \cat{Fam}_{\mathrm{red}}
                \left(
                    \sum_{y\in\pts{G}}
                    \cat{HTh}^{\mathrm{red}}_{m-1}
                    (G,\beta\cdot x\cdot y)
                \right)
            \right)
            &\text{otherwise.}
        \end{cases}
    \end{equation}
\end{definition}

\begin{proposition}\label{prop: set-theoretic and diagrammatic strategies agree}
    For every structure $G$, every $0\leq m\leq\ell$, and every partial assignment $\beta\colon[2(\ell-m)]\to\pts{G}$, there are a canonical bijections
    \begin{align*}
        \operatorname{Strat}_m(G,\beta)
        &\cong
        \operatorname{Ob}
        \bigl(
            \cat{Strat}^{\mathrm{red}}_m(G,\beta)
        \bigr) \text{ and}\\
        \operatorname{HTh}_m(G,\beta)
        &\cong
        \operatorname{Ob}
        \bigl(
            \cat{HTh}^{\mathrm{red}}_m(G,\beta)
        \bigr).
    \end{align*}
\end{proposition}
\begin{proof}
    The proof is routine and basically the same in both cases, so we just show it for strategies (below the isomorphisms follow because the sets have the same size). 
    
    If $m=0$, then $\operatorname{Ob} \bigl( \cat{Strat}^{\mathrm{red}}_0(G,\beta) \bigr) = \operatorname{Ob}(\mathbf{1}) \cong \operatorname{Strat}_0(G,\beta).$ Now, inducting on $m$, if $m>0$, then
    \[
        \operatorname{Ob}
        \bigl(
            \cat{Strat}^{\mathrm{red}}_m(G,\beta)
        \bigr)
        =
        \prod_{x\in\pts{G}}
        2^{
            \sum_{y\in\pts{G}}
            \operatorname{Ob}
            \bigl(
                \cat{Strat}^{\mathrm{red}}_{m-1}
                (G,\beta\cdot x\cdot y)
            \bigr)
        }\\
        \cong 
        \prod_{x\in\pts{G}}
        2^{
            \sum_{y\in\pts{G}}
            \operatorname{Strat}_{m-1}
            (G,\beta\cdot x\cdot y)
        }\\
        =:
        \operatorname{Strat}_m(G,\beta).
    \]
\end{proof}

\body{
    In what follows, it turns out to be easier to forget the requirement of being reduced and work directly with diagrammatic strategies and hypotheticals. For these, we will show that pullbacks exist and then note that pullbacks of reduced strategies (respectively hypotheticals) are themselves reduced. 
}

\subsubsection{Underlying Strategies and Induced Hypotheticals: A Useful Adjunction}

\begin{definition}
    Just as hypotheses have underlying assignments, there is a canonical \define{underlying strategy functor} $U\colon \cat{HTh}_m(G,\beta) \to \cat{Strat}_m(G,\beta)$ for every assignment $\beta$ on a structure $G$, which recursively forgets the hypothetical data and flattens out the family of personalities into a single family of responses. Formally, we have:
    \begin{equation}\label{eqn: underlying strategy functor}
        U :=
        \begin{cases}
            
            \cat{Fam}(\cat{Hyp}(G,\beta))
            \xrightarrow{!}
            \mathbf{1}
            &\text{if }m=0,\\[1.5ex]
            \displaystyle
            \prod_{x\in\pts{G}}
            \left(
                {\mu} \circ
                \cat{Fam}\left(
                    \cat{Fam}\left(
                        \sum_{y\in\pts{G}} U
                    \right)
                \right)
            \right)
            &\text{otherwise,}
        \end{cases}
    \end{equation}
    where ${ \mu \colon \cat{Fam}(\cat{Fam}(\cat C)) \to \cat{Fam}(\cat C)}$ is the canonical flattening functor (a family of a family is just a family in the same way that a list of lists can be flattened to a list). Functoriality follows immediately by induction from the functoriality of sums, products, $\cat{Fam}$, and $\mu$.
\end{definition}

\body{
The following proposition provides us with an analogue to induced hypotheses in our more general setting. They turn out to be right adjoints to the underlying strategy functor. Intuitively, we want to ``induce'' a hypothetical on a strategy by populating it with induced hypotheses.
}

\begin{proposition}\label{prop: induced hypothetical is right adjoint}
    The forgetful functor $U$ from hypotheticals to strategies has as right adjoint\footnote{Categorical aside: one can easily see that this induced hypothetical functor is monic and hence $\cat{Strat}_m(G, \beta)$ is a coreflective subcategory of $\cat{HTh}_{m}(G, \beta)$.} the induced hypothetical functor, namely we have the following adjunction for any $G$ and $\beta$. 
    \[\begin{tikzcd}[row sep=small]
	{\cat{HTh}_{m}(G, \beta)} & {\cat{Strat}_m(G, \beta)}
	\arrow[""{name=0, anchor=center, inner sep=0}, "U", shift left=2, from=1-1, to=1-2]
	\arrow[""{name=1, anchor=center, inner sep=0}, "{\_[G]}", shift left=2, from=1-2, to=1-1]
	\arrow["\dashv"{anchor=center, rotate=-90}, draw=none, from=0, to=1]
\end{tikzcd}\]
On objects, the right adjoint is defined as follows: for a strategy $S$, the induced hypothetical $S[G]$ is the hypothetical whose hypotheses are all induced on $G$. 
\end{proposition}
\begin{proof}
    By Proposition~\ref{prop: Universal property of induced hypotheses}, it is easy to see that induced hypotheses are right adjoint to the forgetful functor from hypotheses to assignments, summarized below. 
    \[\begin{tikzcd}[row sep=small]
	{\cat{Hyp}_m(G, \beta)} & {\pts{G}^{[2\ell - m]}}
	\arrow[""{name=0, anchor=center, inner sep=0}, "U", shift left=2, from=1-1, to=1-2]
	\arrow[""{name=1, anchor=center, inner sep=0}, "{\_[G]}", shift left=2, from=1-2, to=1-1]
	\arrow["\dashv"{anchor=center, rotate=-90}, draw=none, from=0, to=1]
\end{tikzcd}\]
    We shall ``lift'' this adjunction through the recursive definition of hypotheticals, so induct on $m$. For $m = 0$, $\cat{HTh}_{m}(G, \beta) = \cat{Fam}(\cat{Hyp}(G, \beta))$. For any family $S \xrightarrow{[N_s]_{s \in S}} \cat{Hyp}(G, \beta)$, applying the counit $\eta$ of the right adjoint $\_[G]$ pointwise, gives us a morphism $(\mathsf{id}_S, (\eta_s \colon N_s \to \beta[G]))$. The adjunction follows immediately. 

    Suppose $m > 0$. Equation~\ref{eqn: canonical equivalence} gives us: 
    \[
    \cat{HTh}_m(G,\beta) \simeq
        \prod_{x\in\pts{G}}
        \cat{Fam}\left(
            \prod_{y\in\pts{G}}
            \cat{Fam}\left(
                \cat{HTh}_{m-1}(G,\beta\cdot x\cdot y)
            \right)
        \right).
    \]
    From which one can see (similarly to how we will prove Lemma~\ref{lemma: hypothetical projection adjunction}) that the adjunction 
    \[\begin{tikzcd}[row sep=small]
	{\cat{HTh}_{m-1}(G,\beta\cdot x\cdot y)} & {\cat{Strat}_{m-1}(G,\beta\cdot x\cdot y),}
	\arrow[""{name=0, anchor=center, inner sep=0}, "U", shift left=2, from=1-1, to=1-2]
	\arrow[""{name=1, anchor=center, inner sep=0}, "R", shift left=2, from=1-2, to=1-1]
	\arrow["\dashv"{anchor=center, rotate=-90}, draw=none, from=0, to=1]
    \end{tikzcd}\]
    which we get by induction, can be lifted pointwise through families and products. 
\end{proof}

\subsubsection{Projections of Strategies and Hypotheticals}

\body{
    The proof of the following functoriality result uses standard facts in fibered category theory. More formally, one can find all of this, albeit not applied in this setting, in Jacobs' textbook where he covers families and reindexing~\cite[Section~1.2]{JacobsCLTT} and fibred coproducts~\cite[Section~1.9]{JacobsCLTT}. These results are not difficult to see on their own, however, given how nested the construction is, we prefer to point it all out below.
}

\begin{proposition}\label{prop:functoriality-of-projection}
    For any map $f\colon G\to H$ and partial assignment
    $\beta\colon[2(\ell-m)]\to\pts{G}$, projection along~$f$ defines the functors
    \[
        f_*\colon
        \cat{Strat}_m(G,\beta)
        \to
        \cat{Strat}_m(H,f\beta)
    \]
    and
    \[
        f_*\colon
        \cat{HTh}_m(G,\beta)
        \to
        \cat{HTh}_m(H,f\beta).
    \]
    We abuse notation by denoting both functors by $f_*$, since it will always be clear from context which is meant.
\end{proposition}
\begin{proof}
    We prove both statements simultaneously by induction on $m$. If $m=0$, the statement for strategies is immediate, while postcomposition with $f$ defines a functor 
    \[
        f^{\cat{Hyp}}_*\colon
        \cat{Hyp}(G,\beta)
        \to
        \cat{Hyp}(H,f\beta),
    \]
    and applying $\cat{Fam}$ gives the required functor on $0$-hypotheticals.

    Suppose $m>0$. Apply the projection functors supplied by the induction hypothesis to the $(m-1)$-strategies and $(m-1)$-hypotheticals, respectively, and then apply the categorical operations appearing in their recursive definitions---coproducts, reindexing along the fibers of $f$, products, and the relevant applications of $\cat{Fam}$. Since all of these operations are functorial, the result follows.
\end{proof}

\begin{example}\label{example: prop:functoriality-of-projection}
    Write $[C_i]_{i\in I}$ for an indexed family $I^*\to\cat{C}$ which sends each $i$ to some $C_i\in \cat{C}$. In particular, $[C]$
    denotes a singleton family. We also note that repeated entries are not
    identified.

    We will consider $1$-hypotheticals relative to some assignment $\beta$ on a structure $G$. For each $x,y\in\pts{G}$, choose a hypothesis $\nu_{x,y}\colon N_{x,y}\to G$ whose underlying assignment is
    $\beta\cdot x\cdot y$. Consider the hypothetical that, at each
    universal choice $x$, has one personality offering every
    $y\in\pts{G}$ as a response, with singleton continuation
    $[\nu_{x,y}]$; this looks like $\mathbb{H}:=\left( \Bigl[ \bigl[ [N_{x,y}\xrightarrow{\nu_{x,y}}G] \bigr]_{y\in\pts{G}} \Bigr] \right)_{x\in\pts{G}}\in\cat{HTh}_1(G,\beta)$. Let $f\colon G\to H$ be a homomorphism to some other structure $H$. The projection of $\mathbb{H}$ is $f_*\mathbb{H}\in\cat{HTh}_1(H,f\beta)$. Explicitly:
    \[
\Biggl(
    \;\;\;\underbrace{
        \Bigl[
            \;\;\underbrace{
                \bigl[
                        \;\underbrace{
                            [N_{x,y}\xrightarrow{\nu_{x,y}}G]
                        }_{\substack{\text{singleton family}\\
                                     \text{of hypotheses}}}\;
                \bigr]_{y\in\pts{G}}
            }_{\text{family of responses}}\;\;
        \Bigr]
    }_{\text{one personality}}\;\;\;
\Biggr)_{x\in\pts{G}}
\xmapsto{f_*}
\Biggl(
    \;\;\;\underbrace{
        \Bigl[
            \;\;\underbrace{
                \bigl[
                        \;\underbrace{
                            [N_{x,y}\xrightarrow{\nu_{x,y}}G
                                      \xrightarrow{f}H]
                        }_{\substack{\text{singleton family}\\
                                     \text{of hypotheses}}}\;
                \bigr]_{y\in\pts{G}}
            }_{\text{family of responses}}\;\;
        \Bigr]_{x\in f^{-1}(x')}
    }_{\text{family of personalities}}\;\;\;
\Biggr)_{x'\in\pts{H}}.
\]
    Projection postcomposes each hypothesis with $f$ and changes
    its response tag from $y$ to $f(y)$. For each universal choice
    $x'\in\pts{H}$, the outer family then contains one personality
    for each $x\in f^{-1}(x')$.

    Crucially, the response family remains indexed by $\pts{G}$:
    distinct responses with the same image under $f$ remain
    distinct entries. Likewise, personalities coming from
    different universal choices are concatenated, not merged.
    For example, when $G=\{a,b\}$ and $H=\{\ast\}$, projection
    produces two personalities at $\ast$, each containing two
    responses, all tagged by $\ast$.
\end{example}

\begin{remark}
        Given that our intended application is dynamic programming, it might seem troubling to be indexing part of the hypothetical in $H$ using using points from $G$ (see Example~\ref{example: prop:functoriality-of-projection}). However, this is precisely what reduced families (c.f.~\cref{def: reduced families}) are for: since the number of hypotheses on $H$ does \textit{not} depend on the size of $G$, we will have many duplicate copies which we can throw away using reduced families. In particular,  reduction proceeds from the inside out: first reduce the continuations, then remove duplicate responses within each personality, and finally remove duplicate personalities. The relevance of the duplicate entries will be felt when we start looking at morphisms between hypotheticals.
    \end{remark}

\subsubsection{Total categories of strategies and hypotheticals}

\body{
    Proposition~\ref{prop:functoriality-of-projection} tells us that the categories
    $\cat{Strat}_m(G,\beta)$ and $\cat{HTh}_m(G,\beta)$ can be compared along maps
    $G\to H$ and appropriate partial assignments. This motivates the definition of
    ``flattened out'' categories of strategies and hypotheticals, whose objects remember
    the structure and assignment over which they are defined, and whose morphisms allow
    us to transport strategies and hypotheticals between different structures. To do
    this, we first define the category of positions\footnote{This is nothing other than
    a slice category over $[2(\ell-m)]$.}.
}

\begin{definition}
    \label{def:category-of-positions}
    Fix $0\leq m\leq\ell$. The category $\cat{Pos}_m^\sigma$ of \define{$m$-positions} has as its objects pairs $(G,\beta)$, where $G$ is a $\sigma$-structure and $\beta\colon[2(\ell-m)]\to\pts{G}$ is an assignment. A morphism $f\colon(G,\beta)\to(H,\gamma)$ is a morphism $f\colon G\to H$ satisfying $f\beta=\gamma$, as in the following
    commuting triangle.
    \[
    \begin{tikzcd}[sep=small]
        && G \\
        {[2(\ell-m)]} \\
        && H
        \arrow["f", from=1-3, to=3-3]
        \arrow["\beta", from=2-1, to=1-3]
        \arrow["\gamma"', from=2-1, to=3-3]
    \end{tikzcd}
    \]
    Identities and composition are inherited from the category of $\sigma$-structures.
\end{definition}

\begin{figure}
    \centering
    \input{figures/grothendieck_construction}
    \caption{A morphism
    $({\color{red} f},{\color{red} \eta})\colon
    (G,\beta,\mathbb{A})\to(H,\gamma,\mathbb{B})$
    in a total category of strategies or hypotheticals
    (c.f.~\cref{def: total categories of strategies and hypotheticals}).}
    \label{fig: grothendieck construction}
\end{figure}
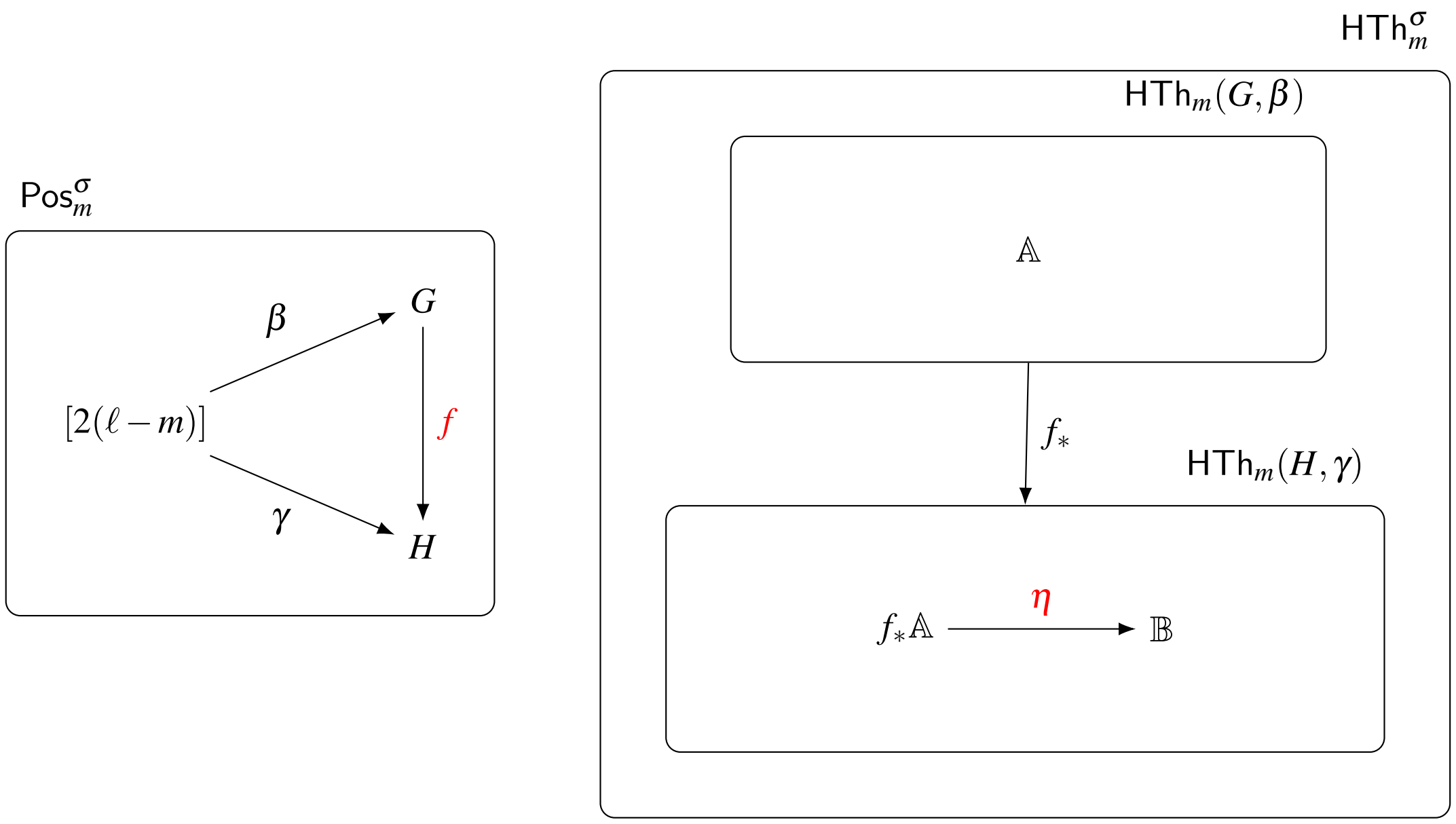

\body{
    One should think of each category $\cat{Strat}_m(G,\beta)$ and
    $\cat{HTh}_m(G,\beta)$ as ``sitting over'' the $m$-position~$(G,\beta)$. We now
    assemble each of these collections of categories into a single ambient category.
    We call the resulting categories the \textit{total category of strategies} and
    the \textit{total category of hypotheticals}. Their objects are strategies or
    hypotheticals equipped with the positions over which they are defined, while
    their morphisms keep track both of a change of position and of a morphism after
    transporting along that change. These are the categories in which we shall take
    pullbacks. We recommend consulting Figure~\ref{fig: grothendieck construction} while reading the following definition. 
}

\begin{definition}\label{def: total categories of strategies and hypotheticals}
    Fix $0\leq m\leq\ell$. For $\cat{X}\in\{\cat{Strat},\cat{HTh}\},$ define the \define{total category $\cat{X}_m^\sigma$} as follows.\footnote{ \textbf{Categorical aside.} We include this footnote for readers coming from a category-theoretic background; anyone unfamiliar with Grothendieck constructions or pseudofunctors may safely ignore it. Proposition~\ref{prop:functoriality-of-projection} can be strengthened by observing that, for both strategies and hypotheticals, there are canonical natural isomorphisms $(\mathsf{id}_G)_* \cong \mathsf{id}_{\cat{X}_m(G,\beta)}$ and $(gf)_*\cong g_*f_*$ for composable maps $f$ and $g$. Thus, for $\cat{X}\in\{\cat{Strat},\cat{HTh}\}$, the assignment $(G,\beta)\mapsto\cat{X}_m(G,\beta)$ together with the projection functors $f_*$ defines a pseudofunctor $\cat{Pos}_m^\sigma\to\cat{Cat}$. The  category $\cat{X}_m^\sigma$ defined here is precisely its Grothendieck  construction. In particular, there is a canonical projection functor  $p_m^{\cat{X}}\colon \cat{X}_m^\sigma\to\cat{Pos}_m^\sigma$ whose fiber over $(G,\beta)$ is $\cat{X}_m(G,\beta)$. Thus the total categories of strategies and hypotheticals are canonically opfibred over the category of positions. As usual, we suppress the canonical coherence isomorphisms in the formulas for identities and composition below. 
    }\vspace{-0.5em}
    \begin{itemize}
        \item An object of $\cat{X}_m^\sigma$ is a triple $(G,\beta,\mathbb{A}),$ where $(G,\beta)\in\cat{Pos}_m^\sigma$ and $\mathbb{A}\in
            \operatorname{Ob}\bigl(\cat{X}_m(G,\beta)\bigr)$.\vspace{-0.5em}
        \item A morphism $(f,\eta)\colon (G,\beta,\mathbb{A}) \to (H,\gamma,\mathbb{B})$ consists of a morphism of positions $f\colon(G,\beta)\to(H,\gamma)$  together with a morphism $\eta\colon f_*\mathbb{A}\to\mathbb{B}$ in $\cat{X}_m(H,\gamma)$.\vspace{-0.5em}
    \end{itemize}

    Its identities are
    $\bigl(\mathsf{id}_G,\mathsf{id}_{\mathbb{A}}\bigr)$ and it has composition $(g,\theta)\circ(f,\eta) := \bigl(gf,\theta\circ g_*(\eta)\bigr)$. We call $\cat{Strat}_m^\sigma$ the
    \define{total category of $m$-strategies} and we call $\cat{HTh}_m^\sigma$ the
    \define{total category of $m$-hypotheticals}.
\end{definition}

\body{
    A morphism in either total category is best visualized as in Figure~\ref{fig: grothendieck construction}: first ship the strategy or hypothetical along the underlying morphism of positions, and then compare the
    transported object with the target inside the category sitting over the codomain position. To interpret this in our language of co-decompositions, one should think of there being an auxiliary co-decomposition of strategies (resp. hypotheticals) whose bags and adhesions are the categories of $m$-strategies (resp. $m$-hypotheticals) sitting over the bags and adhesions of our co-decomposition of structures $\Ta$.
}

\body{
    The point of introducing these total categories is so that \textit{one can take pullbacks of strategies and hypotheticals}! The proof of the following is included in the next subsection.
}

\begin{theorem}\label{thm: total category has pullbacks}
    Fix $0\leq m\leq\ell$ and suppose that the category of positions $\cat{Pos}_m^\sigma$ has pullbacks. Then the total categories $\cat{Strat}_m^\sigma$ and $\cat{HTh}_m^\sigma$ have pullbacks, and their projection functors to $\cat{Pos}_m^\sigma$ preserve them.
\end{theorem}

\body{
The following two subsections are organised as follows: first, we prove \autoref{thm: total category has pullbacks}; and then we spell out the pullback object guaranteed by the theorem.
}

\subsubsection{Proof of \cref{thm: total category has pullbacks}: Pullbacks of Strategies and Hypotheticals}
\label{sec: pullbacks of strategies and hypotheticals}

\body{
    Establishing pullbacks in $\cat{HTh}_m^\sigma$ from first principles is nontrivial and indeed the construction (see Figure~\ref{fig: pullback in total category}) can be an eye-full. 
    \begin{figure}
        \centering
        \input{figures/pullback_in_total_category}
        \caption{Schematic view of the pullback construction in $\cat{HTh}_m^\sigma$. The lower diagram is the pullback of the underlying positions in $\cat{Pos}_m^\sigma$, while above each position sits its corresponding category of diagrammatic hypotheticals (the ``fibers''). The maps of positions induce the pushforward functors ($p_{1*}^{12}, p_{2*}^{12}, \pi_{1*}$ and $\pi_{2*}$) between the fibres. Each one of these functors has a right adjoint. We use these right adjoints to transport the cospan $p_{2*}^{12}\mathbb A_2 \xrightarrow{\zeta_1^{12}} \mathbb A_{12} \xleftarrow{\zeta_2^{12}}p_{1*}^{12}\mathbb A_1$ in $\cat{HTh}_m(G_{12},\beta_{12})$ back up to $\cat{HTh}_m(G _1 \times_{G_{12}} G_2,\beta_1 \times_{\beta_{12}} \beta_2)$. Using the fact that $\cat{HTh}_m(G _1 \times_{G_{12}} G_2,\beta_1 \times_{\beta_{12}} \beta_2)$ has pullbacks (c.f.~Proposition~\ref{lemma: fibers of HTh have pullbacks}), we can compute the pullback of this resulting cospan. This gives the pullback in the total category $\cat{HTh}_m^\sigma$.}
        \label{fig: pullback in total category}
    \end{figure}
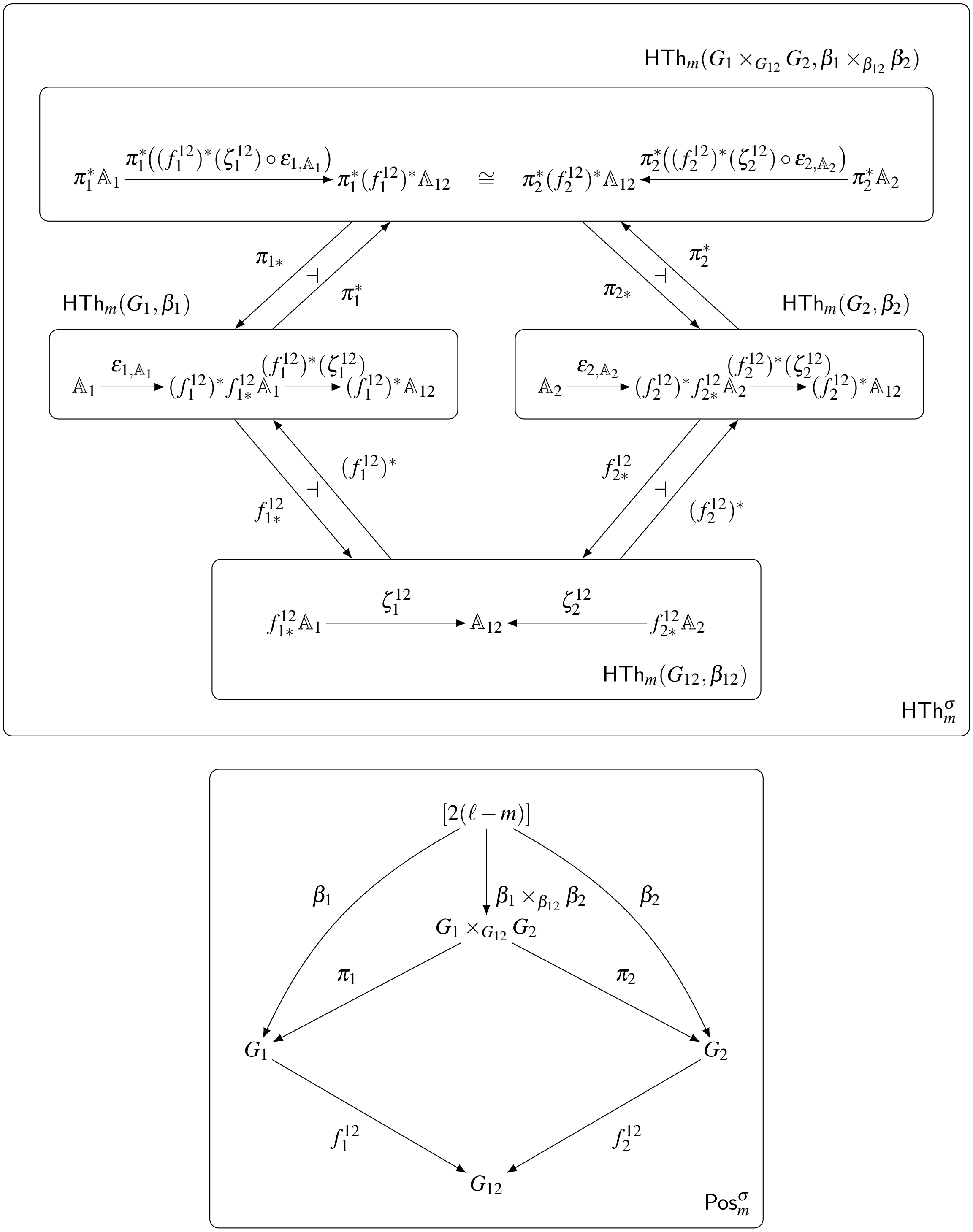
    Fortunately, much of the required machinery is known from fibered category theory and can be found in Jacobs' textbook~\cite{JacobsCLTT}. We will draw upon it as needed. The same argument applies to both strategies and hypotheticals. One can think of the total categories $\cat{Strat}_m^\sigma$ and $\cat{HTh}_m^\sigma$ as being ``fibered'' over $\cat{Pos}_m^\sigma$, with fibers $\cat{Strat}_m(G,\beta)$ and $\cat{HTh}_m(G,\beta)$ respectively. We will need to show two ingredients:\vspace{-0.5em}
    \begin{enumerate}
        \item the pushforward functor $f_*$ (Prop.~\ref{prop:functoriality-of-projection}) of a morphisms of positions $f$ admits a right adjoint $f^*$; and \vspace{-0.5em}
        \item the categories $\cat{Strat}_m(G,\beta)$ and $\cat{HTh}_m(G,\beta)$ for a position $(G,\beta)$, admit pullbacks. \vspace{-0.5em}
    \end{enumerate}
    From here, Theorem~\ref{thm: total category has pullbacks} will follow. To that end, we begin working toward proving the first point. Its proof will combine the known coproduct-reindexing adjunction for families~\cite[Lemma~1.9.5]{JacobsCLTT} with induction on the depth of a strategy or hypothetical.
}

\begin{lemma}\label{lemma: hypothetical projection adjunction}
    Let $\cat{X}\in\{\cat{Strat},\cat{HTh}\}$. Given any morphism of positions $f\colon(G,\beta)\to(H,f\beta)$, we have the following adjunction \vspace{-0.5em}
    \[\begin{tikzcd}
	{\cat{X}_m(G,\beta)} & {\cat{X}_m(H,f\beta)}
	\arrow[""{name=0, anchor=center, inner sep=0}, "{f_*}", shift left=2, from=1-1, to=1-2]
	\arrow[""{name=1, anchor=center, inner sep=0}, "{f^*}", shift left=2, from=1-2, to=1-1]
	\arrow["\dashv"{anchor=center, rotate=-90}, draw=none, from=0, to=1]
\end{tikzcd}\]
    where $f_*$ is the projection functor established in Proposition~\ref{prop:functoriality-of-projection}.
\end{lemma}

\begin{proof}
    Recall that applying $\cat{Fam}$ to an adjunction again gives
    an adjunction: the original adjunction is applied componentwise,
    without changing the function between the indexing sets. We prove both statements simultaneously by induction on $m$. 
    
    \textbf{Base case.} For $m=0$, the statement for strategies is immediate since
    $\cat{Strat}_0(G,\beta)=\mathbf{1}$. For hypotheticals, the existence of pullbacks of hypotheses (Lemma~\ref{lemma: limits of hypotheses}) gives a right adjoint $f^*$ to postcomposition $f_*\colon\cat{Hyp}(G,\beta)\to\cat{Hyp}(H,f\beta)$, namely $f^*(\nu) := \beta[G]\times_{(f\beta)[H]}\nu.$ From the universal property of this pullback, we have a canonical isomorphism of hom-sets
    \[
        {\cat{Hyp}(G,\beta)\bigl(N,f^*M\bigr)} \cong 
        {\cat{Hyp}(H,f\beta)\bigl(f_*N,M\bigr)}
    \]
    which establishes the desired adjunction. Applying $\cat{Fam}$, with unit and counit defined componentwise, gives the required adjunction on $0$-hypotheticals.

    \textbf{Inductive case.} Suppose $m>0$. Let $\cat{X}\in\{\cat{Strat},\cat{HTh}\}$. We will show that the adjunctions supplied by the inductive hypothesis survive the recursive construction of $f_*$. To that end, for each $x\in\pts{G}$, define the functor
    \begin{align}
        T_x\colon \sum_{y\in\pts{G}}\cat{X}_{m-1}(G,\beta\cdot x\cdot y)
        &\to \sum_{y'\in\pts{H}} \cat{X}_{m-1}(H,f\beta\cdot f(x)\cdot y') \\
        (y,\mathbb A)\quad  &\mapsto \quad (f(y),f_*\mathbb A).
    \end{align}
    We claim that
    \[
        \cat{Fam}(T_x)\colon
        \cat{Fam}\left(\sum_{y\in\pts{G}}\cat{X}_{m-1}(G,\beta\cdot x\cdot y)\right)
        \to
        \cat{Fam}\left(\sum_{y'\in\pts{H}}\cat{X}_{m-1}(H,f\beta\cdot f(x)\cdot y')\right)
    \]
    has a right adjoint. By  Proposition~\ref{prop: fam of sum is product of fam} we may use the canonical equivalence
    \[
        \cat{Fam}\left(\sum_{i\in I}\cat{C}_i\right)
        \simeq
        \prod_{i\in I}\cat{Fam}(\cat{C}_i)
    \]
    to equivalently rewrite $\cat{Fam}(T_x)$ as 
    \[
        \cat{Fam}(T_x)\colon
        \prod_{y\in\pts{G}}\cat{Fam}\bigl(\cat{X}_{m-1}(G,\beta\cdot x\cdot y)\bigr)
        \to
        \prod_{y'\in\pts{H}}\cat{Fam}\bigl(\cat{X}_{m-1}(H,f\beta\cdot f(x)\cdot y')\bigr).
    \]
    Under this equivalent interpretation, $\cat{Fam}(T_x)$ first applies the projection functors
    \begin{equation}\label{eqn: prop: hypothetical projection adjunction}
        L_y := \cat{Fam}(f_*)\colon
        \cat{Fam}\bigl(\cat{X}_{m-1}(G,\beta\cdot x\cdot y)\bigr)
        \to
        \cat{Fam}\bigl(\cat{X}_{m-1}(H,f\beta\cdot f(x)\cdot f(y))\bigr)
    \end{equation}
    componentwise and then, for each $y'\in\pts{H}$, concatenates the resulting families over $f^{-1}(y')$. Explicitly, it sends 
    \[
        (\mathbb A_y)_{y\in\pts{G}}
        \mapsto
        \left(
            \sum_{y\in f^{-1}(y')} L_y(\mathbb A_y)
        \right)_{y'\in\pts{H}}.
    \]
    Here the coproduct is taken in the category of families; it is given by co-producting (summing) indexing sets.

    By the induction hypothesis and preservation of adjunctions by $\cat{Fam}$, each $L_y$ has a right adjoint $R_y:=\cat{Fam}(f^*)$. For each $y'\in\pts{H}$, taking the coproduct over $f^{-1}(y')$ is left adjoint to the diagonal functor, which repeats a family once for each $y\in f^{-1}(y')$. This is the coproduct-reindexing adjunction for families~\cite[Lemma~1.9.5]{JacobsCLTT}. Consequently, under the same equivalences, a right adjoint to $\cat{Fam}(T_x)$ is given by
    \[
        (\mathbb B_{y'})_{y'\in\pts{H}}
        \mapsto
        \bigl(R_y(\mathbb B_{f(y)})\bigr)_{y\in\pts{G}}.
    \]
    That is, the pointwise application of the right-adjoints $R_y$ to the hypotheticals at the image indexes $f(y)$.

    For strategies, the projection functor on $m$-strategies applies these functors componentwise and then concatenates the resulting families over each preimage $f^{-1}(x')$. Applying the same coproduct-reindexing argument therefore supplies a right adjoint to $f_*$.

    For hypotheticals, applying $\cat{Fam}$ once more shows that $\cat{Fam}(\cat{Fam}(T_x))$ has a right adjoint. The projection functor on $m$-hypotheticals then applies these functors componentwise and concatenates the resulting outer families over each fibre $f^{-1}(x')$. Again applying the same coproduct-reindexing argument supplies a right adjoint to $f_*$.
\end{proof}

\body{
    The following example illustrates the right adjoint $f^*$ in the case of hypotheticals.
}

\begin{example}[The right adjoint to projection]
    \label{example: hypothetical projection adjunction}
    We use the notation for indexed families from
    Example~\ref{example: prop:functoriality-of-projection}.
    Take $m=\ell=1$, so that the initial assignment is empty and there is one universal choice followed by one existential response. Work in the empty relational signature (i.e. we work with \textit{sets}) and consider the unique map
    \[
        f\colon G=\{a,b\}\to H=\{\ast\}.
    \]

    Let $N$ be the structure with point set $[2]$. Since the signature is empty, there is only one hypothesis over each assignment, so let the unique hypothesis with head $N$ be the unique map $\nu\colon N\to H$ and, for each $x,y\in\pts{G}$, let $\nu_{x,y}\colon N\to G$ be the map sending $\nu_{x,y}(1)=x$ and $\nu_{x,y}(2)=y$. 
    (Note that $f\circ\nu_{x,y}=\nu$ for every $x,y$).

    Consider the hypothetical $\mathbb B$ on $H$ which, at its only universal choice $\ast$, has one personality offering the only existential response $\ast$, with singleton continuation $[N\xrightarrow{\nu}H]$. The right adjoint acts
    as follows:
    \[
    \Biggl(
        \;\;\;\underbrace{
            \Bigl[
                \;\;\underbrace{
                    \bigl[
                        \;\underbrace{
                            [N\xrightarrow{\nu}H]
                        }_{\substack{\text{singleton family}\\
                                     \text{of hypotheses}}}\;
                    \bigr]
                }_{\text{one response }\ast}\;\;
            \Bigr]
        }_{\text{one personality}}\;\;\;
    \Biggr)_{x'\in\pts{H}}
    \xmapsto{f^*}
    \Biggl(
        \;\;\;\underbrace{
            \Bigl[
                \;\;\underbrace{
                    \bigl[
                        \;\underbrace{
                            [N\xrightarrow{\nu_{x,y}}G]
                        }_{\substack{\text{singleton family}\\
                                     \text{of hypotheses}}}\;
                    \bigr]_{y\in\pts{G}}
                }_{\text{both lifts of the response }\ast}\;\;
            \Bigr]
        }_{\text{still one personality}}\;\;\;
    \Biggr)_{x\in\pts{G}}.
    \]
    To read this expression, first fix a universal choice $x\in\pts{G}$. Since $f(x)=\ast$, we take the family of personalities belonging to $\mathbb B(\ast)$. Its indexing set does not change: there is still exactly one personality at this fixed $x$.

    Within that personality, however, the response $\ast$ has two possible lifts, namely $a$ and $b$. The right adjoint retains both. For the lift $y$, it applies the depth-zero right adjoint to the continuation. That is, we obtain the family of pullbacks of $\nu$ with the induced hypothetical on the lifted assignment $x\cdot y$: giving $[N\xrightarrow{\nu_{x,y}}G]$. Spelling out the two universal choices separately, we obtain
    \[
    \begin{aligned}
        (f^*\mathbb B)(a)
        &=
        \Bigl[
            \bigl[
                [N\xrightarrow{\nu_{a,a}}G],\,
                [N\xrightarrow{\nu_{a,b}}G]
            \bigr]
        \Bigr],
        \\[1.5ex]
        (f^*\mathbb B)(b)
        &=
        \Bigl[
            \bigl[
                [N\xrightarrow{\nu_{b,a}}G],\,
                [N\xrightarrow{\nu_{b,b}}G]
            \bigr]
        \Bigr].
    \end{aligned}
    \]

    For example, the second response at universal choice $a$ chooses $b$ and leaves the completed assignment~$(a,b)$.

    The two universal choices therefore each have one personality, and that personality offers two responses. The lifts of the existential response become entries of the inner family. Also, these two response entries are distinct: their underlying assignments have different last coordinates, so reduction does not identify them.
\end{example}

\body{
    Now that we have proven that there is a right adjoint to the projection functors, we are ready to establish the easier fact that the fibers sitting over each position, that is categories of the form $\cat{Strat}_m(G,\beta)$ and $\cat{HTh}_m(G,\beta)$, have pullbacks. The proof will once again be an induction on $m$.
}

\begin{lemma}\label{lemma: fibers of HTh have pullbacks}
    For every structure $G$, every $0\leq m\leq\ell$, and every partial assignment $\beta\colon[2(\ell-m)] \to \pts{G}$, the categories $\cat{Strat}_m(G,\beta)$ and $\cat{HTh}_m(G,\beta)$ admit pullbacks.
\end{lemma}

\begin{proof}
    We prove both statements simultaneously by induction on $m$. For $m=0$, $\cat{Strat}_0(G,\beta)=\mathbf{1}$ has pullbacks, while the hypothetical case follows from Lemma~\ref{lemma: limits of hypotheses} and Proposition~\ref{prop:fam-pullbacks}.

    For $m>0$, Proposition~\ref{prop: fam of sum is product of fam} gives the canonical equivalences
    \begin{equation}
        \cat{Strat}_m(G,\beta)
        := \prod_{x\in\pts{G}}
        \cat{Fam}\left(
            \sum_{y\in\pts{G}}
            \cat{Strat}_{m-1}(G,\beta\cdot x\cdot y)
        \right) 
        \simeq
        \prod_{x\in\pts{G}}
        \prod_{y\in\pts{G}}
        \cat{Fam}\left(
            \cat{Strat}_{m-1}(G,\beta\cdot x\cdot y)
        \right),
    \end{equation} \label{eqn: canonical equivalence strategies}
    
    and
    \begin{align}
        \cat{HTh}_m(G,\beta)
        &:= \prod_{x\in\pts{G}} \cat{Fam}\left( \cat{Fam}\left( \sum_{y\in\pts{G}} \cat{HTh}_{m-1} (G,\beta\cdot x\cdot y) \right) \right) \nonumber \\ 
        &\simeq
        \prod_{x\in\pts{G}}
        \cat{Fam}\left(
            \prod_{y\in\pts{G}}
            \cat{Fam}\left(
                \cat{HTh}_{m-1}(G,\beta\cdot x\cdot y)
            \right)
        \right). \label{eqn: canonical equivalence}
    \end{align}
    By induction, each $\cat{Strat}_{m-1}(G,\beta\cdot x\cdot y)$ and $\cat{HTh}_{m-1}(G,\beta\cdot x\cdot y)$ has pullbacks. Products of categories with pullbacks have pullbacks componentwise, and $\cat{Fam}$ preserves the existence of pullbacks by Proposition~\ref{prop:fam-pullbacks}. 
\end{proof}

\begin{remark}
    The reader will notice in the above proof that the awful set-theoretic bookkeeping becomes a pointwise computation of pullbacks: after all, this was the whole point of passing to categories of families. 
\end{remark}

\body{
    At this point we are ready to prove Theorem~\ref{thm: total category has pullbacks}. As we mentioned earlier, by known results in fibered category theory~\cite[Theorem~8.3]{StreicherFibredCategories}, once one establishes lemmas of the form of Lemma~\ref{lemma: hypothetical projection adjunction} and Lemma~\ref{lemma: fibers of HTh have pullbacks}, existence of pullbacks in the total categories is guaranteed. 
}

\begin{theorem}[Theorem~8.3,~\cite{StreicherFibredCategories}] \label{thm: Streicher pullbacks}
    Let $p\colon\cat{E}\to\cat{B}$ be a fibration, where $\cat{B}$ has pullbacks. Then $\cat{E}$ has pullbacks and $p$ preserves them if and only if every fibre has pullbacks and every reindexing functor preserves them.
\end{theorem}

\begin{proof}[Proof of Theorem~\ref{thm: total category has pullbacks}]
    By Lemma~\ref{lemma: hypothetical projection adjunction}, the opfibrations defining $\cat{Strat}_m^\sigma$ and $\cat{HTh}_m^\sigma$ are also fibrations whose reindexing functors $f^*$ preserve pullbacks, being right adjoints. Furthermore, their fibres have pullbacks by Lemma~\ref{lemma: fibers of HTh have pullbacks}; hence Theorem~\ref{thm: Streicher pullbacks} applies to both.
\end{proof}

\body{
    Despite having proven the desired theorem and despite the fact that we have already given the high-level picture in Figure~\ref{fig: pullback in total category}, we find it clarifying to spell out the construction explicitly, step by step. The reader familiar with fibered categories and who is satisfied with the black-box proof of Theorem~\ref{thm: total category has pullbacks} can safely skip ahead past Section~\ref{sec: explaining pullbacks in total category}. The following section will focus only on the total category of hypotheticals, since the very same reasoning can be applied to the total category of strategies.
}

\subsubsection{Spelling out the Pullback Object Guaranteed by Theorem~\ref{thm: total category has pullbacks}}\label{sec: explaining pullbacks in total category}

\body{
    Throughout we recommend `zooming out' by glancing sporadically at Figure~\ref{fig: pullback in total category}.
    For any $m$, consider a cospan in the total category of $m$-hypotheticals; the following one, say: 
    \begin{equation}\label{eqn:  spelling out pullbacks in the total category}
        (G_1, \beta_1, \mathbb{A}_2) \xrightarrow{(f_1^{12}, \xi_{2}^{12})} (G_{12}, \beta_{12}, \mathbb{A}_{12}) \xleftarrow{(f_1^{12}, \xi_{2}^{12})} (G_2, \beta_2, \mathbb{A}_2).
    \end{equation}
     By how morphisms are defined in the total category of hypotheticals (c.f.~Figure~\ref{fig: grothendieck construction}), this determines the following cospan in the category $\cat{HTh}_m(G_{12}, \beta_{12})$: 
    \[(f_1^{12})_*(\mathbb{A}_1) \xrightarrow{\xi_{1}^{12}} A_{12} \xleftarrow{\xi_2^{12}} (f_2^{12})_*(\mathbb{A}_2).\]
    Our goal is to transport this pullback back up into the category $\cat{HTh}_m(G_1 \times_{G_{12}} G_2, \beta_1 \times_{\beta_{12}} \beta_2)$ and, as it turns out, the fact that the pushforward (i.e.~projection) maps between the categories of $m$-hypotheticals have right adjoints is crucial to be able to do this.  
}

\body{
    By using the right adjoint $(f_1^{12})^*$ of $(f_1^{12})_*$, we can lift the arrow $\xi_1^{12}$ to an arrow $(f_1^{12})^*(\xi_{1}^{12})$. However this arrow does not have $\mathbb{A}_1$ as its domain. This is where the fact that $(f_1^{12})^*$ is right adjoint comes handy: the unit of the adjunction, call it $\varepsilon_{1, \mathbb{A}_1}$, supplies the required map into the domain of $(f_1^{12})^*(\xi_{1}^{12})$. Diagrammatically we obtain the following composite morphism in the category $\cat{HTh}_m(G_1, \beta_1)$: 
    \[ \mathbb{A}_1 \xrightarrow{\varepsilon_{1, \mathbb{A}_1}} (f_1^{12})^*((f_1^{12})_*(\mathbb{A}_1)) \xrightarrow{(f_1^{12})^*(\xi_{1}^{12})} (f_1^{12})^*(G_{12}).\]
    Applying the very same reasoning, but to $(f_2^{12})_*$ and its right adjoint, we obtain the composite 
    \[ \mathbb{A}_2 \xrightarrow{\varepsilon_{2, \mathbb{A}_2}} (f_2^{12})^*((f_2^{12})_*(\mathbb{A}_2)) \xrightarrow{(f_2^{12})^*(\xi_{2}^{12})} (f_2^{12})^*(G_{12})\]
    in $\cat{HTh}_m(G_2, \beta_2)$. Both of these can be lifted to obtain the following diagram in $\cat{HTh}_m(G_1 \times_{G_{12}} G_2, \beta_1 \times_{\beta_{12}} \beta_2)$, using the right adjoints $\pi_1^*$ and $\pi_2^*$ of $(\pi_1)_*$ and $(\pi_2)_*$, respectively:
    \[\pi_1^*\mathbb{A}_1 \xrightarrow{\pi_1^*\left( (f_1^{12})^*(\xi_{1}^{12}) \circ \varepsilon_{1, \mathbb{A}_1} \right)} \pi_1^*(f_1^{12})^*(\mathbb{A}_{12}) \quad \quad \pi_2^*(f_2^{12})^*(\mathbb{A}_{12}) \xleftarrow{\pi_2^*\left( (f_2^{12})^*(\xi_{2}^{12}) \circ \varepsilon_{2, \mathbb{A}_2} \right)} \pi_2^*\mathbb{A}_2.\]
    This diagram forms a cospan in $\cat{HTh}_m(G_1 \times_{G_{12}} G_2, \beta_1 \times_{\beta_{12}} \beta_2)$ since $\pi_1^*(f_1^{12})^*(\mathbb{A}_{12})$ and $\pi_2^*(f_2^{12})^*(\mathbb{A}_{12})$ are isomorphic (to see this, notice that, since we started with a commuting square of pseudofunctors, their adjoints also form a commuting square of pseudofunctors).
}

\body{
    Now, since $\cat{HTh}_m(G_1 \times_{G_{12}} G_2, \beta_1 \times_{\beta_{12}} \beta_2)$ has pullbacks (c.f.~Lemma~\ref{lemma: fibers of HTh have pullbacks}), we can compute the pullback of
    \[\pi_1^*\mathbb{A}_1 \xrightarrow{\pi_1^*\left( (f_1^{12})^*(\xi_{1}^{12}) \circ \varepsilon_{1, \mathbb{A}_1} \right)} \pi_1^*(f_1^{12})^*(\mathbb{A}_{12}) \cong \pi_2^*(f_2^{12})^*(\mathbb{A}_{12}) \xleftarrow{\pi_2^*\left( (f_2^{12})^*(\xi_{2}^{12}) \circ \varepsilon_{2, \mathbb{A}_2} \right)} \pi_2^*\mathbb{A}_2\]
    to obtain a pullback object $P$ with pullback legs (i.e. pullback projections) $p_1$ and $p_2$. Theorem~\ref{thm: total category has pullbacks} ensures that the object \[(G_1 \times_{G_{12}} G_2, \beta_1 \times_{\beta_{12}} \beta_2, P)\] is the pullback of Diagram~\ref{eqn:  spelling out pullbacks in the total category}. The legs of its cone come once again from the adjunctions. By passing the leg $p_1$ through $(\pi_1)_*$ and by applying the counit of the adjunction, we obtain a morphism \[(\pi_1)_*(P) \xrightarrow{\eta_1} \mathbb{A}_1\]
    and hence we have that the leg of the cospan in the total category is $(\pi_1, \eta_1)$. Similarly we produce the other leg $(\pi_2, \eta_2)$. By similar reasoning and since the projections form a comutative square of pseudofunctors, we have the following comutative square in $\cat{HTh}_m(G_{12} , \beta_{12})$: 
    \[\begin{tikzcd}[row sep=small]
	& {(\pi_1)_*(f_1^{12})_* \mathbb{A}_1} && {(\pi_2)_*(f_2^{12})_* \mathbb{A}_2} & \\
	{(f_1^{12})_* \mathbb{A}_1} &&&& {(f_2^{12})_* \mathbb{A}_2} \\
	&& { A_{12}}
	\arrow["\sim", no head, from=1-2, to=1-4]
	\arrow[from=1-2, to=2-1]
	\arrow[from=1-4, to=2-5]
	\arrow["{\xi_{1}^{12}}"', from=2-1, to=3-3]
	\arrow["{\xi_2^{12}}", from=2-5, to=3-3]
\end{tikzcd}\]
    This defines a commuting square in the total category which, using the definition of arrows in the total category and the universal propertty of $P$, is easily seen to satisfy the universal property of pullbacks.
}

\subsection{A delicious analogue of Lemma~\ref{lem:lim_env_is_env_lim}}

\body{
    Let $f\colon G\to H$ be a homomorphism and $\beta\colon[\ell]\to\pts{G}$ an assignment. In Lemma~\ref{lem:helper_1}, we prove that $f$ is a morphism of hypotheses $\beta[G]\to (f\beta)[G]$. Now we prove an analogue of strategies and hypotheticals.
}

\begin{lemma}[Analogue of Lemma~\ref{lem:helper_1}]\label{lemma: canonical comparison from projected induced to induced projected}
    Let $f$ be a morphism and $S$ be any strategy over the position $(G, \beta)$, then there is a canonical morphism from the projection to $H$ of the induced hypothetical on $G$ of $S$ (namely $f_*(S[G])$) to the induced hypothetical on $H$ of the projection of $S$ (namely $f_*(S)[H])$. This morphism is denoted 
    \begin{equation}\label{eqn: canonical comparison from projected induced to induced projected}
        S[f] \colon f_*(S[G]) \to f_*(S)[H].
    \end{equation}
\end{lemma}
\begin{proof}
    Consider now any morphism $f \colon G \to H$. Using the adjunctions of Proposition~\ref{prop: induced hypothetical is right adjoint}, we can form the following picture. 
    \begin{equation}\label{eqn: comparing induced hypotheticals after projection}
        \begin{tikzcd}[sep=small]
	G && {\cat{HTh}_m(G, \beta)} && {\cat{Strat}_m(G, \beta)} \\
	\\
	H && {\cat{HTh}_m(H, f\beta)} && {\cat{Strat}_m(H, f\beta)}
	\arrow["f", from=1-1, to=3-1]
	\arrow[""{name=0, anchor=center, inner sep=0}, "U", shift left=2, from=1-3, to=1-5]
	\arrow["{f_*}", color={rgb,255:red,214;green,92;blue,92}, from=1-3, to=3-3]
	\arrow[""{name=1, anchor=center, inner sep=0}, "{\_[G]}", shift left=2, color={rgb,255:red,214;green,92;blue,92}, from=1-5, to=1-3]
	\arrow["{f_*}", color={rgb,255:red,92;green,92;blue,214}, from=1-5, to=3-5]
	\arrow[""{name=2, anchor=center, inner sep=0}, "U", shift left=2, from=3-3, to=3-5]
	\arrow[""{name=3, anchor=center, inner sep=0}, "{\_[H]}", shift left=2, color={rgb,255:red,92;green,92;blue,214}, from=3-5, to=3-3]
	\arrow["\dashv"{anchor=center, rotate=-90}, draw=none, from=0, to=1]
	\arrow["\dashv"{anchor=center, rotate=-90}, draw=none, from=2, to=3]
\end{tikzcd}
    \end{equation}
    This gives us a canonical way of comparing the image of an induced hypothetical on $G$ with the induced hypothetical on its projected strategy. In symbols, starting with a strategy $S \in \cat{Strat}_M(G, \beta)$ and staring at Diagram~\ref{eqn: comparing induced hypotheticals after projection}, we have two options: 
    \begin{enumerate}\vspace{-0.5em}
        \item we move left and then down, i.e.~we form an induced hypothetical on $G$ and then project it to $H$ (the {\color{red}red} path in Diagram~\ref{eqn: comparing induced hypotheticals after projection});\vspace{-0.5em}
        \item we move down and then left, i.e.~we project $S$ to a strategy on $H$ and then produce its induced hypothetical (the {\color{blue}blue} path in Diagram~\ref{eqn: comparing induced hypotheticals after projection}).\vspace{-0.5em}
    \end{enumerate}
    Since the underlying strategy functors $U$ commute (as pseudofunctors) with the pushforwards $f_*$ of $f$, we have that $U(f_*(S[G]))$ is isomorphic to $f_*(S)$ via an isomorphism $\iota \colon U(f_*(S[G])) \to f_*(S)$. The unit of the adjunction $U \dashv \_[H]$, gives us a canonical morphism $\eta_{f_*(S[G])} \colon f_*(S[G]) \to U(f_*(S[G]))[H]$, which, composed with the image of $\iota$ under $\_[H]$ gives us the desired comparison map of Equation~\ref{eqn: canonical comparison from projected induced to induced projected} as the following composite: \[ f_*(S[G]) \xrightarrow{\color{blue} \eta_{f_*(S[G])}} U(f_*(S[G]))[H] \xrightarrow{\color{red} \sim} f_*(S)[H].\]
    Diagrammatically, all of this can be summarized as follows. 
    \[\begin{tikzcd}[sep=small]
	& {S[G]} &&&& S \\
	G && {\cat{HTh}_m(G, \beta)} && {\cat{Strat}_m(G, \beta)} \\
	\\
	H && {\cat{HTh}_m(H, f\beta)} && {\cat{Strat}_m(H, f\beta)} \\
	& {f_*(S[G])} &&& {U(f_*(S[G]))} \\
	& {U(f_*(S[G]))[H]} &&&& {f_*(S)} \\
	&& {f_*(S)[H]}
	\arrow["\in"{description}, color={rgb,255:red,59;green,59;blue,59}, dotted, no head, from=1-2, to=2-3]
	\arrow[color={rgb,255:red,107;green,107;blue,107}, maps to, from=1-2, to=5-2]
	\arrow[color={rgb,255:red,107;green,107;blue,107}, maps to, from=1-6, to=1-2]
	\arrow["\ni"{description}, color={rgb,255:red,59;green,59;blue,59}, dotted, no head, from=1-6, to=2-5]
	\arrow[color={rgb,255:red,107;green,107;blue,107}, maps to, from=1-6, to=6-6]
	\arrow["f", from=2-1, to=4-1]
	\arrow[""{name=0, anchor=center, inner sep=0}, "U", shift left=2, from=2-3, to=2-5]
	\arrow["{f_*}", from=2-3, to=4-3]
	\arrow[""{name=1, anchor=center, inner sep=0}, "{\_[G]}", shift left=2, from=2-5, to=2-3]
	\arrow["{f_*}", from=2-5, to=4-5]
	\arrow[""{name=2, anchor=center, inner sep=0}, "U", shift left=2, from=4-3, to=4-5]
	\arrow[""{name=3, anchor=center, inner sep=0}, "{\_[H]}", shift left=2, from=4-5, to=4-3]
	\arrow["\in"{description}, color={rgb,255:red,59;green,59;blue,59}, dotted, no head, from=4-5, to=5-5]
	\arrow["\in"{description}, color={rgb,255:red,59;green,59;blue,59}, dotted, no head, from=4-5, to=6-6]
	\arrow["\in"{description, pos=0.7}, color={rgb,255:red,59;green,59;blue,59}, dotted, no head, from=5-2, to=4-3]
	\arrow[color={rgb,255:red,107;green,107;blue,107}, maps to, from=5-2, to=5-5]
	\arrow["\eta", color={rgb,255:red,92;green,92;blue,214}, from=5-2, to=6-2]
	\arrow[color={rgb,255:red,107;green,107;blue,107}, maps to, from=5-5, to=6-2]
	\arrow["\sim"{description}, color={rgb,255:red,214;green,92;blue,92}, from=5-5, to=6-6]
	\arrow["\in"{description, pos=0.8}, color={rgb,255:red,59;green,59;blue,59}, dotted, no head, from=6-2, to=4-3]
	\arrow["\sim"{description}, color={rgb,255:red,214;green,92;blue,92}, from=6-2, to=7-3]
	\arrow[color={rgb,255:red,107;green,107;blue,107}, maps to, from=6-6, to=7-3]
	\arrow["\in"{description, pos=0.9}, color={rgb,255:red,59;green,59;blue,59}, dotted, no head, from=7-3, to=4-3]
	\arrow["\dashv"{anchor=center, rotate=-90}, draw=none, from=0, to=1]
	\arrow["\dashv"{anchor=center, rotate=-90}, draw=none, from=2, to=3]
\end{tikzcd}\]

\end{proof}
    
\begin{theorem}\label{thm: limit of strategies is strategy of limit}
    Let $\Ta=(T,\{G_a\},\{g_t^e\})$ be a tree co-decomposition and for each $a\in V(T)\cup E(T)$ denote by $\pi_a$ the projection $\lim\Ta\to G_a$. Let\vspace{-0.5em}
    \begin{itemize}
        \item $\beta \colon [2\ell - m] \to \lim\Ta$; \vspace{-0.5em}
        \item $S\in\cat{Strat}_m(\lim\Ta,\beta)$ be am $m$-strategy on $\lim\Ta$ with respect to $\beta$; and\vspace{-0.5em}
        \item for each $a\in V(T)\cup E(T)$, let ${\pi_a}_*S \in\cat{Strat}_m(\lim\Ta,\pi_a \circ \beta)$ be the projection of $S$ to $(G_a, \pi_a \circ \beta)$.\vspace{-0.5em}
    \end{itemize}
    For each $t\in V(T)$ and incident $e\in E(T)$, consider the induced morphism of strategies (c.f.~Lemma~\ref{prop:functoriality-of-projection}) in the total category\vspace{-0.5em}
    \[
        (g_t^e, g_{t*}^e) \colon \left(G_t,\;\pi_t\beta,\;\pi_{t*} S \right)\to \left(  G_e,\; \pi_e \beta,\;\pi_{e*} S\right), 
    \]
    then the tuple\vspace{-0.5em}
    \[
    \Ta_{*}S := \; \left(T,\;\underset{a\in V(T)\cup E(T)}{\{\pi_{a*}S\}},\;\underset{t\in e\in E(T)}{\{(g_t^e, g_{t*}^e)\}}\right)
    \]
    is a co-decomposition valued in the total category of strategies. Furthermore, there is a canonical map\vspace{-0.5em}
    \begin{equation}\label{eqn: canonical map of limit of strategies is strategy of limit}
        \left(\lim \Ta, \beta, S \right) \to \lim (\Ta_{*}S).
    \end{equation}
\end{theorem}
\begin{proof}
    Observe that $\Ta_{*}S$ is a codecompostion by definition. Recalling that the total category of strategies has pullbacks (c.f.~Theorem~\ref{thm: total category has pullbacks}) and noticing that $\left(\lim \Ta, \beta, S \right)$ forms the apex of a cone sitting over $\Ta_{*}S$, we get the desired morphism~\ref{eqn: canonical map of limit of strategies is strategy of limit} as the unique limit arrow (recall that we compute the limit of $\Ta_{*}S$ recursively by successive pullbacks).
\end{proof}

\begin{corollary}[Analogue for Lemma~\ref{lem:lim_env_is_env_lim}]\label{corollary: deliscious}
    Let $\Ta=(T,\{G_a\},\{g_t^e\})$ be a tree co-decomposition and for each $a$ denote by $\pi_a$ the projection $\lim\Ta\to G_a$. Let $S$ be a strategy for $\lim\Ta$.

    For each $a\in V(T)\cup E(T)$, consider the hypothetical $\mathbb{H}_a:=({\pi_a}_*S)[G_a]$ induced by the projected strategy ${\pi_a}_*S$. For each $t\in V(T)$ and incident $e\in E(T)$, consider the induced morphism of induced hypotheticals (Lemma~\ref{lemma: canonical comparison from projected induced to induced projected}) $\pi_{t*}(S)[g_t^e]\colon{{g_t^e}_*(\;\pi_{t*}(S)[G_t]\;)} \to  {({g_t^e}_*\circ\pi_{t*}(S))[G_e]}$ and let $\rho_t^e\colon \mathbb{H}_t\to\mathbb{H}_e$ be the composite defined as follows 
    \[\pi_{t*}(S)[G_t] 
        \xrightarrow{ {g_t^e}_* } 
            {g_t^e}_*(\;\pi_{t*}(S)[G_t]) 
                \xrightarrow{\pi_{t*}(S)[g_t^e] } 
                    ({g_t^e}_*\circ\pi_{t*}(S))[G_e]                
                        \xrightarrow{\sim} 
                        \pi_{e*}(S)[G_e].
    \]
    The tuple $S[\Ta]:=(T,\{\mathbb{H}_a\},\{\rho_e^t\})$ is a co-decomposition of hypotheticals. Furthermore, there is a canonical map $S[\lim\Ta] \to \lim S[\Ta]$.
\end{corollary}
\begin{proof}
    Notice that $S[\Ta]$ is none other than the codecompostion obtained by passing all the bags and adhesion maps of the codecompostion $\Ta_* S$ of Theorem~\ref{thm: limit of strategies is strategy of limit} through the adjunction $U_{G_a} \dashv \; \_[G_a]$ of Proposition~\ref{prop: induced hypothetical is right adjoint}. By this very adjunction, since induced hypotheticals are right adjoint to the underlying strategy functor and since right adjoints preserve limits, there is a canonical map from $S[\lim\Ta] \to \lim S[\Ta]$. 
\end{proof}

\begin{addendum}[Addendum to the proof of Corollary~\ref{corollary: deliscious}]
    Those who prefer reasoning purely via abstract nonsense can safely skip this addendum. However, there may be readers who might find the above proof a little psychologically unsatisfying since it is nontrivial to see that all the required morphisms do indeed exists and that they are in the right places. For such readers, we will spell things out a little more below. To that end consider 
    
\[\begin{tikzcd}[sep=small]
	& {G_1} &&&& {(G_1, \pi_1\beta, \pi_{1*}S)} \\
	\\
	&&&&& {(G_{12}, g_1^{12} \circ \pi_{1} \circ \beta, (g_{1}^{12})_*(\pi_{{1}*}(S)))} \\
	G & {G_{12}} && {(G, \beta, S)} & {(G, \beta, P)} & {(G_{12}, \pi_{12} \circ \beta, \pi_{{12}*}(S))} \\
	&&&&& {(G_{12}, g_2^{12} \circ \pi_{2} \circ \beta, (g_{2}^{12})_*(\pi_{{2}*}(S)))} \\
	\\
	& {G_1} &&&& {(G_2, \pi_2\beta, \pi_{2*}S)}
	\arrow["{g_1^{12}}", from=1-2, to=4-2]
	\arrow["{(g_1^{12})_*}", from=1-6, to=3-6]
	\arrow["\sim", no head, from=3-6, to=4-6]
	\arrow["{\pi_1}", curve={height=-6pt}, from=4-1, to=1-2]
	\arrow["{\pi_2}"', curve={height=6pt}, from=4-1, to=7-2]
	\arrow[curve={height=-30pt}, from=4-4, to=1-6]
	\arrow[curve={height=30pt}, from=4-4, to=7-6]
	\arrow[curve={height=-12pt}, from=4-5, to=1-6]
	\arrow[curve={height=12pt}, from=4-5, to=7-6]
	\arrow["\sim", no head, from=5-6, to=4-6]
	\arrow["{g_2^{12}}"', from=7-2, to=4-2]
	\arrow[from=7-6, to=5-6]
\end{tikzcd}\]
where $P$ is the pullback of strategies $\pi_{1*}S \times_{\pi_{12*} S} \pi_{2*}S$. Applying the `induced hypothetical' functor ${}\_[G]$ gives us
\[\begin{tikzcd}[sep=small]
	& {(G_{1}, \pi_{1} \circ \beta, \;\pi_{{1}*}S[G_1])} & {(G_1, \pi_1\beta,   \;  (\pi_{1*}S)[G_1])} \\
	\\
	& {(G_{12}, g_1^{12} \circ \pi_{1} \circ \beta, \; \left((g_{1}^{12})_*\pi_{{1}*}S\right)[G_1]))} & {(G_{12}, g_1^{12} \circ \pi_{1} \circ \beta,  \;  \left((g_{1}^{12})_*\pi_{{1}*}S\right)[G_{12}])} \\
	{(G, \beta, S[G])} && {(G_{12},  \pi_{12} \circ \beta,  \;  \left(\pi_{{12}*}S \right)[G_{12}])} \\
	& {(G_{12}, g_2^{12} \circ \pi_{2} \circ \beta,   \;   \left((g_{2}^{12})_*\pi_{{2}*}S\right)[G_2]))} & {(G_{12}, g_2^{12} \circ \pi_{2} \circ \beta,  \;  \left((g_{2}^{12})_*\pi_{{2}*}S\right)[G_{12}])} \\
	\\
	& {(G_{1}, \pi_{1} \circ \beta, \; \pi_{{1}*}S[G_1])} & {(G_2, \pi_2\beta,   \;  (\pi_{2*}S)[G_2])}
	\arrow["{S[\pi_1]}", from=1-2, to=1-3]
	\arrow["{(g_1^{12})_*}"{description}, from=1-3, to=3-2]
	\arrow["{(\pi_{1*}S)[g_1^{12}] := \eta_{\pi_{1*}S} \circ \cong}"', shift right=3, from=3-2, to=3-3]
	\arrow["\sim", no head, from=3-3, to=4-3]
	\arrow[curve={height=-30pt}, from=4-1, to=1-2]
	\arrow[curve={height=30pt}, from=4-1, to=7-2]
	\arrow["{(\pi_{2*}S)[g_2^{12}] := \eta_{\pi_{2*}S} \circ \cong}", shift left=3, from=5-2, to=5-3]
	\arrow["\sim"', no head, from=5-3, to=4-3]
	\arrow["{S[\pi_2]}"', from=7-2, to=7-3]
	\arrow["{(g_2^{12})_*}"{description}, from=7-3, to=5-2]
\end{tikzcd}\]
from which we deduce, by preservation of limits (since ${}\_[G]$ is a right-adjoint), that
    \begin{align*}
        (G, \beta, P[G]) &:= \left( (G_1, \; \pi_1\beta, \; \pi_{1*}S) \times_{(G_{12}, \; \pi_{12} \circ \beta, \; \pi_{{12}*}(S))} (G_{12}, \; \pi_{12} \circ \beta, \; \pi_{{12}*}(S)) \right)[G_{12}] \\
        &\cong (G_1, \; \pi_1\beta, \;\;(\pi_{1*}S)[G_1]) \times_{(G_{12}, \;  \pi_{12} \circ \beta,  \;\;  \left(\pi_{{12}*}S \right)[G_{12}])} (G_2, \; \pi_2\beta,   \;\;  (\pi_{2*}S)[G_2]).
    \end{align*}
\end{addendum}

\body{The distinction between projecting an induced hypothetical and inducing a projected strategy is subtle and, perhaps, depending on the reader's psychology, easier to see through abstraction than combinatorially. Below, we attempt to illustrate it concretely with a small example on a ``cross graph'' (four points, two diagonal edges; i.e.~$2K_2$). The example also also shows why the additional outer family in the definition of hypotheticals is necessary when computing pullbacks.
}

\begin{example}
    Let $A$ and $B$ each be copies of the two-point structure with a single symmetric edge $0\sim1$, and let $o$ be the one-point structure with a loop. Then $X:=A\times_o B$ is the graph \[\begin{tikzcd}[row sep=small]
	00 & 01 \\
	10 & 11
	\arrow[no head, from=1-1, to=2-2]
	\arrow[no head, from=1-2, to=2-1]
\end{tikzcd}\] with has points $\{00,01,10,11\}$, and edges between points at Hamming distance $2$ (both coordinates are different). Consider the strategy $S$ with universal/existential play/response given by
    $ 
                00\mapsto01,\;
                01\mapsto11,\;
                11\mapsto10,\;
                10\mapsto00,
    $
    visualised below as a clockwise turn around the square:
    \[
    \begin{tikzcd}[row sep=1.5em, column sep=2.5em]
        00 \arrow[r] & 01 \arrow[d] \\
        10 \arrow[u] & 11 \arrow[l]
    \end{tikzcd}
    \]
    In this example, hypotheses can only record whether there is an edge or a non-edge, so, for brevity, write $\top$ for the induced two-point hypothesis in which the selected pair is related, and $\bot$ for the one in which it is not. Since each $x$ and $S(x)$ have Hamming distance $1$, the induced hypothetical will always contain $\bot$ as the singleton hypothesis; in particular, the induced hypothetical is
    \[
        S[X](x)=\bigl[[\,x\cdot S(x)\mapsto\{\bot\}\,]\bigr].
    \]

    Let $\pi_A\colon X\to A$, $\pi_B\colon X\to B$, and $\pi_o\colon X\to o$ be the projections. We write an assignment $[2]\to A$ sending $1\mapsto a_1$ and $2\mapsto a_2$ as $a_1\cdot a_2$. With this notation, the projection of $S[X]$ to $A$ is given below (labels on the underbraces indicate which universal choice in $X$ the corresponding personality was projected from: we add these for legibility):
    \[
    \begin{aligned}
    (\pi_{A*}S[X])(0)
    &=
    \left[
    \underbrace{[\,0\cdot0\mapsto\{\bot\}\,]}_{00},
    \underbrace{[\,0\cdot1\mapsto\{\bot\}\,]}_{01}
    \right],\\
    (\pi_{A*}S[X])(1)
    &=
    \left[
    \underbrace{[\,1\cdot0\mapsto\{\bot\}\,]}_{10},
    \underbrace{[\,1\cdot1\mapsto\{\bot\}\,]}_{11}
    \right].
    \end{aligned}
    \]
    We read this as follows: the assignment $0\cdot1$, for example, is carrying the hypothesis $\bot$, meaning the no-edge hypothesis. Notice that this hypothesis need not be induced here (indeed, in this case it isn't! --- $0$ is adjacent to $1$ in $A$): it is simply the projection of what we had in $S[X]$ at that coordinate.

    Before moving on and having a look at the other projections, it is worth contrasting this with first projecting the strategy $S$ and then taking its induced hypothetical on $A$. In this case we obtain
    \[
    \begin{aligned}
    ((\pi_{A*}S)[A])(0)
    &=
    \left[
        [\,0\cdot0\mapsto\{\bot\},\;
           0\cdot1\mapsto\{\top\}\,]
    \right],\\
    ((\pi_{A*}S)[A])(1)
    &=
    \left[
        [\,1\cdot0\mapsto\{\top\},\;
           1\cdot1\mapsto\{\bot\}\,]
    \right].
    \end{aligned}
    \]
    Thus projecting the induced hypothetical and inducing the projected strategy are not the same operation, that's why we need the canonical comparison map of Lemma~\ref{lemma: canonical comparison from  projected induced to induced projected}. In this example, when   focusing on $A$, the map is $S[\pi_A]\colon \pi_{A*}(S[X]) \to (\pi_{A*}S)[A]$ which merges the two projected personalities over each universal choice and sends their carried hypotheses to the corresponding induced hypotheses in $A$.

    Returning to the projections, shipping $S[X]$ to $B$, produces 
    \[
    \begin{aligned}
    (\pi_{B*}S[X])(0)
    &=
    \left[
    \underbrace{[\,0\cdot1\mapsto\{\bot\}\,]}_{00},
    \underbrace{[\,0\cdot0\mapsto\{\bot\}\,]}_{10}
    \right],\\
    (\pi_{B*}S[X])(1)
    &=
    \left[
    \underbrace{[\,1\cdot1\mapsto\{\bot\}\,]}_{01},
    \underbrace{[\,1\cdot0\mapsto\{\bot\}\,]}_{11}
    \right].
    \end{aligned}
    \]

    \textbf{N.B.} It is important to notice that, whether we project to $A$ or to $B$, each projection has two distinct personalities over each universal choice. The point of the outer family in \cref{def: hypothetical} is to keep these behaviors separate rather than merging them into a single nondeterministic response. In contrast, projected strategies lose this distinction, since their definition involves only a single application of $\cat{Fam}$.

    Projecting further to $o$ gives four personalities, which, for the same reason we just pointed out, despite having identical underlying behavior, remain distinct entries of the family. This projection is:
    \[
    (\pi_{o*}S[X])(\bullet)
    =
    \left[
    \underbrace{[\,\bullet\cdot\bullet\mapsto\{\bot\}\,]}_{00},
    \underbrace{[\,\bullet\cdot\bullet\mapsto\{\bot\}\,]}_{01},
    \underbrace{[\,\bullet\cdot\bullet\mapsto\{\bot\}\,]}_{10},
    \underbrace{[\,\bullet\cdot\bullet\mapsto\{\bot\}\,]}_{11}
    \right].
    \]
    In the reduced version of the above, the personalities merge into one. It is the maps between hypotheticals that `remember' that there should be distinctions: In the cospan $\pi_{A*}S[X]\to\pi_{o*}S[X] \leftarrow \pi_{B*}S[X],$ the personality arising from $x$ on the $A$-side and the personality arising from $x$ on the $B$-side are both sent to the personality labelled by $x$ over $o$. 

    Consequently, the pullback matches precisely the two personalities arising from the same $x\in X$. Their responses pull back to
    \[
        \bigl(\pi_A(S(x)),\pi_B(S(x))\bigr)=S(x),
    \]
    and their hypotheses pull back to $x\cdot S(x)\mapsto\{\bot\}$. Hence the pullback recovers exactly the original hypothetical $S[X]$.

    The tricky part, which is worth being explicit about, is that, although the two projections to $A$ and $B$ individually forget the ``direction'' of the spiral, it is the cospan that they form over $o$ that ``remembers'' how to pair them up appropriately so that the pullback reconstructs the original spiral. The comparison maps $S[\pi_A]$, $S[\pi_B]$, and $S[\pi_o]$ then relate these projected induced hypotheticals to the induced hypotheticals of the projected strategies appearing above.
\end{example}

\subsection{Proof of Theorem~\ref{thm: FO model checking on codecompositions}: FO-Model-Checking in FPT-Time on Co-Decompositions}

\body{
    At this point, the dynamic programming argument for full FO follows exactly the same pattern as the one given for $\exists$-FO. This is the payoff of going through the trouble of proving the existence of pullbacks of hypotheicals: we took care in finding canonical definitions of the objects used in the dynamic programming so that we could reveal the satisfying (to us) fact that, structurally, from an algorithmic point of view, there is no distinction between $\exists$-FO and full FO (despite the latter being much more complicated). 
}

\body{
    The arguments we employed to establish the algorithm for $\exists$-FO used \textit{no special property of assignments and hypotheses} except for:\vspace{-0.5em}
    \begin{enumerate}
        \item the existence of an induced-object construction,\vspace{-0.5em}
        \item functorial projections to the bags,\vspace{-0.5em}
        \item pullbacks for combining compatible local states,\vspace{-0.5em}
        \item the fact that the collection of states on each bag has size bounded by the parameter.\vspace{-0.5em}
    \end{enumerate}
    All the machinery we built so far established the corresponding ingredients for strategies and hypotheticals. Consequently, once we define supported hypotheticals and write down the corresponding supportedness recurrence, the algorithm of Theorem~\ref{thm:exists-fo} can be lifted to the full FO model checking case by employing the dictionary we gave earlier in Table~\ref{tab:existential-to-full-fo}.
}

\body{
    Now, we will fill-in the only two missing points: supported hypotheticals (third from bottom in Table~\ref{tab:existential-to-full-fo}) and their corresponding ``head condition'' (last in Table~\ref{tab:existential-to-full-fo}). We will deliberately follow the structure of Paragraph~\ref{para: def: supported hypoethis} so that the analogy is unmistakable. Let $\Ta=(T,\{G_a\},\{g_t^e\})$ be a rooted co-decomposition, $i$ be a node threin, and $\pi_i\colon G[i]\rightarrow G_i$ be the canonical projection from the partial limit at $i$ to its root bag. We say that an $m$-hypothetical $\mathbb{H}_i\in\cat{HTh}_m(G_i,\beta_i)$ is \define{supported at $i$} if $\mathbb{H}_i=(\pi_i)_*(S[G[i]])$ for some $m$-strategy $S\in\cat{Strat}_m(G[i],\beta)$, where $\beta$ projects to $\beta_i$. Given an $m$-strategy $S_i\in\cat{Strat}_m(G_i,\beta_i)$, the \define{supported hypotheticals for $S_i$} are those supported hypotheticals whose underlying strategy is $S_i$. Again, as was the case for $\exists$-FO, what is important is which hypotheticals are supported at the root. This is the following lemma.
}

\begin{lemma}[Analogue of Lemma~\ref{lemma: support at root is enough}]\label{lemma: satisfiability at root main}
    Take an FO-formula $\phi = \forall x_1 \exists y_1, \dots, \forall x_\ell \exists y_\ell \psi(\vec{x}, \vec{y})$ in alternating prenex form over a signature $\sigma$ and consider a tree co-decomposition $\Ta$ rooted at some node~$r$. Then $\lim \Ta \models \phi$ if and only if there is a supported $\ell$-hypothetical $\mathbb{H}$ ar $r$ such that the head of every hypothesis occurring in a $0$-hypothetical of $\mathbb{H}$ satisfies $\phi$.   
\end{lemma}
\begin{proof}
    Thinking of satisfaction as a game between a universal player and an existential player, $\lim \Ta \models \phi$ if and only if there is a there is an $\ell$-strategy $S$ on $\lim \Ta$ such that every hypothesis in the hypothetical $S[\lim \Ta]$ induced on $\lim \Ta$ has head satisfying $\phi$. Since $G[r] \cong \lim \Ta$, the supported hypotheticals at the root $r$ are precisely the projections of these induced hypotheticals to the bag $G_r$. Since projection doesn't change the heads of hypotheses, the result follows.  
\end{proof}

\begin{observation}
    \label{prop: computing supported hypotheticals}
    The dynamic programming recurrences of Lemmas~\ref{lem:running_env} and~\ref{lem:two_children} remain valid for supported hypotheticals under the dictionary of Table~\ref{tab:existential-to-full-fo}. In particular, the supported hypotheticals at every node can be computed bottom-up using only projections, the canonical comparison maps $S[f]\colon f_*(S[G])\longrightarrow f_*(S)[H]$, and pullbacks in the total category of hypotheticals.
\end{observation}

\body{
    The proof of the above observation is the same as in the existential case, so we just sketch it. To see that it does indeed hold, notice that at a leaf, supported hypotheticals are simply induced hypotheticals on the bag. At an internal node, one combines the supported hypotheticals coming from the relevant subtree or subtrees with compatible local strategies on the current bag, using the canonical comparison maps $S[f]$ to situate the resulting induced hypotheticals over the appropriate adhesion and then taking their pullback in $\cat{HTh}_m^\sigma$. (This is exactly what we did in the case of hypotheses.) The universal property of these pullbacks gives exactly the hypotheticals arising from strategies on the corresponding partial limit. Thus the bottom-up recurrence of \cref{lem:running_env,lem:two_children} carries over verbatim.
}

\begin{lemma}\label{lem: bounded hypothetical state space}
    For fixed signature $\sigma$, formula $\phi$, and width $w$, the number of reduced hypotheticals that can occur over a bag of size at most $w$ is bounded by a function of $w$ and $|\phi|$.
\end{lemma}

\body{
    \cref{lem: bounded hypothetical state space} follows directly from the recursive definition using \textit{reduced} families. Perhaps more interesting is to see explicitly just how horrible the explicit value of the above bound can get. At each internal node, we may need to consider every pair of supported hypotheticals arising from the relevant child tables. Thus, if $h_\ell$ bounds the number of reduced $\ell$-hypotheticals stored at any bag, then the number of candidate pairs considered at a node is at most $h_\ell^2.$ If $\mathsf{compPB}(w,|\phi|)$ denotes an upper bound on the time required to compute and reduce a single pullback of hypotheticals over bags of width at most $w$, then the work performed at each node is bounded by $O\!\left(h_\ell^2\,\mathsf{compPB}c_{\mathrm{pb}}(w,|\phi|)\right)$. Thus, over a co-decomposition with $b$ bags, the total running time is bounded by $O\!\left( b\,h_\ell^2\,c_{\mathrm{pb}}(w,|\phi|) \right)$. Of course, the worst term is $h_\ell$ which, staring at the recursive definition of hypotheticals, amounts to: 
\[
    h_\ell
    \leq
    \underbrace{
    2^{\,w\,2^{\,w\,2^{\,w\,2^{\cdots^{\,w\,2^{\,w\,2^{2^q}}}}}}}
    }_{2\ell+2\text{ exponential layers}},
    \quad \text{where} \quad 
    q:=\sum_{R\in\sigma}(2\ell)^{\operatorname{ar}(R)}.
\]
}\label{para: explicit bounds}

\begin{proof}[Proof of Theorem~\ref{thm: FO model checking on codecompositions}]
    Put $\phi$ into alternating prenex normal form and process the co-decomposition from the leaves towards the root, maintaining at each node the reduced supported hypotheticals. By Proposition~\ref{prop: computing supported hypotheticals}, these tables are computed by exactly the same bottom-up procedure as the supported-hypothesis tables of Theorem~\ref{thm:exists-fo}, replacing hypotheses by hypotheticals and taking the corresponding pullbacks in $\cat{HTh}_m^\sigma$.

    By Lemma~\ref{lem: bounded hypothetical state space}, the number and size of the reduced hypotheticals considered at any one node are bounded by a function of $w$ and $|\phi|$. Since each node is processed once, the total running time is $f_\sigma(w,|\phi|)\cdot b.$ Finally, at the root we apply Lemma~\ref{lemma: satisfiability at root main} and accept when there is a supported hypothetical all of whose leaf hypotheses satisfy the quantifier-free part of $\phi$.
\end{proof}

%% file: figures/grothendieck_construction.tex
\begin{tikzpicture}[
    >=Latex,
    every node/.style={font=\small},
    basebox/.style={draw, rounded corners, inner sep=12pt},
    fiberbox/.style={draw, rounded corners, inner sep=12pt},
    totalbox/.style={draw, rounded corners, inner sep=20pt}
]

\node (I) at (0,0) {$[2(\ell-m)]$};
\node (G) at (2.8,1.2) {$G$};
\node (H) at (2.8,-1.2) {$H$};

\draw[->] (I) -- node[above left] {$\beta$} (G);
\draw[->] (I) -- node[below left] {$\gamma$} (H);
\draw[->] (G) -- node[right] {\color{red} $f$} (H);

\node[
    basebox,
    fit=(I)(G)(H),
    label={[anchor=south west]north west:{$\cat{Pos}_m^\sigma$}}
] (PosBox) {};

\node (A) at (8.7,1.7) {$\mathbb{A}$};

\node[
    fiberbox,
    fit=(A),
    minimum width=5.8cm,
    minimum height=2.2cm,
    label={[anchor=south east, xshift=-2pt, yshift=2pt]north east:{$\cat{HTh}_m(G,\beta)$}}
] (TopBox) {};

\node (fA) at (7.5,-2.0) {$f_*\mathbb{A}$};
\node (B)  at (10.0,-2.0) {$\mathbb{B}$};

\draw[->] (fA) -- node[above] {\color{red}$\eta$} (B);

\node[
    fiberbox,
    fit=(fA)(B),
    minimum width=7.0cm,
    minimum height=2.4cm,
    label={[anchor=south east, xshift=-2pt, yshift=2pt]north east:{$\cat{HTh}_m(H,\gamma)$}}
] (BotBox) {};

\draw[->] (TopBox.south) -- node[right] {$f_*$} (BotBox.north);

\node[
    totalbox,
    fit=(TopBox)(BotBox),
    inner sep=18pt,
    label={[anchor=south east, xshift=-2pt, yshift=2pt]north east:{$\cat{HTh}_m^\sigma$}}
] (TotBox) {};

\end{tikzpicture}

%% file: figures/pullback_in_total_category.tex
\begin{tikzpicture}[
    >=Latex,
    font=\small,
    posbox/.style={draw, rounded corners, inner sep=10pt},
    fibbox/.style={draw, rounded corners, inner sep=10pt},
    totbox/.style={draw, rounded corners, inner sep=16pt},
    obj/.style={inner sep=1pt},
    labNE/.style={anchor=south east, xshift=-2pt, yshift=2pt},
    labNW/.style={anchor=south west, xshift=2pt, yshift=2pt}
]


\node[
    fibbox,
    minimum width=14.0cm,
    minimum height=2.1cm,
    label={[labNE, name=FibPBLabel]north east:{
        $\cat{HTh}_m
        (G_1\times_{G_{12}}G_2,
        \beta_1\times_{\beta_{12}}\beta_2)$
    }}
] (FibPB) at (0,6.65) {};

\node[obj] (pi1A1) at (-6.1,6.25)
    {$\pi_1^*\mathbb A_1$};

\node[obj] (pi1A12) at (-1.45,6.25)
    {$\pi_1^*(f_1^{12})^*\mathbb A_{12}$};

\node[obj] (pi2A12) at (1.45,6.25)
    {$\pi_2^*(f_2^{12})^*\mathbb A_{12}$};

\node[obj] (pi2A2) at (6.1,6.25)
    {$\pi_2^*\mathbb A_2$};

\node at (0,6.25) {$\cong$};

\draw[->]
    (pi1A1) --
    node[above, inner sep=2pt] {
        $\begin{gathered} 
            \pi_1^*\!\bigl((f_1^{12})^*(\zeta_1^{12})\circ{}
            \varepsilon_{1,\mathbb A_1}\bigr)
        \end{gathered}$
    }
    (pi1A12);

\draw[->]
    (pi2A2) --
    node[above, inner sep=2pt] {
        $\begin{gathered}
            \pi_2^*\!\bigl((f_2^{12})^*(\zeta_2^{12})\circ{}
            \varepsilon_{2,\mathbb A_2}\bigr)
        \end{gathered}$
    }
    (pi2A12);


\node[
    fibbox,
    minimum width=6.4cm,
    minimum height=1.4cm,
    label={[labNW, name=Fib1Label]north west:
        {$\cat{HTh}_m(G_1,\beta_1)$}}
] (Fib1) at (-3.65,3.2) {};

\node[obj] (A1) at (-6.3,3.0)
    {$\mathbb A_1$};

\node[obj] (unitTarget1) at (-4.1,3.0)
    {$(f_1^{12})^*f_{1*}^{12}\mathbb A_1$};

\node[obj] (target1) at (-1.5,3.0)
    {$(f_1^{12})^*\mathbb A_{12}$};

\draw[->]
    (A1) --
    node[above] {$\varepsilon_{1,\mathbb A_1}$}
    (unitTarget1);

\draw[->]
    (unitTarget1) --
    node[above] {$(f_1^{12})^*(\zeta_1^{12})$}
    (target1);

\node[
    fibbox,
    minimum width=6.4cm,
    minimum height=1.4cm,
    label={[labNE, name=Fib2Label]north east:
        {$\cat{HTh}_m(G_2,\beta_2)$}}
] (Fib2) at (3.65,3.2) {};

\node[obj] (A2) at (1.0,3.0)
    {$\mathbb A_2$};

\node[obj] (unitTarget2) at (3.2,3.0)
    {$(f_2^{12})^*f_{2*}^{12}\mathbb A_2$};

\node[obj] (target2) at (5.8,3.0)
    {$(f_2^{12})^*\mathbb A_{12}$};

\draw[->]
    (A2) --
    node[above] {$\varepsilon_{2,\mathbb A_2}$}
    (unitTarget2);

\draw[->]
    (unitTarget2) --
    node[above] {$(f_2^{12})^*(\zeta_2^{12})$}
    (target2);


\node[
    fibbox,
    minimum width=8.6cm,
    minimum height=2.2cm,
    label={[labNE]south east:
        {$\cat{HTh}_m(G_{12},\beta_{12})$}}
] (Fib12) at (0,-0.8) {};

\node[obj] (p1A1) at (-3.0,-0.7)
    {$f_{1*}^{12}\mathbb A_1$};

\node[obj] (A12) at (0,-0.7)
    {$\mathbb A_{12}$};

\node[obj] (p2A2) at (3.0,-0.7)
    {$f_{2*}^{12}\mathbb A_2$};

\draw[->]
    (p1A1) --
    node[above] {$\zeta_1^{12}$}
    (A12);

\draw[->]
    (p2A2) --
    node[above] {$\zeta_2^{12}$}
    (A12);


\draw[->]
    ([xshift=-2.1cm]FibPB.south) --
    node[above left] {$\pi_{1*}$}
    coordinate[midway] (pi1Down)
    ([xshift=-0.3cm]Fib1.north);

\draw[->]
    ([xshift=0.3cm]Fib1.north) --
    node[below right] {$\pi_1^*$}
    coordinate[midway] (pi1Up)
    ([xshift=-1.5cm]FibPB.south);

\node at ($(pi1Down)!0.5!(pi1Up)$) {$\dashv$};

\draw[->]
    ([xshift=1.5cm]FibPB.south) --
    node[below left] {$\pi_{2*}$}
    coordinate[midway] (pi2Down)
    ([xshift=-0.3cm]Fib2.north);

\draw[->]
    ([xshift=0.3cm]Fib2.north) --
    node[above right] {$\pi_2^*$}
    coordinate[midway] (pi2Up)
    ([xshift=2.1cm]FibPB.south);

\node at ($(pi2Down)!0.5!(pi2Up)$) {$\dashv$};

\draw[->]
    ([xshift=-0.3cm]Fib1.south) --
    node[below left] {$f_{1*}^{12}$}
    coordinate[midway] (p1Down)
    ([xshift=-2.1cm]Fib12.north);

\draw[->]
    ([xshift=-1.5cm]Fib12.north) --
    node[above right] {$(f_1^{12})^*$}
    coordinate[midway] (p1Up)
    ([xshift=0.3cm]Fib1.south);

\node at ($(p1Down)!0.5!(p1Up)$) {$\dashv$};

\draw[->]
    ([xshift=-0.3cm]Fib2.south) --
    node[above left] {$f_{2*}^{12}$}
    coordinate[midway] (p2Down)
    ([xshift=1.5cm]Fib12.north);

\draw[->]
    ([xshift=2.1cm]Fib12.north) --
    node[below right] {$(f_2^{12})^*$}
    coordinate[midway] (p2Up)
    ([xshift=0.3cm]Fib2.south);

\node at ($(p2Down)!0.5!(p2Up)$) {$\dashv$};


\node[
    totbox,
    fit=(FibPB)(FibPBLabel)
        (Fib1)(Fib1Label)
        (Fib2)(Fib2Label)
        (Fib12),
    label={[labNE]south east:{$\cat{HTh}_m^\sigma$}}
] (TotBox) {};


\begin{scope}[yshift=-5.5cm]

\node[obj] (I) at (0,1.8)
    {$[2(\ell-m)]$};

\node[obj] (PB) at (0,0)
    {$G_1\times_{G_{12}}G_2$};

\node[obj] (G1) at (-3.6,-1.9) {$G_1$};
\node[obj] (G2) at (3.6,-1.9) {$G_2$};
\node[obj] (G12) at (0,-4.0) {$G_{12}$};

\draw[->]
    (I) --
    node[pos=0.8,right]
        {$\beta_1\times_{\beta_{12}}\beta_2$}
    (PB);

\draw[->]
    (I) to[bend right=18]
    node[above left] {$\beta_1$}
    (G1);

\draw[->]
    (I) to[bend left=18]
    node[above right] {$\beta_2$}
    (G2);

\draw[->]
    (PB) --
    node[above left] {$\pi_1$}
    (G1);

\draw[->]
    (PB) --
    node[above right] {$\pi_2$}
    (G2);

\draw[->]
    (G1) --
    node[below left] {$f_1^{12}$}
    (G12);

\draw[->]
    (G2) --
    node[below right] {$f_2^{12}$}
    (G12);

\node[
    posbox,
    fit=(I)(PB)(G1)(G2)(G12),
    inner sep=14pt,
    label={[labNE]south east:{$\cat{Pos}_m^\sigma$}}
] (PosBox) {};

\end{scope}

\end{tikzpicture}

%% file: sections/discussion.tex
\section{Discussion \& Open Problems}\label{sec: discussion}

\body{Here we introduced tree co-width, a new graph invariant which is bounded on paths and trees, but also on some dense and symmetric classes of graphs, such as products of small graphs. The theory of tree co-width can be developed purely graph-theoretically. However, category theory provides conceptual guidance and structural insights that allow us to connect decomposition methods to stratified models in applications such as epidemiology, and design efficient algorithms that leverage universal properties. In particular, the category-theoretic treatment turns out to be essential, at least to the authors, in showing that first-order model checking can be performed directly on tree co-decompositions, without constructing their limits, in FPT-time parameterized by the width and the formula (with linear dependence on the number of bags).} 

\subsubsection{Relationships to Feferman-Vaught}

\body{
    It is worth comparing our results to general Feferman-Vaught results. For this we turn to Hodges' textbook~\cite{Hodges_1993} (specifically Equations 6.3, 6.4 and Theorem 9.6.2) where a general Feferman-Vaught theorem is given in terms of certain substructures of products $\prod_{i \in I}A_i$ called $B$-valued substructure, where $B$ is a Boolean algebra whose points are certain subsets of the index set $I$. At first glance, being substructures of products, these may seem similar to pullbacks, but they are in fact not. To see this, first observe that, if  $C\subseteq A_1\times A_2$ is a substructure and let $I=\{1,2\}$, then there are precisely two Boolean algebras of subsets of $I$: $B_1=\{\varnothing,I\}$ and $B_2= 2^I$. Thus $C$ admits gluing over $B_2$ if and only if $C=\pi_{A_1}(C)\times\pi_{A_2}(C)$ (this is already remarked by Hodges~\cite{Hodges_1993}). Furthermore, $C$ is $B_1$-valued if and only if $C \cong A_1 \cong A_2$.  Consequently, arbitrary pullbacks $A_1\times_A A_2$ need not satisfy Hodges' gluing and valuedness conditions for any Boolean algebra on $I$. To make this very concrete, consider any pair of surjections $X\twoheadrightarrow 2\twoheadleftarrow Y$ of sets: their pullback can satisfy Hodges' conditions only in the degenerate case where both maps are bijections.
}

\body{
    What the previous discussion \textit{does} establish is that our \cref{thm: FO model checking on codecompositions} does not follow \textit{verbatim} from what Hodges calls `{generalized Feferman-Vaught}'. However, it doesn't rule out the possibility of first hacking the structure -- by taking FO-interpretations or transductions, say -- and \textit{then} applying known results. Also, perhaps more interestingly, meta-theorems of the kind discussed above are somewhat hazy, as Hodges points out~\cite{Hodges_1993}, in the sense that many different kinds of results get lumped together under the umbrella term `Feferman-Vaught'. This presents a difficulty in evaluating where and how our result fits into that literature and indeed we cannot rule out the possibility obtaining an alternative proof of our algorithmic meta-theorems by leveraging some other kind of Feferman-Vaught-like result.
}

\body{
    Either way, there are some structural differences between typical proof techniques in this area and those we employed here. Whereas typical model-theoretic arguments rely on accounting for logical types, we explicitly manipulate proof-relevant data (meaning strategies and hypotheticals) in our dynamic programming arguments. These may be two sides of the same coin: just as interpretations allow us to reason either by transforming structures or by translating formulas in the opposite direction, preserving  satisfaction~\cite[Proposition 3.5]{immerman2012descriptive}, there may be a similar correspondence between composing logical types and composing the structures that witness their satisfaction. Earnestly engaging with these complexities and providing a more unified treatment of these techniques would be an interesting research program. It is, however, beyond the scope of the present paper. 
}

\subsection{Future Directions}
\body{Tree co-decompositions appear to be a rich topic which already straddles a few different research areas such as graph theory, finite model theory, category theory (as developed in this paper) and other areas including scientific modeling (albeit, framed in a different language; c.f. Libkind et al.~\cite{LibkindEtAl2022Epidemic}) and applied topology (under the guise of cellular shaves; c.f. Curry~\cite{Curry2014}, Hansen and Ghrist~\cite{HansenGhrist2019}, Robinson~\cite{Robinson2017} and Ribeiro, Ten{\'o}rio, Belieni Souza and Mesquita~\cite{ribeiro2026cooperative}). Given this range, it no surprise that further questions abound. Below we highlight several conjectures and open research directions that stem from the current paper and which suit the interests of the authors; this is by no means an exhaustive list.}

\paragraph{Some Questions in Graph Theory}

\body{
    Despite being closely related to tree-width by virtue of categorical duality, tree co-decompositions provide a notion of width that is genuinely distinct from tree-width. Moreover, classes of bounded co-width are distinct from those of bounded twin-width: as shown earlier, there are classes of graphs, such as powers of irreflexive complete graphs, which have bounded tree co-width but unbounded twin-width. The converse direction, however, is unclear.
}

\begin{question}
    Are there classes of bounded twin-width that have unbounded tree co-width?
\end{question}

\body{Although it can be shown that binary trees have bounded path co-width (in contrast to what happens when considering classical path-width), it is not clear whether all trees have bounded path co-width. Our characterization of bounded path co-width in terms of graph grammars (c.f. Apendix~\ref{appendix: grammar}) may provide an approach to this question. We believe that this class may be a candidate for separating path co-width from tree co-width asymptotically.}

\begin{conjecture}
    Trees -- even allowing loops -- have unbounded path co-width.
\end{conjecture}

\body{Although all categorical limits commute with each other, it is still interesting to investigate how tree co-width interacts with other algebraic graph operations.}
\begin{question}
    What can be said about the tree co-width of prime graphs (those not admitting nontrivial factorizations along categorical products)? Furthermore, can the tree co-width of a box-product of graphs be determined from the tree co-width of its factors?
\end{question}

\begin{figure}
    \centering
    \includegraphics[width=0.5\linewidth, page=16]{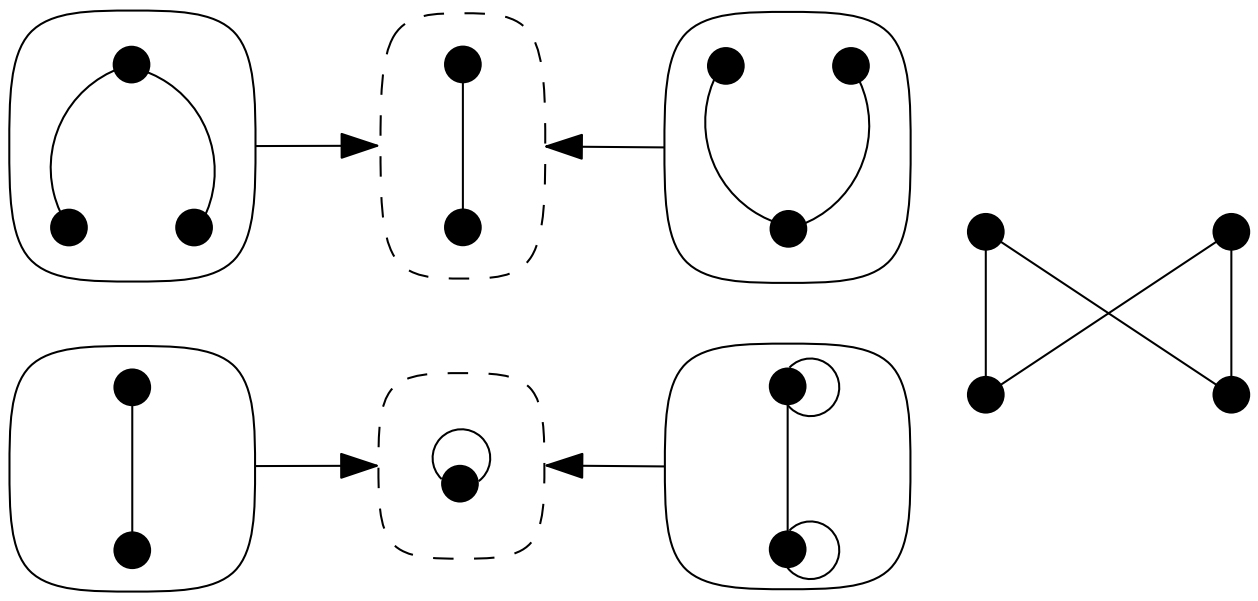}
    \caption{Two genuinely distinct co-decompositions (left) which build the same graph (right).}
    \label{fig:distinct}
\end{figure}

\body{The uniqueness and canonicity of decompositions is another compelling direction. For instance, in \cref{fig:distinct}, we see genuinely different co-decompositions which build the same graph.
}

\body{
It is natural to ask whether there are canonical decompositions into prime pieces, for example an algebraic dual to the Tutte decomposition by Cunningham and Edmonds~\cite{CunninghamEdmonds1980}.}
\begin{question}
    Is there a canonical way to decompose graphs into prime pieces using tree co-decompositions with adhesions of size at most two?
\end{question}

\body{Co-decompositions of graphs pose some unique challenges when compared to classical graph decompositions. Some very intuitive properties of graph decompositions no longer hold in the dual setting. For example, given a pair of graphs $H$ and $G$ and a graph homomorphism $f \colon H\rightarrow G$, we have that any tree decomposition of $G$ can be `pulled back' along $f$ to a tree decomposition of $H$~\cite{BumpusKocsis2023}. More precisely, given a tree decomposition $(T,\{B_t\})$ of $G$, one can define a graph decomposition $(T,\{C_t\})$ of $H$ by setting each of the bags $C_t$ equal to the preimage $f^{-1}(B_t)$ of the bag $B_t$ in $H$. Importantly, the shape of the graph decomposition is preserved. And when $f$ is an inclusion, we get that the size of each bag decreases or stays the same: in other words, tree-width is subgraph monotone. One can view the preimage $f^{-1}(B_t)$ as the categorical pullback of the inclusion $B_t\hookrightarrow G$ and $f$.}

\body{In the classical setting, one can go further. In~\cite{bumpus2024pushing}, Bumpus, Fairbanks and Turner introduced \emph{lassos}, categorical structures which give rise to generalizations of graph contraction maps, and showed that tree-shaped structured decompositions can always be pushed forward along lasso-induced contractions, without increasing their width or changing their shape. In particular, this generalizes the familiar monotonicity of tree-width under graph contractions to a broad categorical setting. However, the arguments rely on the colimit structure of ordinary decompositions, and do not directly extend to co-decompositions, where the underlying object is instead reconstructed as a limit.}

\body{Indeed, it is natural to ask whether the machinery dualises: what happens when we push a co-decomposition forward along a graph homomorphism? Unfortunately, even the most direct construction fails. Pushing the constituent objects of a co-decomposition forward along a homomorphism does not, in general, produce a co-decomposition of the target object, as illustrated in~\autoref{fig:cant_copullback_codecomps}.
\begin{figure}
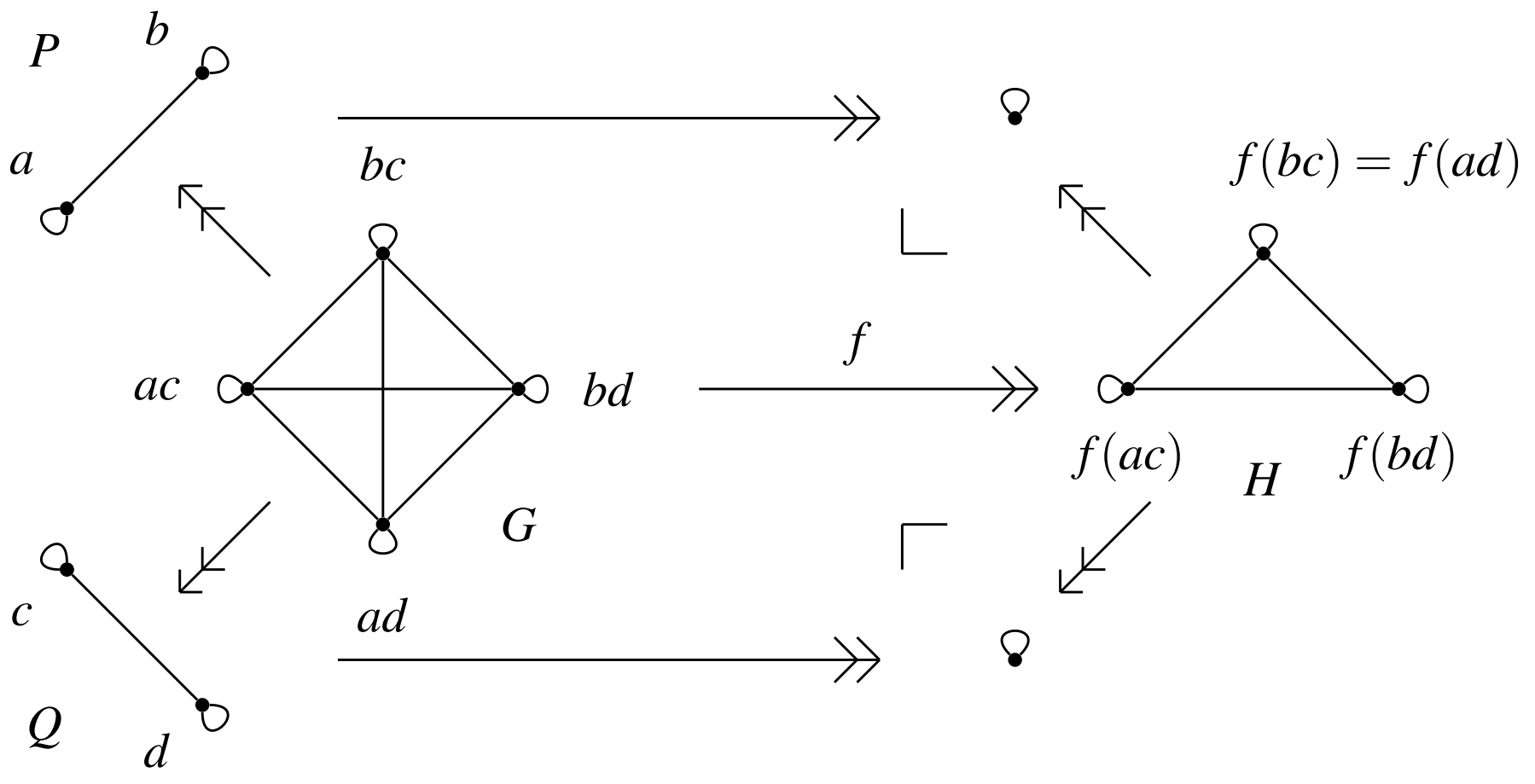

    \centering
    \include{figures/no_pushout}
    \caption{The graph $A$ is obtained from the graphs $P$ and $Q$ by pullback over the terminal graph $o$. Illustrated is the epimorphism $f:A\twoheadrightarrow B$ and the respective componentwise pushouts of the factors of $A$. The resulting diagram does not have $B$ as its pullback.}
    \label{fig:cant_copullback_codecomps}
\end{figure}
}

\body{Nevertheless, this does not rule out the existence of a more suitable class of morphisms under which co-decompositions behave well. Just as graph minors provide a natural containment relation with respect to which tree-width is monotone, it is reasonable to seek an analogous theory for tree co-width. Such a theory might also provide new tools for understanding obstructions to bounded co-width.}
\begin{question}
    Is there a natural analogue of graph minors for which tree co-width is monotone? More strongly, is there a class of morphisms along which tree co-decompositions can be pushed forward without increasing their width or changing their shape?
\end{question}

\paragraph{Some Questions in Algorithmics}

\body{Although it has not been our focus here, computing tree co-decompositions is an important task. In our motivating example in the introduction, Letícia built a stratified-model by specifying a co-decomposition. This means that the co-decomposition \textit{is} the input, while the larger model, the limit of the co-decomposition, is not (also because it can be huge! -- Letícia wouldn't want to write this down). However, from a purely graph-theoretic and algorithmic perspective, it is natural to ask how to compute a co-decomposition from an explicitly given graph. To that end, the computational complexity of determining co-width is largely unexplored; we conjecture it to be hard:}
\begin{conjecture}
    It is NP-hard to decide if a graph has tree co-width at most $k$.
\end{conjecture}

\body{Then again, factoring integers represented in unitary can be done in polynomial time and the same is true for factoring graphs as products~\cite{handbookproductgraphs}, so perhaps the conjecture above might fail. Either way, from a parameterized perspective, we expect it to be possible to efficiently compute minimum-width co-decompositions.}

\begin{conjecture}
    Deciding whether a graph has path (respectively, tree) co-width at most $k$ is in FPT parameterized by path (resp., tree) co-width.
\end{conjecture}

\begin{question}
    Which natural graph problems beyond first-order logic admit FPT algorithms parameterized by tree co-width when a co-decomposition is given? In particular, for which problems can hardness be established already at bounded tree co-width?
\end{question}

\body{Our algorithmic results also suggest several natural extensions. Perhaps the most immediate is counting first-order logic: rather than merely deciding whether a first-order formula is satisfied, can we count the number of assignments satisfying it? We expect that our methods can be extended to address this problem, with appropriate modifications to the dynamic programming tables. The main difficulty is to ensure that a satisfying assignment is counted exactly once, even when it admits multiple witnesses or strategies.}

\begin{conjecture}
    There is an algorithm which, given a tree co-decomposition $\mathcal T$ of width at most $w$ with $b$ bags, and a first-order formula $\phi(\bar x)$, computes
    \[
        \left|\left\{
            \bar a \mid
            \lim\mathcal T \models \phi(\bar a)
        \right\}\right|
    \]
    using $f_\sigma(w,|\phi|)\cdot b$ arithmetic operations, for some computable function $f_\sigma$ depending on the signature.
\end{conjecture}

\body{A related direction is extending our model-checking results to first-order logic with counting quantifiers, or to weighted variants of first-order logic. More broadly, one may ask how far the strategy-based methods developed here extend to other logics. The hardness of \textsc{Path} on implicitly represented graphs suggests that the boundary between tractability and intractability deserves particular attention, although it does not by itself rule out fixed-parameter tractability for more expressive logics.}
\begin{question}
    Which extensions of first-order logic admit fixed-parameter tractable model-checking algorithms on tree co-decompositions?
\end{question}

\paragraph{Categorical connections.}
\body{The construction of strategies and hypotheticals also suggests connections with other categorical approaches to logic. In particular, our use of families carrying explicit witnesses of logical compatibility is reminiscent of proof-relevant logical relations and Artin gluing, as developed, for example, by Sterling and Spitters~\cite{SterlingSpitters2018} and Sterling and Harper~\cite{SterlingHarper2021}. It would be interesting to investigate whether the categories of hypotheticals constructed here admit a natural description in terms of gluing or descent, and whether this perspective yields further algorithmic insights.}

\body{
    Our results also belong to the broader tradition of Feferman--Vaught--Mostowski composition methods. Recent categorical approaches, particularly the game-comonadic work of Jakl, Marsden and Shah~\cite{JaklMarsdenShah2022,JaklMarsdenShah2025}, provide general frameworks for logical composition and its algorithmic applications. An interesting direction is to clarify the relationship between these methods and our computation of supported hypotheticals over pullback-based tree co-decompositions.
    }

\body{
    Finally, consider contravariant powerset construction. Applying this construction to the constituent sets of a co-decomposition yields diagrams of Boolean algebras with compatibility maps running in the opposite direction -- specifically, these are structured decompositions~\cite{StructuredDecompositions} of Boolean Algebras. Taking colimits of these diagrams in categories of partial Boolean algebras, as studied by van den Berg and Heunen~\cite{vanDenBergHeunen2012} and Abramsky and Barbosa~\cite{AbramskyBarbosa2021} and building upon the work of Kochen and Specker~\cite{kochen2011problem}, suggests a possible connection with contextuality and the theory of partial algebras. Understanding this connection more precisely, and whether our algorithmic methods might find applications in that setting, seems a promising direction for future work.
}

\section*{Acknowledgements}

We thank \textit{Fábio Botler} for helpful discussions, \textit{Zoltan Kocsis} for pointing out the connection between our constructions and Artin gluing and for making us aware of categorical accounts of partial Boolean algebras, \textit{Johannes Carmesin} for pointing us towards comparisons with twin-width and \textit{Vinicius Fernandes dos Santos} for motivating us to establish the relationship between tree-width and co-width. 

Benjamin Merlin Bumpus acknowledges support of the \textit{São Paulo Research Foundation (FAPESP), grant 2025/16921-5}. Laura Larios-Jones acknowledges support of the \textit{University of Glasgow Jim Gatheral Travel Scholarship}. Will J.\ Turner acknowledges support by \textit{DFG project number 546892829}.

\section*{AI Methodology}
We used OpenAI's ChatGPT to assist with technical details and occasionally exposition of proofs, specifically when when drafting notation and wrangling index bookkeeping. AI use was also helpful for the construction of figures and illustrative examples; literature review; and brainstorming organizational options for the manuscript. We also used it as an effective sounding board for discussing and searching for counterexamples to conjectures and proposed arguments. All mathematical results and proof strategies were developed by the authors, who, of course, take full responsibility for the contents of the paper. The writing is also our own -- all em-dashes, of which the authors are fond of, were typed by hand.

%% file: figures/no_pushout.tex
\begin{tikzpicture}[scale=0.7]
	\begin{pgfonlayer}{nodelayer}
		\node [style=black dot] (9) at (11, 1.5) {};
		\node [style=black dot] (10) at (9.5, 0) {};
		\node [style=black dot] (11) at (11, -1.5) {};
		\node [style=black dot] (12) at (12.5, 0) {};
		\node [style=black dot] (13) at (7.5, 2) {};
		\node [style=black dot] (14) at (9, 3.5) {};
		\node [style=black dot] (15) at (7.5, -2) {};
		\node [style=black dot] (16) at (7.5, -2) {};
		\node [style=black dot] (17) at (9, -3.5) {};
		\node [style=none] (18) at (8.75, 2.25) {};
		\node [style=none] (19) at (9, 2) {};
		\node [style=none] (20) at (9.75, 1.25) {};
		\node [style=none] (21) at (8.75, -2.25) {};
		\node [style=none] (22) at (9, -2) {};
		\node [style=none] (23) at (9.75, -1.25) {};
		\node [style=none] (24) at (8.75, -2) {};
		\node [style=none] (25) at (9, -1.75) {};
		\node [style=none] (26) at (9, -2.25) {};
		\node [style=none] (27) at (9.25, -2) {};
		\node [style=none] (28) at (8.75, 2) {};
		\node [style=none] (29) at (9, 2.25) {};
		\node [style=none] (30) at (9, 1.75) {};
		\node [style=none] (31) at (9.25, 2) {};
		\node [style=black dot] (32) at (20.75, 1.5) {};
		\node [style=black dot] (33) at (19.25, 0) {};
		\node [style=black dot] (35) at (22.25, 0) {};
		\node [style=black dot] (36) at (18, 3) {};
		\node [style=black dot] (39) at (18, -3) {};
		\node [style=none] (41) at (18.5, 2.25) {};
		\node [style=none] (42) at (18.75, 2) {};
		\node [style=none] (43) at (19.5, 1.25) {};
		\node [style=none] (44) at (18.5, -2.25) {};
		\node [style=none] (45) at (18.75, -2) {};
		\node [style=none] (46) at (19.5, -1.25) {};
		\node [style=none] (47) at (18.5, -2) {};
		\node [style=none] (48) at (18.75, -1.75) {};
		\node [style=none] (49) at (18.75, -2.25) {};
		\node [style=none] (50) at (19, -2) {};
		\node [style=none] (51) at (18.5, 2) {};
		\node [style=none] (52) at (18.75, 2.25) {};
		\node [style=none] (53) at (18.75, 1.75) {};
		\node [style=none] (54) at (19, 2) {};
		\node [style=none] (55) at (14.5, 0) {};
		\node [style=none] (56) at (18.25, 0) {};
		\node [style=none] (58) at (18, -0.25) {};
		\node [style=none] (59) at (18, 0.25) {};
		\node [style=none] (60) at (18, 0) {};
		\node [style=none] (61) at (17.75, 0.25) {};
		\node [style=none] (62) at (17.75, -0.25) {};
		\node [style=none] (63) at (10.5, 3) {};
		\node [style=none] (64) at (16.5, 3) {};
		\node [style=none] (65) at (16.25, 2.75) {};
		\node [style=none] (66) at (16.25, 3.25) {};
		\node [style=none] (67) at (16.25, 3) {};
		\node [style=none] (68) at (16, 3.25) {};
		\node [style=none] (69) at (16, 2.75) {};
		\node [style=none] (70) at (10.5, -3) {};
		\node [style=none] (71) at (16.5, -3) {};
		\node [style=none] (72) at (16.25, -3.25) {};
		\node [style=none] (73) at (16.25, -2.75) {};
		\node [style=none] (74) at (16.25, -3) {};
		\node [style=none] (75) at (16, -2.75) {};
		\node [style=none] (76) at (16, -3.25) {};
		\node [style=none] (77) at (7, 2.5) {$a$};
		\node [style=none] (78) at (8.5, 4) {$b$};
		\node [style=none] (79) at (7, -2.5) {$c$};
		\node [style=none] (80) at (8.5, -4) {$d$};
		\node [style=none] (81) at (8.5, 0) {$ac$};
		\node [style=none] (82) at (11, 2.5) {$bc$};
		\node [style=none] (83) at (13.5, 0) {$bd$};
		\node [style=none] (84) at (11, -2.5) {$ad$};
		\node [style=none] (85) at (22, 2.5) {$f(bc)=f(ad)$};
		\node [style=none] (86) at (19.25, -0.75) {$f(ac)$};
		\node [style=none] (87) at (22.25, -0.75) {$f(bd)$};
		\node [style=none] (88) at (16.25, 0.5) {$f$};
		\node [style=none] (89) at (16.75, 2) {};
		\node [style=none] (90) at (16.75, 1.5) {};
		\node [style=none] (91) at (17.25, 1.5) {};
		\node [style=none] (92) at (17.25, -1.5) {};
		\node [style=none] (93) at (16.75, -1.5) {};
		\node [style=none] (94) at (16.75, -2) {};
		\node [style=none] (95) at (12.5, -1.5) {$G$};
		\node [style=none] (96) at (20.75, -1) {$H$};
		\node [style=none] (97) at (7.25, 3.75) {$P$};
		\node [style=none] (98) at (7.25, -3.75) {$Q$};
	\end{pgfonlayer}
	\begin{pgfonlayer}{edgelayer}
		\draw [style=c0] (14) to (13);
		\draw [style=c0] (10) to (9);
		\draw [style=c0] (9) to (12);
		\draw [style=c0] (12) to (11);
		\draw [style=c0] (11) to (10);
		\draw [style=c0] (10) to (12);
		\draw [style=c0] (9) to (11);
		\draw [style=c0] (17) to (16);
		\draw [style=c0, in=90, out=0, loop] (14) to ();
		\draw [style=c0, in=-180, out=-90, loop] (13) to ();
		\draw [style=c0, in=-180, out=90, loop] (16) to ();
		\draw [style=c0, in=-90, out=0, loop] (17) to ();
		\draw [style=c0, in=-135, out=135, loop] (10) to ();
		\draw [style=c0, in=135, out=45, loop] (9) to ();
		\draw [style=c0, in=-45, out=45, loop] (12) to ();
		\draw [style=c0, in=-45, out=-135, loop] (11) to ();
		\draw [style=c0] (20.center) to (18.center);
		\draw [style=c0] (18.center) to (28.center);
		\draw [style=c0] (18.center) to (29.center);
		\draw [style=c0] (30.center) to (19.center);
		\draw [style=c0] (19.center) to (31.center);
		\draw [style=c0] (21.center) to (24.center);
		\draw [style=c0] (21.center) to (26.center);
		\draw [style=c0] (21.center) to (23.center);
		\draw [style=c0] (25.center) to (22.center);
		\draw [style=c0] (22.center) to (27.center);
		\draw [style=c0] (33) to (32);
		\draw [style=c0] (32) to (35);
		\draw [style=c0] (33) to (35);
		\draw [style=c0, in=45, out=135, loop] (36) to ();
		\draw [style=c0, in=45, out=135, loop] (39) to ();
		\draw [style=c0, in=-135, out=135, loop] (33) to ();
		\draw [style=c0, in=135, out=45, loop] (32) to ();
		\draw [style=c0, in=-45, out=45, loop] (35) to ();
		\draw [style=c0] (43.center) to (41.center);
		\draw [style=c0] (41.center) to (51.center);
		\draw [style=c0] (41.center) to (52.center);
		\draw [style=c0] (53.center) to (42.center);
		\draw [style=c0] (42.center) to (54.center);
		\draw [style=c0] (44.center) to (47.center);
		\draw [style=c0] (44.center) to (49.center);
		\draw [style=c0] (44.center) to (46.center);
		\draw [style=c0] (48.center) to (45.center);
		\draw [style=c0] (45.center) to (50.center);
		\draw [style=c0] (55.center) to (56.center);
		\draw [style=c0] (56.center) to (59.center);
		\draw [style=c0] (56.center) to (58.center);
		\draw [style=c0] (62.center) to (60.center);
		\draw [style=c0] (60.center) to (61.center);
		\draw [style=c0] (63.center) to (64.center);
		\draw [style=c0] (64.center) to (66.center);
		\draw [style=c0] (64.center) to (65.center);
		\draw [style=c0] (69.center) to (67.center);
		\draw [style=c0] (67.center) to (68.center);
		\draw [style=c0] (70.center) to (71.center);
		\draw [style=c0] (71.center) to (73.center);
		\draw [style=c0] (71.center) to (72.center);
		\draw [style=c0] (76.center) to (74.center);
		\draw [style=c0] (74.center) to (75.center);
		\draw [style=c0] (89.center) to (90.center);
		\draw [style=c0] (90.center) to (91.center);
		\draw [style=c0] (92.center) to (93.center);
		\draw [style=c0] (93.center) to (94.center);
	\end{pgfonlayer}
\end{tikzpicture}

%% file: sections/relational_structures.tex
\section{Relational Structures}\label{appendix:RS}

\body{
    A \define{signature} $\sigma$ is a set $S$ of relation symbols and an \define{arity} function $\text{ar}:S\rightarrow \mathbb{N}_0$.
    A (relational) \define{structure} $A$ over a signature $\sigma=(S,\text{ar})$ is a set $A$ along with an \define{interpretation function} $I$ which assigns to each symbol $R\in S$ a relation $R^A$ over $A$ of arity $\text{ar}(R)$. When the context is clear, we abuse notation and write $R$ for the relation $R^A$.
}

\begin{example}
    A directed graph can be encoded as a relational structure in several ways. For a simple directed graph $G$ with vertex set $V(G)$, we may use the signature $\sigma_{\mathrm{E}}=\{\mathsf{E}\}$, where $\operatorname{ar}(\mathsf{E})=2$, and interpret it on $V(G)$ by $\mathsf{E}^G=\{(u,v)\in V(G)^2 : \text{there is an edge from $u$ to $v$ in $G$}\}$. This encoding records adjacency but does not distinguish parallel edges.
\end{example}

\begin{example}
    Alternatively, a directed graph $G=(V,E,s,t)$ with loops and parallel edges may be encoded over the signature $\sigma_{\mathrm{inc}}=\{\mathsf{Vertex},\mathsf{Edge},\mathsf{Source},\mathsf{Target}\}$, where $\mathsf{Vertex}$ and $\mathsf{Edge}$ are unary and $\mathsf{Source}$ and $\mathsf{Target}$ are binary. The underlying set is the coproduct (i.e. disjoint union) $V+E$, with $\mathsf{Vertex}^G=V$, $\mathsf{Edge}^G=E$, $\mathsf{Source}^G=\{(e,s(e)):e\in E\}$, and $\mathsf{Target}^G=\{(e,t(e)):e\in E\}$. This encoding retains the identity of each edge and therefore distinguishes parallel edges.
\end{example}


\body{
    Given two structures $A$ and $B$ over a signature $\sigma=(S,\text{ar})$, a \define{homomorphism} from $A$ to $B$ is a function $f:A\rightarrow B$ so that for each relation $R\in S$ and each tuple $\vec{x}$ from $A$ with compatible arity, we have $R^A(\vec{x})\implies R^B(f(\vec{x}))$. Here $f(\vec{x})$ is computed pointwise. Relational structures and their homomorphisms form a category (c.f. Proposition~\ref{prop: relational structures category}). This category has all limits~\cite[Paragraph~3.10]{parker2022exponentiability}, which is enough for our purposes here.
}

\begin{proposition}\label{prop: relational structures category}
    The structures over a given signature $\sigma$, along with the corresponding homomorphisms between them, form a category with all limits.
\end{proposition}
\begin{proof}
    It is evident it forms a category: homomorphisms of structures consist of certain kinds of functions on the ground sets, so, once one realizes that relations are preserved under composition, associativity and identities are immediate. One can consult Parker~\cite[Paragraph~3.10]{parker2022exponentiability} to see it has limits; for completeness we will give details here. Let $D$ be a small diagram of $\sigma$-structures. Working in $\cat{Set}$, take the limit $L$ of the underlying point sets with projections $\pi_j\colon L\to\pts{D(j)}$. For each relation symbol $R\in S$ of arity $k$, define
    \[
        R^L=\{\vec{x}\in L^k:R^{D(j)}(\pi_j(\vec{x}))\text{ for every }j\}.
    \]
    By construction, the projections $\pi_j:L\rightarrow D(j)$ are homomorphisms and hence form a cone over $D$. Let $(X,\{q_j\})$ be any other cone over $D$. The universal property of the limit of the underlying point sets gives a unique function $q\colon\pts{X}\rightarrow L$ such that $\pi_j\circ q=q_j$ for every $j$. If $R^X(\vec{x})$, then $R^{D(j)}(q_j(\vec{x}))$ for every $j$, since each $q_j$ is a homomorphism. It follows from the definition of $R^L$ that $R^L(q(\vec{x}))$, and hence $q$ is a homomorphism. Therefore $L$ is the limit of $D$.
\end{proof}

\subsection{Modeling Equality}
\label{subsec:model_eq}

\body{In the framework of relational structures, we cannot ask directly about equality of points in a first-order formula. However, our methods extend in the following way: We can model equality using an additional binary relation that appears as a loop on every point in the input co-decomposition/structure. We note that it is possible for there to be structures in our calculations (in particular inside the heads of hypotheses), for which this binary relation relates two distinct points. This technique is made precise in Observation~\ref{loop_trick}.}

\body{For a structure $G$, let $G^\circ$ be the structure obtained from $G$ by adding a new relation $\circ\subseteq\pts{G}^2$ which is reflexive on every point but relates no distinct points. That is, it is a loop everywhere. For a co-decomposition $\Ta$, we denote by $\Ta^\circ$ the co-decomposition obtained from $\Ta$ by replacing each structure by its variant with this new relation.}

\begin{observation}\label{loop_trick}
   Let $\Ta$ be a co-decomposition. Then $\lim(\Ta^\circ)=(\lim(\Ta))^\circ$. 
\end{observation}

\body{We can use our machinery on $\Ta^\circ$ to ask questions about equality on $\lim\Ta$. That is, evaluating FO with equality requires only a linear amount of pre-processing on the co-decomposition.}

%% file: sections/appendix_binary.tex
\section{Making Tree Co-Decompositions Binary}\label{appendix:binary}

\begin{definition}
    A tree co-decomposition $\Ta=(T,\{G_a\},\{g_t^e\})$ is \define{binary tree co-decomposition} if $T$ is a binary tree.
\end{definition}

\begin{lemma}
    Let $\Ta$ be a tree co-decomposition. Then there is a binary tree co-decomposition $\Ta^*$ so that $\lim \Ta=\lim\Ta^*$ as well as $w(\Ta)=w(\Ta^*)$ and, if $b$ is the number of bags in $\Ta$, then the number of bags in~$\Ta^*$ is $O(b\log b)$.
\end{lemma}
\begin{proof}
Write $\Ta=(T,\{G_a\},\{g_t^e\})$. For each node $t\in V(T)$, let $\Ta^t=(T^t,\{G^t_a=G_t\},\{\text{id}_{G_t}\})$ be a binary tree co-decomposition whose auxiliary tree has its leaves labelled $v_e^t$ for each edge $e$ incident with $t$ in $T$, whose bags and adhesions are all equal to $G_t$ and whose maps are all the identity on $G_t$. Note that it is easy to choose~$\Ta^t$ so that it has at most $O(\log d(t))$ bags. 

Let $T^*$ be the tree obtained from the disjoint union of the trees $T^t$ by adding for each $e\in E(T)$, with endvertices $s$ and $t$ in $T$, an edge labelled by $e$ between the leaf $v_e^t$ of $T^t$ and the leaf $v_e^s$ of $T^s$. We let $\Ta'$ be the tree co-decomposition whose auxiliary tree is $T^*$, whose bags and adhesions are \[
\{G_e\}_{e\in E(T)}\cup \bigcup_{t\in V(T)}\{G^t_a:=G_t\}_{a\in V(T^t)\cup E(T^t)}\]
and whose maps are
\[
\bigcup_{t\in V(T)}\{g_{v_e^t}^e:=g_{t}^e\,|\,\text{$e\in E(T)$ incident in $T$ with $t$}\}
\]
along with the identity maps within each subtree $T^t$.

Observe that the original tree co-decomposition $\Ta$ is exactly the tree co-decomposition obtained as a partial evaluation of $\Ta^*$ by contracting for each $t\in V(T)$ the subtree $T^t$ in $T^*$. Hence $\lim \Ta^*=\lim \Ta$. The other required conditions follow by construction.
\end{proof}

%% file: sections/appendix_tree_vs_co.tex
\section{Tree Co-Width vs Tree-width}\label{appendix:tree-co-vs-tree}
In this subsection we prove the following.
\begin{lemma}
    Let $G$ be a graph that admits a tree decomposition with bags of size at most $w\in\mathbb N$. Then $G$ admits a tree co-decomposition of width at most $w+3$. In formulas:
    \[
    \text{tree co-width}(G)\le \text{treewidth}(G)+2
    \]
\end{lemma}
\begin{proof}
    Let $(T,\{B_a\}_{a\in V(T)\cup E(T)})$ be a tree decomposition of $G$ of width at most $w-1$. That is, $T$ is a tree and the collection $\{B_a\}$ is a family of induced subgraphs indexed by the nodes and edges\footnote{We aknowledge that this is not the usual definition of tree-decomposition that appears in the literature however, for our purposes, the equivalence is not hard to see.} of $T$ so that the following criteria hold:\vspace{-.25cm}
    \begin{enumerate}
        \item together the bags $\{B_a\}$ cover $G$; and\vspace{-.25cm}
        \item for each vertex $x\in V(G)$, the indices $a\in V(T)\cup E(T)$ for which $x\in V(B_a)$ form a contiguous subtree of $T$; and\vspace{-.25cm}
        \item for every $a\in V(T)\cup E(T)$, we have $|B_a|\le w$.\vspace{-.25cm}
    \end{enumerate}
    We assume using standard methods that that the maximum degree of $T$ is $3$. We will construct a tree co-decomposition $\Ta:=(T,\{G_a\},\{g_e^t\})$, having the same auxiliary tree as the input tree-decomposition, as follows: For each $a\in V(T)\cup E(T)$, let $V_a$ be the partition of $V(G)$ with the following parts:\vspace{-.25cm}
    \begin{enumerate}
        \item for each $x\in B_a$, we have $\{x\}\in V_a$; and\vspace{-.25cm}
        \item for each connected component $C$ of $T-a$, let $P_C\subseteq V(G)$ be the set of all vertices of $G$ in bags indexed in $C$ that are not in $V(B_a)$, and include $P_C\in V_a$ if $P_C$ is nonempty.\vspace{-.25cm}
    \end{enumerate}
    Let $G_a$ be the graph with vertex set $V_a$ (whose elements we call \emph{part-vertices}) and an edge between two part-vertices $X\subseteq V(G)$ and $Y\subseteq V(G)$ if there is an edge in $G$ between some $x\in X$ and $y\in Y$. Note that there are at most $|B_a|+3$ part-vertices in $G_a$; on for each each $B_a$ and one vertex for each neighbour of $a$.

    Now let $t\in V(T)$ and $e\in E(T)$ be an incident edge. Observe that for each part-vertex $X\in V_t$, there is a unique part-vertex $Y\in V_e$ so that $X\subseteq Y$. Let $g_t^e\colon V_t\to V_e$ denote the map defined this way, and observe that it defines an epimorphism (surjective homomorphism) $G_t\to G_e$.

    We will show that the limit of $\Ta$ is $G$. To that end, for each $a\in V(T)\cup E(T)$, let $g_a\colon G\to G_a$ be the epimorphism sending each vertex to the part-vertex containing it, and observe that $(G,\{g_a\})$ is a cone over $\Ta$ with $G$ as its apex. It remains to prove that $G$ is universal.
    
    Let $(H,\{h_a\})$ be another cone over $\Ta$. Let $u\in V(H)$. We make the following claims:

    \textbf{Claim.} There exists a node $t\in V(T)$ for which $h_t(u)$ is a singleton.
    
    \textit{Proof of claim.} Let $t\in V(T)$ be arbitrary. If $h_t(u)$ is a singleton, then we are done, so suppose this is not the case. Then $h_t(u)=P_C$ for some connected component $C$ of $T-t$. Let $tt'$ be the unique edge of $T$ incident with $C$ but not contained within it. Then $P_C$ is also one of the parts of $G_{tt'}$ and $h_{tt'}(u)=P_C$ by the commutativity of the cone $(H,\{h_a\})$. Since otherwise we would be done, assume that $h_{t'}(u)$ is not a singleton. Then $h_{t'}(u)=P_D$ for some $D\subseteq C$ a connected component of $T-t'$. Starting again with $P_D$ in place of $P_C$, we iterate. Since $T$ is finite, this process must terminate in finding our required singleton.
    \qed

    \textbf{Claim.} There exists a unique vertex $x\in V(G)$ so that $x\in h_a(u)$ for all $a\in V(T)\cup E(T)$.
    
    \textit{Proof of claim.} By the first claim, $h_a(u)$ is a singleton $\{x\}$ for some $a\in V(T)$ and $x\in V(G)$. Observe that by the commutativity of the cones, $h_a(u)=\{x\}$ for each $a$ in the contiguous subtree $T[x]$ of $T$ whose bags contain $x$. For any $a\in V(T)\cup E(T)$ not in $T[x]$, it is clear by induction that $h_a(u)$ is the unique part representing the connected component of $T-a$ containing $T[x]$. That is, the part of $V_a$ containing $x$. Hence the choice of $x$ satisfies our claim.
    \qed

    Observe that the map that sends each $u\in V(H)$ to $x\in V(G)$ gives a function $f\colon H\to G$ that commutes with both cones, and that furthermore this function is unique with this property.
    It remains to show that $f$ is a graph homomorphism. So let $(u,v)$ be an edge in $H$. Observe that there is an edge between from $h_a(u)=g_a(f(u))$ to $h_a(v)g_a(f(v))$ in $G_a$ for each $a\in V(T)\cup E(T)$. Hence $(f(u),f(v))$ is an edge in $G$. Hence $f$ is a homomorphism and we conclude that $G$ is the limit of $\Ta$.

    Our proof is completed by observing that for every $a\in V(T)\cup E(T)$, the forest $T-a$ has at most three connected components, and hence the graph $G_a$ has at most $|B_a|+3\le w-1+3=w+2$ vertices.
\end{proof}

%% file: sections/appendix_3_build.tex
\section{Every Graph is the Limit of Graphs with at most $3$ Vertices}\label{appendix: tree shape important}

\begin{figure}
    \centering
    \input{graph_co.tikz}
    \caption{An illustration of the construction in Lemma~\ref{lemm:thin_build}}
    \label{fig:graph_co}
\end{figure}
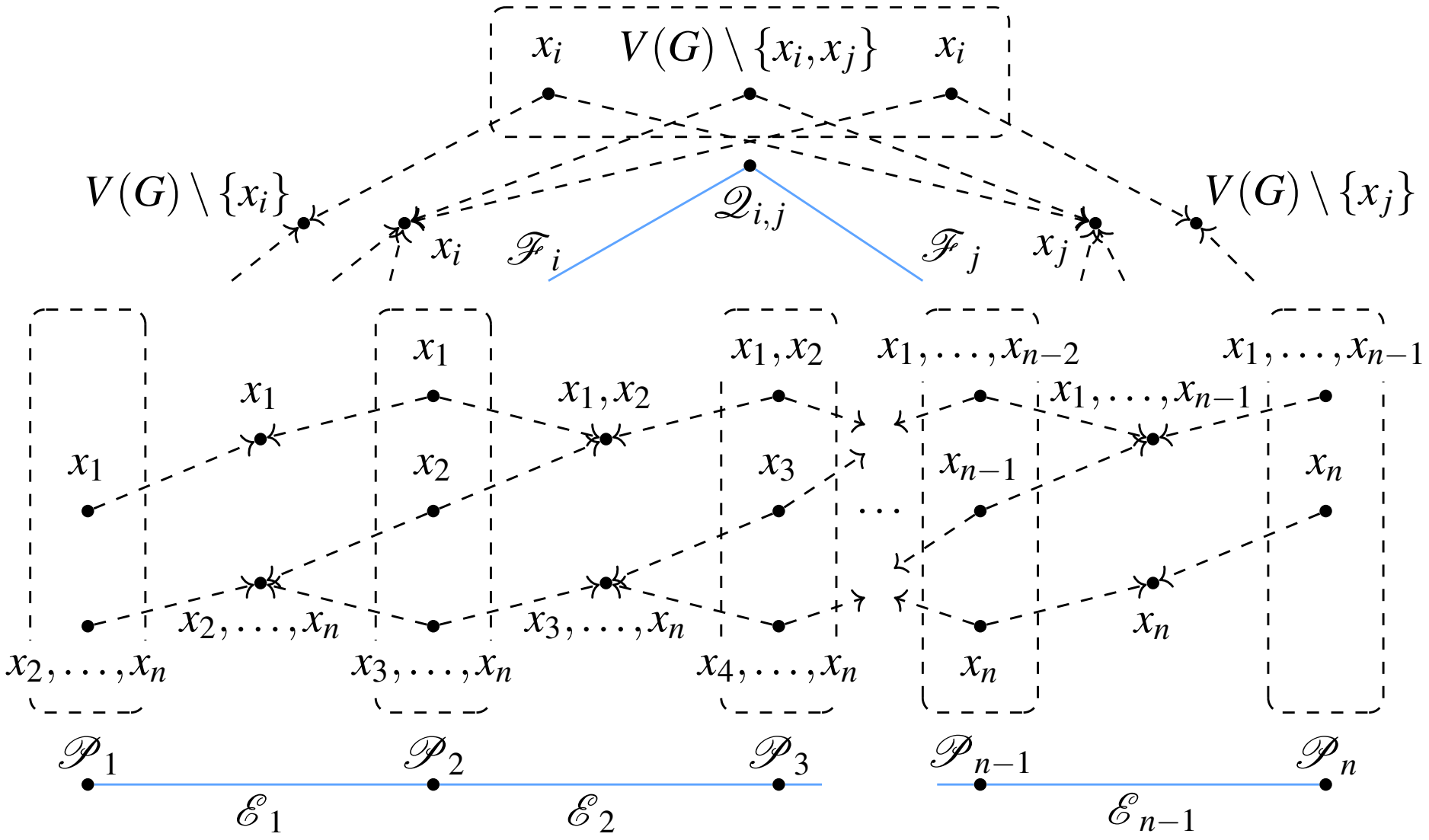

\begin{lemma}\label{lemm:thin_build}
    Every graph is isomorphic to a limit of a graph co-decomposition whose bags and adhesions each contain at most three vertices.
\end{lemma}
\begin{proof}
    See Figure~\ref{fig:graph_co}.
    Let $G$ be a graph. Let $(x_1,\ldots,x_n)$ be an ordering of the vertex set of $G$. Since otherwise it is trivial, suppose that $n\ge 3$. Define the following partitions of $V(G)$:
    
    Let $\Pa_i=\{\;\{x_1,\ldots,x_{i-1}\},\;\{x_{i}\},\;\{x_{i+1},\ldots,x_n\}\;\}$ for each $i\in\{2,\ldots,n-1\}$,
    
    $\Pa_1=\{\;\{x_1\},\;\{x_{2},\;\ldots,x_n\}\;\}\quad$ and $\quad\Pa_n=\{\;\{x_1,\ldots,x_{n-1}\},\;\{x_{n}\}\;\}$.
    
    Let $\Ea_i=\{\;\{x_1,\ldots,x_{i}\},\;\{x_{i+1},\ldots,x_n\}\;\}$ for each $i=[n-1]$.
    
    Let $\Da_i=\{\;\{x_i\},\; V(G)\setminus\{x_{i}\}\;\}$ for each $i\in[n]$.
    
    Let $\Qa_{i,j}=\{\;\{x_i\},\;\{x_j\},\; V(G)\setminus\{x_i,x_j\}\;\}$ for each $i,j\in[n]$.
    
    \noindent For a partition $\Xa$ of $V(G)$, denote by $G/\Xa$ the graph with vertices the parts of~$\Xa$ and an edge between two parts (perhaps as a loop between the same two parts) when there is an edge in $G$ between a vertex in each of the parts. There is a homomorphism from $G$ to $G/\Xa$ taking each vertex to its part, which we call here the \define{projection to $G/\Xa$}.

    Let $H$ be the graph with the following vertices: $p_i$ for $i\in[n]$, $q_{i,j}$ for $i,j\in[n]$ and edges $e_i$ spanning $(p_i,p_{i+1})$ for $i\in[n-1]$ and $d_{i,j}^k$ for $k\in\{i,j\}$ and $i,j\in[n]$ spanning $(p_k,q_{i,j})$. For $i\in [n]$, let $G_{p_i}$ be a copy of $G/\Pa_i$, let $G_{q_{i,j}}$ be a copy of $G/\Qa_{i,j}$, let $G_{e_i}$ be a copy of $G/\Ea_i$ and $G_{d_{i,j}^k}$ be a copy of $G/\Da_k$. Let $\Ha:=(H,\{G_a\},\{g_t^e\})$ be the graph co-decomposition given by these definitions with maps defined so that a vertex $x$ in $G/\Xa$ is mapped to a vertex $y$ in $G/\Ya$ if the part of $\Xa$ corresponding to $x$ is a subset of the part of $\Ya$ corresponding to $y$. By definition these maps commute with the projections to each quotient, which we denote by $g_t:G\rightarrow G_t$ for each $t\in V(H)$. Observe that the pair $(G,\{g_t\})$ is a cone over $\Ha$. We will prove that it is a limit cone.

    Let $(G',\{f_t\})$ an arbitrary cone over $\Ha$. We are required to show that there is a unique homomorphism $\phi:G'\rightarrow G$ which commutes with the maps $\{g_t\}$ and $\{f_t\}$.

    \noindent\textbf{Claim.} For every $u\in V(G')$, there is a unique $i\in [n]$ so that the image of $u$ in $G_a$ is the partition containing $x_i$ for each $a\in V(H)\cup E(H)$.
    
    \noindent\textit{Proof of claim.} Let $j\in [n]$. If $u$ is mapped by $f_{p_j}$ to $\{x_1,\ldots,x_{j-1}\}$ in $G_{p_j}=\Pa_j$, then it is mapped to $\{x_1,\ldots,x_{j}\}$ in $G_{e_{j}}=\Ea_j$ and hence to $\{x_1,\ldots,x_{j+1}\}$ in $G_{p_{j+1}}=\Pa_{j+1}$. If $u$ is mapped instead by $f_{p_j}$ to $\{x_{j}\}$ in $G_{p_j}=\Pa_j$, then it is mapped to $\{x_1,\ldots,x_{j}\}$ in $G_{e_{j}}=\Ea_j$ and $\{x_j,\ldots,x_n\}$ in $G_{e_{j-1}}=\Ea_{j-1}$ and hence to $\{x_1,\ldots,x_{j+1}\}$ in $G_{p_{j+1}}=\Pa_{j+1}$ and $\{x_j,\ldots,x_{n}\}$ in $G_{p_{j-1}}=\Pa_{j-1}$. Finally, if $u$ is mapped by $f_{p_j}$ to $\{x_{j+1},\ldots,x_n\}$ in $G_{p_j}=\Pa_j$, then it is mapped to $\{x_j,\ldots,x_n\}$ in $G_{e_{j-1}}=\Ea_{j-1}$ and hence to $\{x_j,\ldots,x_n\}$ in $G_{p_{j-1}}=\Pa_{j-1}$.

    We conclude that there is a unique $i\in[n]$ so that $u$ is mapped by $g_{p_i}$ to the singleton $\{x_i\}$ in $G_{p_i}=\Pa_i$ and for $k<i$, to the initial segment in $G_{p_k}=\Pa_k$ and for~$k>i$ to the final segment in $G_{p_k}=\Pa_k$.
    
    It remains to check the bags with indices $q_{k,j}$. This is routine after observing that~$u$ must map to the singleton $\{x_i\}$ in $G_{d_{i,j}^i}$ and the complement of $\{x_j\}$ in $G_{d_{i,j}^j}$ for each~$j\neq i$. This completes the proof of our claim. \hfill$\lozenge$

    Let $\phi:G'\rightarrow G$ be the map sending each vertex $u$ to the vertex $x_i$ for the $i$ given by the above claim. By construction, when we restrict only to vertices, this map commutes with the cones. Furthermore, it is easy to see that it is the unique map with this property. It remains to check $\phi$ defines a graph homomorphism. Suppose there is an edge of $G'$ which spans the vertices $u$ and $v$. Then there are indices~$i$ and~$j$ so that $\phi(u)=x_i$ and $\phi(v)=x_j$. Now consider the map $f_{q_{i,j}}:G'\rightarrow G_{q_{i,j}}$. Since this is a homomorphism, there must be an edge between the singletons $\{x_i\}$ and $\{x_j\}$ in~$G_{q_{i,j}}$. By definition, this only occurs when there is an edge from $x_i$ to $x_j$ in $G$. This completes our proof. 
\end{proof}

\begin{example}
    The `three' in Lemma~\ref{lemm:thin_build} is tight. Indeed, the graph co-decompositions with bags and adhesion of size at most two have limits which are direct products of graphs with at most two vertices each. Basic graphs such as the triangle are not realised as such.
\end{example}

%% file: figures/graph_co.tikz
\begin{tikzpicture}
	\begin{pgfonlayer}{nodelayer}
		\node [style=none] (99) at (1, 1.25) {};
		\node [style=none] (100) at (1.25, 1) {};
		\node [style=none] (101) at (1, 3) {};
		\node [style=none] (102) at (1.25, 3.25) {};
		\node [style=none] (103) at (9.75, 3.25) {};
		\node [style=none] (104) at (10, 3) {};
		\node [style=none] (105) at (10, 1.25) {};
		\node [style=none] (106) at (9.75, 1) {};
		\node [style=none] (172) at (-6, -4.75) {$x_1$};
		\node [style=none] (173) at (-6, -8.25) {$x_2,\ldots,x_n$};
		\node [style=none] (174) at (-7, -8.75) {};
		\node [style=none] (175) at (-6.75, -9) {};
		\node [style=none] (176) at (-7, -2.25) {};
		\node [style=none] (177) at (-6.75, -2) {};
		\node [style=none] (178) at (-5.25, -2) {};
		\node [style=none] (179) at (-5, -2.25) {};
		\node [style=none] (180) at (-5, -8.75) {};
		\node [style=none] (181) at (-5.25, -9) {};
		\node [style=none] (182) at (-3, -3.5) {$x_1$};
		\node [style=none] (183) at (-3, -7.5) {$x_2,\ldots,x_n$};
		\node [style=none] (192) at (0, -2.75) {$x_1$};
		\node [style=none] (193) at (0, -8.25) {$x_3,\ldots,x_n$};
		\node [style=none] (194) at (0, -4.75) {$x_2$};
		\node [style=black dot] (195) at (-6, -5.5) {};
		\node [style=black dot] (196) at (-6, -7.5) {};
		\node [style=black dot] (197) at (0, -3.5) {};
		\node [style=black dot] (198) at (0, -5.5) {};
		\node [style=black dot] (199) at (0, -7.5) {};
		\node [style=black dot] (200) at (-3, -4.25) {};
		\node [style=black dot] (201) at (-3, -6.75) {};
		\node [style=none] (210) at (6, -2.75) {$x_1,x_2$};
		\node [style=none] (211) at (6, -8.25) {$x_4,\ldots,x_n$};
		\node [style=none] (212) at (6, -4.75) {$x_3$};
		\node [style=black dot] (213) at (3, -4.25) {};
		\node [style=black dot] (214) at (6, -3.5) {};
		\node [style=black dot] (215) at (3, -6.75) {};
		\node [style=none] (224) at (3, -3.5) {$x_1,x_2$};
		\node [style=none] (225) at (3, -7.5) {$x_3,\ldots,x_n$};
		\node [style=black dot] (226) at (6, -5.5) {};
		\node [style=black dot] (227) at (6, -7.5) {};
		\node [style=none] (231) at (-1, -8.75) {};
		\node [style=none] (232) at (-0.75, -9) {};
		\node [style=none] (233) at (-1, -2.25) {};
		\node [style=none] (234) at (-0.75, -2) {};
		\node [style=none] (235) at (0.75, -2) {};
		\node [style=none] (236) at (1, -2.25) {};
		\node [style=none] (237) at (1, -8.75) {};
		\node [style=none] (238) at (0.75, -9) {};
		\node [style=none] (239) at (5, -8.75) {};
		\node [style=none] (240) at (5.25, -9) {};
		\node [style=none] (241) at (5, -2.25) {};
		\node [style=none] (242) at (5.25, -2) {};
		\node [style=none] (243) at (6.75, -2) {};
		\node [style=none] (244) at (7, -2.25) {};
		\node [style=none] (245) at (7, -8.75) {};
		\node [style=none] (246) at (6.75, -9) {};
		\node [style=none] (249) at (9.5, -2.75) {$x_1,\ldots,x_{n-2}$};
		\node [style=none] (250) at (9.5, -8.25) {$x_n$};
		\node [style=none] (251) at (9.5, -4.75) {$x_{n-1}$};
		\node [style=black dot] (252) at (9.5, -3.5) {};
		\node [style=black dot] (253) at (9.5, -5.5) {};
		\node [style=black dot] (254) at (9.5, -7.5) {};
		\node [style=none] (255) at (8.5, -8.75) {};
		\node [style=none] (256) at (8.75, -9) {};
		\node [style=none] (257) at (8.5, -2.25) {};
		\node [style=none] (258) at (8.75, -2) {};
		\node [style=none] (259) at (10.25, -2) {};
		\node [style=none] (260) at (10.5, -2.25) {};
		\node [style=none] (261) at (10.5, -8.75) {};
		\node [style=none] (262) at (10.25, -9) {};
		\node [style=none] (267) at (15.5, -2.75) {$x_1,\ldots,x_{n-1}$};
		\node [style=none] (268) at (15.5, -4.75) {$x_n$};
		\node [style=none] (269) at (14.5, -8.75) {};
		\node [style=none] (270) at (14.75, -9) {};
		\node [style=none] (271) at (14.5, -2.25) {};
		\node [style=none] (272) at (14.75, -2) {};
		\node [style=none] (273) at (16.25, -2) {};
		\node [style=none] (274) at (16.5, -2.25) {};
		\node [style=none] (275) at (16.5, -8.75) {};
		\node [style=none] (276) at (16.25, -9) {};
		\node [style=black dot] (277) at (15.5, -3.5) {};
		\node [style=black dot] (278) at (15.5, -5.5) {};
		\node [style=black dot] (279) at (12.5, -4.25) {};
		\node [style=black dot] (280) at (12.5, -6.75) {};
		\node [style=none] (281) at (12.5, -3.5) {$x_1,\ldots,x_{n-1}$};
		\node [style=none] (282) at (12.5, -7.5) {$x_n$};
		\node [style=none] (283) at (-6, -9.75) {$\Pa_1$};
		\node [style=none] (284) at (0, -9.75) {$\Pa_2$};
		\node [style=none] (285) at (-3, -10.75) {$\Ea_1$};
		\node [style=none] (286) at (6, -9.75) {$\Pa_3$};
		\node [style=none] (287) at (2.75, -10.75) {$\Ea_2$};
		\node [style=none] (290) at (9.5, -9.75) {$\Pa_{n-1}$};
		\node [style=none] (291) at (15.5, -9.75) {$\Pa_n$};
		\node [style=none] (292) at (12.5, -10.75) {$\Ea_{n-1}$};
		\node [style=black dot] (293) at (-2.25, -0.5) {};
		\node [style=black dot] (294) at (-0.5, -0.5) {};
		\node [style=black dot] (295) at (13.25, -0.5) {};
		\node [style=black dot] (296) at (11.5, -0.5) {};
		\node [style=none] (297) at (0.25, -1) {$x_i$};
		\node [style=none] (298) at (-4.25, 0) {$V(G)\setminus\{x_i\}$};
		\node [style=none] (299) at (10.75, -1) {$x_j$};
		\node [style=none] (300) at (15.25, 0) {$V(G)\setminus\{x_j\}$};
		\node [style=black dot] (301) at (2, 1.75) {};
		\node [style=black dot] (302) at (5.5, 1.75) {};
		\node [style=none] (303) at (2, 2.5) {$x_i$};
		\node [style=none] (304) at (5.5, 2.5) {$V(G)\setminus\{x_i,x_j\}$};
		\node [style=black dot] (305) at (9, 1.75) {};
		\node [style=none] (306) at (9, 2.5) {$x_i$};
		\node [style=none] (307) at (-3.5, -1.5) {};
		\node [style=none] (308) at (-1.75, -1.5) {};
		\node [style=none] (309) at (-0.75, -1.5) {};
		\node [style=none] (310) at (11.25, -1.5) {};
		\node [style=none] (311) at (12, -1.5) {};
		\node [style=none] (312) at (14.25, -1.5) {};
		\node [style=none] (313) at (1.75, -1) {$\Fa_i$};
		\node [style=none] (314) at (9, -1) {$\Fa_j$};
		\node [style=none] (315) at (5.5, -0.25) {$\Qa_{i,j}$};
		\node [style=none] (316) at (7.5, -4) {};
		\node [style=none] (317) at (7.5, -4.5) {};
		\node [style=none] (318) at (7.5, -7) {};
		\node [style=none] (319) at (8, -6.5) {};
		\node [style=none] (320) at (8, -7) {};
		\node [style=none] (321) at (8, -4) {};
		\node [style=none] (322) at (7.75, -5.5) {$\ldots$};
		\node [style=black dot] (323) at (-6, -10.25) {};
		\node [style=black dot] (324) at (0, -10.25) {};
		\node [style=black dot] (325) at (6, -10.25) {};
		\node [style=black dot] (326) at (9.5, -10.25) {};
		\node [style=black dot] (327) at (15.5, -10.25) {};
		\node [style=black dot] (328) at (5.5, 0.5) {};
		\node [style=none] (329) at (8.5, -1.5) {};
		\node [style=none] (330) at (2, -1.5) {};
		\node [style=none] (331) at (6.75, -10.25) {};
		\node [style=none] (332) at (8.75, -10.25) {};
		\node [style=none] (334) at (-7, -7.75) {};
		\node [style=none] (335) at (-1, -7.75) {};
		\node [style=none] (336) at (1, -7.75) {};
		\node [style=none] (337) at (8.5, -3.25) {};
		\node [style=none] (338) at (10.5, -3.25) {};
		\node [style=none] (339) at (14.5, -3.25) {};
		\node [style=none] (340) at (16.5, -3.25) {};
		\node [style=none] (341) at (-5, -7.75) {};
		\node [style=none] (343) at (5, -3.25) {};
		\node [style=none] (344) at (7, -3.25) {};
		\node [style=none] (345) at (5, -7.75) {};
		\node [style=none] (346) at (7, -7.75) {};
	\end{pgfonlayer}
	\begin{pgfonlayer}{edgelayer}
		\draw [style=d1] (102.center)
			 to [bend right=45] (101.center)
			 to (99.center)
			 to [bend right=45] (100.center)
			 to (106.center)
			 to [bend right=45] (105.center)
			 to (104.center)
			 to [bend right=45] (103.center)
			 to cycle;
		\draw [style=d1, bend right=45] (177.center) to (176.center);
		\draw [style=d1, bend right=45] (174.center) to (175.center);
		\draw [style=d1] (175.center) to (181.center);
		\draw [style=d1, bend right=45] (181.center) to (180.center);
		\draw [style=d1, bend right=45] (179.center) to (178.center);
		\draw [style=d1] (178.center) to (177.center);
		\draw [style=d1 arrow] (195) to (200);
		\draw [style=d1 arrow] (196) to (201);
		\draw [style=d1 arrow] (198) to (201);
		\draw [style=d1 arrow] (199) to (201);
		\draw [style=d1 arrow] (197) to (200);
		\draw [style=d1 arrow] (197) to (213);
		\draw [style=d1 arrow] (198) to (213);
		\draw [style=d1 arrow] (199) to (215);
		\draw [style=d1 arrow] (226) to (215);
		\draw [style=d1 arrow] (214) to (213);
		\draw [style=d1 arrow] (227) to (215);
		\draw [style=d1, bend right=45] (234.center) to (233.center);
		\draw [style=d1, bend right=45] (231.center) to (232.center);
		\draw [style=d1] (232.center) to (238.center);
		\draw [style=d1, bend right=45] (238.center) to (237.center);
		\draw [style=d1, bend right=45] (236.center) to (235.center);
		\draw [style=d1] (235.center) to (234.center);
		\draw [style=d1, bend right=45] (242.center) to (241.center);
		\draw [style=d1, bend right=45] (239.center) to (240.center);
		\draw [style=d1] (240.center) to (246.center);
		\draw [style=d1, bend right=45] (246.center) to (245.center);
		\draw [style=d1, bend right=45] (244.center) to (243.center);
		\draw [style=d1] (243.center) to (242.center);
		\draw [style=d1, bend right=45] (258.center) to (257.center);
		\draw [style=d1, bend right=45] (255.center) to (256.center);
		\draw [style=d1] (256.center) to (262.center);
		\draw [style=d1, bend right=45] (262.center) to (261.center);
		\draw [style=d1, bend right=45] (260.center) to (259.center);
		\draw [style=d1] (259.center) to (258.center);
		\draw [style=d1, bend right=45] (272.center) to (271.center);
		\draw [style=d1, bend right=45] (269.center) to (270.center);
		\draw [style=d1] (270.center) to (276.center);
		\draw [style=d1, bend right=45] (276.center) to (275.center);
		\draw [style=d1, bend right=45] (274.center) to (273.center);
		\draw [style=d1] (273.center) to (272.center);
		\draw [style=d1 arrow] (277) to (279);
		\draw [style=d1 arrow] (278) to (280);
		\draw [style=d1 arrow] (252) to (279);
		\draw [style=d1 arrow] (253) to (279);
		\draw [style=d1 arrow] (254) to (280);
		\draw [style=d1 arrow] (301) to (293);
		\draw [style=d1 arrow] (302) to (294);
		\draw [style=d1 arrow] (305) to (294);
		\draw [style=d1 arrow] (301) to (296);
		\draw [style=d1 arrow] (305) to (295);
		\draw [style=d1 arrow] (302) to (296);
		\draw [style=d1 arrow] (307.center) to (293);
		\draw [style=d1 arrow] (308.center) to (294);
		\draw [style=d1 arrow] (309.center) to (294);
		\draw [style=d1 arrow] (310.center) to (296);
		\draw [style=d1 arrow] (311.center) to (296);
		\draw [style=d1 arrow] (312.center) to (295);
		\draw [style=d1 arrow] (214) to (316.center);
		\draw [style=d1 arrow] (226) to (317.center);
		\draw [style=d1 arrow] (227) to (318.center);
		\draw [style=d1 arrow] (252) to (321.center);
		\draw [style=d1 arrow] (253) to (319.center);
		\draw [style=d1 arrow] (254) to (320.center);
		\draw [style=c1] (323) to (324);
		\draw [style=c1] (324) to (325);
		\draw [style=c1] (326) to (327);
		\draw [style=c1] (328) to (329.center);
		\draw [style=c1] (328) to (330.center);
		\draw [style=c1] (325) to (331.center);
		\draw [style=c1] (326) to (332.center);
		\draw [style=d1] (176.center) to (334.center);
		\draw [style=d1] (233.center) to (335.center);
		\draw [style=d1] (236.center) to (336.center);
		\draw [style=d1] (269.center) to (339.center);
		\draw [style=d1] (275.center) to (340.center);
		\draw [style=d1] (255.center) to (337.center);
		\draw [style=d1] (261.center) to (338.center);
		\draw [style=d1] (179.center) to (341.center);
		\draw [style=d1] (343.center) to (345.center);
		\draw [style=d1] (344.center) to (346.center);
	\end{pgfonlayer}
\end{tikzpicture}

%% file: sections/appendix_grammar.tex
\section{Path Co-Width: An Alternative Definition via a Grammar}\label{appendix: grammar}
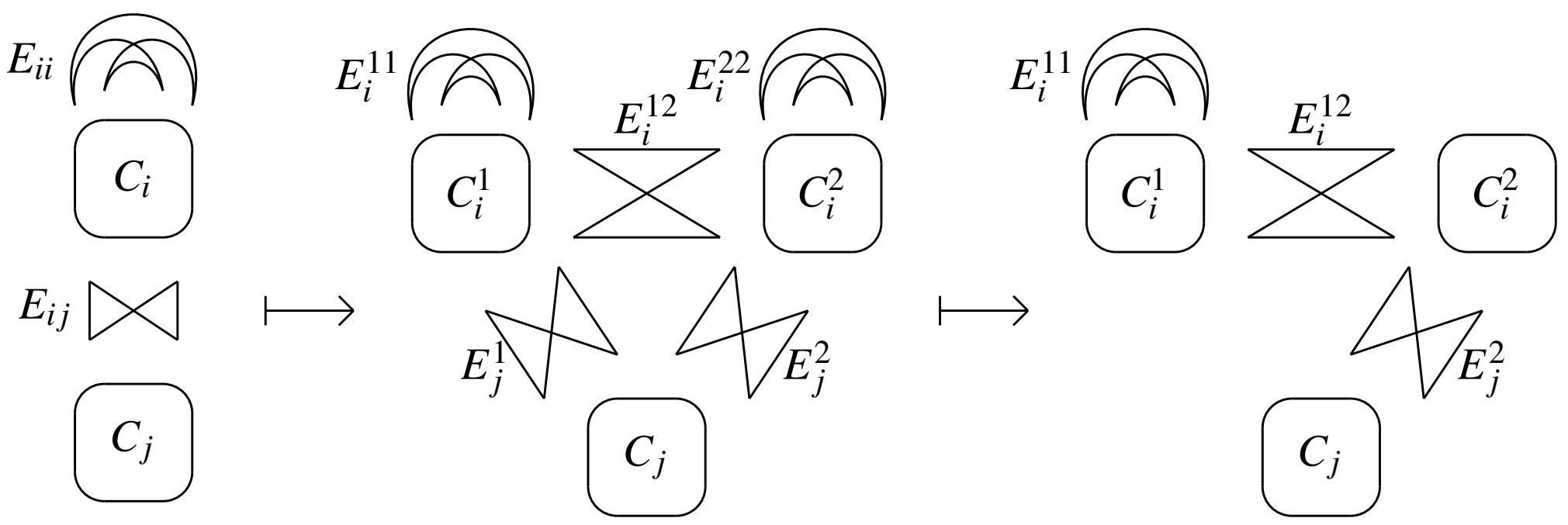
\begin{figure}
    \centering
    \input{grammar.tikz}
    \caption{On the left we have a colored graph $(G, \lambda)$ with two color classes $i$ and $j$. We will apply the color splitting operation at $i$ with matrix 
    \(
    A(i) = \bigl(\begin{smallmatrix}
    1 & 1 & 1 \\
    0 & 1 & 0
    \end{smallmatrix}\bigr)
    \)
    where rows are indexed by $C_i^1$ and $C_i^2$ and columns are indexed by $C_j$, $C_i^1$ and $C_i^2$. We obtain a graph -- the one on the far right -- in two steps: (1) first we duplicate the color class $i$ and copying all edges and (2) we then delete edges between color classes whenever the matrix $A(i)$ has a zero entry.}
    \label{fig:grammar}
\end{figure}

In what follows, we reformulate the definition of path co-width for \textit{symmetric graphs} so that it is easier to analyse inductively.
For $k\in\mathbb{N}$, a \define{$(\le k)$-coloured graph} is a graph $G$ along with an assignment $\lambda \colon V(G) \to \{1, \dots, k\}$ of colours to the vertices (these are not proper colourings). We define the following two operations on $(\le k)$-coloured graphs -- the ideas are simple, but aiming for maximum clarity, slightly verbose, bare with us: 
\begin{enumerate}
    \item \define{Colour merging.} Picking colors $i$ and $j$, define the $(\le k)$-coloured graph $(i \to j)(G, \lambda) := (G, \lambda')$ where $\lambda'$ recolors colour class $i$ to colour class $j$; formally we have
    \[\lambda' \colon x \mapsto \begin{cases} j &\text{if } \lambda(x) = i \\ 
    \lambda(x) &\text{otherwise}.\end{cases}\]
    \item \define{Class splitting.} This operation can only be performed when $G$ is colored with exactly $k'$ colors where $k' < k$. The rough idea is that, picking a color $i$, we will duplicate this color class and then join edges between these copies to the other color classes according to some simple rules (see Figure~\ref{fig:grammar} to get a feeling). Since this involves keeping track of various adjacencies, we will end up storing this data as a matrix. We need some notation. Denote the induced subgraph on $\lambda(i)^{-1}$ as $C_i$. Given a choice of a color $i$, a \define{class-$i$ splitting matrix} $A(i)$ is a binary matrix whose rows and columns are indexed as 
    \[
        \begin{array}{c|cccc|cc|c}
            & C_1 & C_2 & \cdots & C_{i-1} & C_i^1 & C_i^2 & C_{i+1} \; \cdots \; C_{k'} \\
            \hline
            C_i^1   &  &  &  &  & E_{i}^{11} & E_i^{12} &  \\
            C_i^2  &  &  &  &  & E_i^{21} & E_i^{22} &  
        \end{array}
    \]
    and satisfying the following conditions: 
    \begin{enumerate}
        \item every row indexed by some $j \neq i$ is a non-zero vector and 
        \item the block matrix  
            \(
                \begin{pmatrix}
                    \cdots & \big| & E_{i}^{11} & E_i^{12} & \big| & \cdots \\
                    \cdots & \big| & E_i^{21} & E_i^{22} & \big| & \cdots
                \end{pmatrix}
            \) is non-zero and $E_i^{21} = E_i^{12}$. 
    \end{enumerate}
    Given this matrix, we can finally define the class splitting operation, which we denote $A(i)(G, \lambda)$. This is the graph obtained as follows (see Figure~\ref{fig:grammar}):
    \begin{itemize}
        \item let $G'$ be the graph obtained by replacing $C_i$ with two copies $C_i^1$ and $C_i^2$ having the same adjacencies both within and outwith their color class  
        \item for each $h \in \{1, 2\}$, delete all edges between $C_i^h$ and $C_j$ if the         \item for each $h \in \{1, 2\}$, delete all edges between $C_i^h$ and $C_j$ if the $(C_i^h, C_j)$-th entry of $A(i)$ is zero; similarly delete all edges between $C_i^h$ and $C_i^g$ -- where $g \in \{1, 2\}$ -- if the $(C_i^h, C_i^g)$-th entry of $A(i)$ is zero (if $g = h$, then this amounts to removing all edges within the induced subgraph $C_i^h$).
        \item we obtain a $k' + 1$ coloring by adding a new color class for $C_{i}^2$ and leaving all other color classes the same (in particular $C_{i}^1$ retains color $i$).
    \end{itemize}
    \end{enumerate}
Define the \define{class of $(\le k)$-coloured graphs $\mathbb{C}_k$} recursively as follows:
\begin{enumerate}
    \item every $(\le k)$-coloured graph on $\le k$ vertices is in $\mathbb{C}_k$; and
    \item if $(G,\lambda)\in\mathbb{C}_k$, then the graph obtained from $(G,\lambda)$ by merging colour classes is in $\mathbb{C}_k$; and
    \item if $(G,\lambda)\in\mathbb{C}_k$ and $(G,\lambda)$ uses strictly less than $k$ colours, then the graph obtained from $(G,\lambda)$ by class splitting is in $\mathbb{C}_k$.
\end{enumerate}

We will prove that the class $\mathbb{C}_k$ defined by the grammar above captures precisely the graphs of path co-width at most $k$. To do this, we will first need the notion of \textit{gradual} tree co-decompositions which are roughly tree co-decompositions in which each cospan is making only small changes. 

\begin{definition}\label{def: gradual tree co-decomp}
    We say that a tree co-decomposition $\Ta = (T,\{G_a\},\{g_t^e\})$ of a graph $G$ is \define{gradual} if for every cospan $G_i \xrightarrow{g_i^e} G_e \xleftarrow{g_j^e} G_j$ we have that: (1) the difference in the number of vertices of $G_i$ and $G_j$ is at most one and (2) at least one of $g_i^e$ and $g_j^e$ is an identity.
\end{definition}

\begin{lemma}\label{lemm:helper}
    Let $k\in \mathbb{N}$ and $\Ta=(T,\{G_a\},\{g_t^e\})$ be a tree co-decomposition with limit~$G$,  width~$k$ and $n$ bags. Then there exists a gradual tree co-decomposition~$\Ta'$ such that: (1) it has $G$ as its limit (2) it has width at most $k$ and (3) it has at most $2kn$ bags. 
\end{lemma}
\begin{proof}
    Let $G_i\xrightarrow{g_i^e}G_e\xleftarrow{g_j^e} G_j$ be one of the cospans in $\Pa$. Suppose that $G_i$ has at least two vertices more than $G_e$. We observe that we can replace the cospan by the double cospan
    \[
        G_i\xrightarrow{{g_i^e}'} G_i'\xleftarrow{\cat{id}} {G_i}'\xrightarrow{{g_i^e}''}G_e\xleftarrow{g_j^e}G_j
    \]
    where ${g_i^e}'$ performs the identification of a single pair of vertices and ${g_i^e}''$ is such that $g_i^e={g_i^e}''\circ {g_i^e}'$. This new path co-decomposition retains the limit $G$, since pullbacks along identities are identities and thus the original path co-decomposition can be recovered as a partial evaluation at the cospan $G_i\xrightarrow{{g_i^e}'} G_i'\xleftarrow{\cat{id}} G_i'$. Since, by this reasoning, we need to replace each cospan with at most $2k$ cospans, the bound on the number of bags follows.
\end{proof}

We will now show that the graphs obtained our grammar are equivalent to graphs given by gradual path co-decompositions. 

\begin{lemma}
    For $k\in\mathbb{N}$, the class of graphs of path co-width $k$ is the class of underlying (i.e. uncolored) graphs of the graphs in $\mathbb{C}_k$.
\end{lemma}
\begin{proof}
    Let $\Pa=(P,\{G_a\},\{g_t^e\})$ be a path co-decomposition of width $k$. By Lemma~\ref{lemm:helper}, we may assume that $\Pa$ is gradual. Label the vertices of $P$ by $1,\ldots,n$ in the path order. Let $G[i]$ be the limit of the restriction of $\Pa$ to the first $i$ bags and \[g[i] \colon G[i] \twoheadrightarrow G_{i,i+1}\] be the epimorphism produced by the partial evaluation process -- in the diagram below, this is the composite of the $g_i$ (blue) and $g_i^{i+1}$ (orange) arrows below. 
   \begin{center}
\begin{adjustbox}{width=\textwidth}
\begin{tikzcd}[sep=small]
    &&&&& {G[i+1]} &&&&&& \\
    &&& {G[i]} \\
    \\
    {G_1} && {G_2} & \dots & {G_{i-1}} && {G_i} && {G_{i+1}} & \dots & \dots & {G_n} \\
    & {G_{1,2}} && \dots && {G_{i-1, i}} && {G_{i, i+1}} && {G_{i+1, i+2}} & \dots
    \arrow[from=1-6, to=2-4]
    \arrow["{g_{i+1}}", color={rgb,255:red,0;green,153;blue,64}, from=1-6, to=4-9]
    \arrow[color={rgb,255:red,92;green,92;blue,214}, two heads, from=2-4, to=4-1]
    \arrow[color={rgb,255:red,92;green,92;blue,214}, two heads, from=2-4, to=4-3]
    \arrow[color={rgb,255:red,92;green,92;blue,214}, two heads, from=2-4, to=4-5]
    \arrow["{g_i}", color={rgb,255:red,92;green,92;blue,214}, two heads, from=2-4, to=4-7]
    \arrow[two heads, from=4-1, to=5-2]
    \arrow[two heads, from=4-3, to=5-2]
    \arrow[two heads, from=4-3, to=5-4]
    \arrow[two heads, from=4-5, to=5-4]
    \arrow[two heads, from=4-5, to=5-6]
    \arrow[two heads, from=4-7, to=5-6]
    \arrow["{g_i^{i,i+1}}", color={rgb,255:red,214;green,153;blue,92}, two heads, from=4-7, to=5-8]
    \arrow[two heads, from=4-9, to=5-8]
    \arrow["{g_i^{i+1, i+2}}", color={rgb,255:red,120;green,48;blue,192}, from=4-9, to=5-10]
    \arrow[from=4-12, to=5-11]
\end{tikzcd}
\end{adjustbox}
\end{center}
    Let $\Pa[i]=(P,\{G[i],G_{i+1},G_{i+2},\ldots,G_n\},\{g[i]\}\cup\{g_t^e\})$ be the partial evaluation of $\Pa$ obtained by contracting the first $i-1$ edges in the auxiliary path $P$. 
    
    For each $i\in [n-1]$, let $(G[i],g[i])$ be the $(\le k)$-coloured graph obtained from $G[i]$ by colouring each vertex by its image under $g[i]$. Let $g[n]$ simply be the projection of $G[n]=G$ into $G_n$ (as if there is an additional adhesion map $g_{n}^{n,n+1}=\cat{id}_{G_n}$). We will prove by induction on $i\in[n]$ that $(G[i],g[i])$ is in the class $\mathbb{C}_k$.

    For the base case, we observe that $G[1]=G_1$ which has at most $k$ vertices. Hence the $(\le k)$-coloured graph $(G[1],g[1])$ is in $\mathbb{C}_k$ by definition. Now assume that $(G[i],g[i])\in \mathbb{C}_k$ for some $i\in[n]$.

    Suppose first that the adhesion $G_{i,i+1}$ and the bag $G_{i+1}$ are the same size, and hence are isomorphic. Then $G[i+1]=G[i]$ and $(G[i+1],g[i+1])$ is the $(\le k)$-coloured graph obtained from $(G[i],g[i])$ by colour merging according the the composite $g_{i+1}^{i+1,i+2}\circ g[i]$.
    
    Now suppose that $G_{i+1}$ has one more vertex than $G_i$. For the vertices of $G_{i+1}$, we write $v_1,\ldots,v_{m-1},v_m^1,v_m^2$ where $v_m^1$ and $v_m^2$ is the pair mapping to the same vertex in $G_{i,i+1}$.
    We claim that $(G[i+1],g[i+1])$ is the $(\le k)$-coloured graph obtained from $(G[i],g[i])$ by a class splitting, followed by a colour merging analogous to the one in the previous case. In particular, the class splitting is given by the matrix 
    \[
    \begin{array}{c|cccc|cc}
        & C_1 & C_2 & \cdots & C_{m-1} & C_m^1 & C_m^2  \\
        \hline
        C_m^1   &  &  & E_j^{1} &  & E_{m}^{11} & E_m^{12}  \\
        C_m^2  &  &  & E_j^{2}  &  & E_m^{21} & E_m^{22} 
    \end{array}
    \]
    where $E_j^r$ is $1$ if and only if either there is no edge between the images of $v_m^r$ and $v_j$ in $G_i=G_{i,i+1}$ or there is an edge between $v_m^r$ and $v_j$ in $G_{i+1}$; and for each $pq=11,12,21,22$ we set $E_m^{pq}=1$ if and only if there is an edge between $v_m^p$ and $v_m^q$ in $G_{i+1}$ (which might be the same vertex). We observe that the `non-zero' conditions (a) and (b) are satisfied because epimorphism are surjective as edge maps.

    For the converse direction, we will prove the following statement on the number of operations (colour merging and class splitting) $n\in\mathbb{N}$ needed to obtain our $(\le k)$-coloured graph: If $(G,\lambda)\in\mathbb{C}_k$ is constructed with $n$ operations then there is a path co-decomposition $\Pa=(P,\{G_a\},\{g_t^e\})$ so that $G$ is the limit of $\Pa$ and $\lambda$ is the vertex-map of the projection from $G$ into the last bag $G_n$ of $\Pa$ in the path order.

    Suppose that $(G,\lambda')$ is obtained from $(G,\lambda)$ by colour merging. Let $G_n'$ denote the graph obtained from $G_n$ my identifying the vertices representing the colours merged. Let $g_n':G_n\twoheadrightarrow G_n'$ denote the corresponding epimorphism. Let $\Pa'$ denote the path co-decomposition obtained from $\Pa$ by appending the cospan $G_n\xtwoheadrightarrow{g_n'} G_n'\xtwoheadleftarrow{\cat{id}} G_n'$ to the end of $\Pa$. Then $\Pa'$ has $G$ as its limit and the colouring given by the projection to the last graph is $\lambda'$, as required.

    Now suppose that $(G',\lambda')$ is obtained from $(G,\lambda)$ by class splitting with the matrix
    \[
        \begin{array}{c|cccc|cc|c}
            & C_1 & C_2 & \cdots & C_{i-1} & C_i^1 & C_i^2 & C_{i+1} \; \cdots \; C_{k'} \\
            \hline
            C_i^1   &  &  &  &  & E_{i}^{11} & E_i^{12} &  \\
            C_i^2  &  &  &  &  & E_i^{21} & E_i^{22} &  
        \end{array}
    \]
    Label the vertices of $G_i$ by $v_1,\ldots,v_m$, so that they correspond with the colours $1,\ldots,m$. Let $G_{n+1}$ be the graph on vertices $v_1,\ldots,v_{i-1},v_i^1,v_i^2,v_{i+1},\ldots,v_m$ so that it induces $G_n-v_i$ on the vertices aside from $v_i^1$ and $v_i^2$ and so that there is an edge from $v_j$ to $v_i^r$ if and only if matrix entry indexed by $C_j$ and $C_i^r$ is $1$; and for $pq=11,12,21,22$ there is an edge between $v_i^p$ and $v_i^q$ if and only if $E_i^{pq}=1$. Let $g_{n+1}^{n,n+1}:G_n\twoheadleftarrow G_{n+1}$ be the homomorphism between the two graphs. By the non-zero conditions on the matrix ((a) and (b)), we get that the defined homomorphism is an epimorphism. Now let $\Pa'$ be obtained from $\Pa$ by appending on the end the cospan
    \[
    G_n\xtwoheadrightarrow{\cat{id}} G_n\xtwoheadleftarrow{g_{n+1}^{n,n+1}} G_{n+1}
    \]
    and observe that $\Pa'$ has $G'$ as its limit and the projection into $G_{n+1}$ can be interpreted as the colouring~$\lambda'$.
\end{proof}

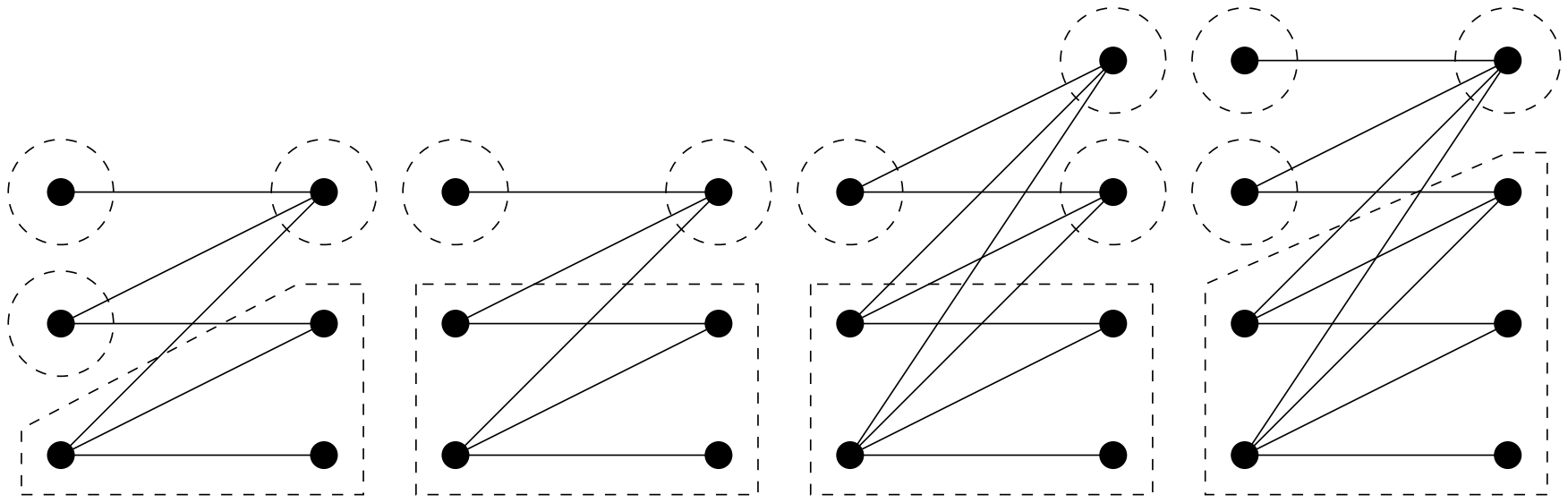
\begin{figure}
  \centering
  \resizebox{0.7\linewidth}{!}{%
    \begin{tikzpicture}
      \tikzstyle{every node}=[font=\fontsize{14.4pt}{18.7pt}\selectfont]
      \draw [ fill={rgb,255:red,0; green,0; blue,0}, fill opacity=1] (7.5,-7) circle (0.125cm);
      \draw [ fill={rgb,255:red,0; green,0; blue,0}, fill opacity=1] (10,-7) circle (0.125cm);
      \draw [ fill={rgb,255:red,0; green,0; blue,0}, fill opacity=1] (10,-8.25) circle (0.125cm);
      \draw [ fill={rgb,255:red,0; green,0; blue,0}, fill opacity=1] (7.5,-5.75) circle (0.125cm);
      \draw [ fill={rgb,255:red,0; green,0; blue,0}, fill opacity=1] (10,-5.75) circle (0.125cm);
      \draw [ fill={rgb,255:red,0; green,0; blue,0}, fill opacity=1] (7.5,-8.25) circle (0.125cm);
      \draw [ fill={rgb,255:red,0; green,0; blue,0}, fill opacity=1] (11.25,-7) circle (0.125cm);
      \draw [ fill={rgb,255:red,0; green,0; blue,0}, fill opacity=1] (13.75,-7) circle (0.125cm);
      \draw [ fill={rgb,255:red,0; green,0; blue,0}, fill opacity=1] (13.75,-8.25) circle (0.125cm);
      \draw [ fill={rgb,255:red,0; green,0; blue,0}, fill opacity=1] (11.25,-5.75) circle (0.125cm);
      \draw [ fill={rgb,255:red,0; green,0; blue,0}, fill opacity=1] (13.75,-5.75) circle (0.125cm);
      \draw [ fill={rgb,255:red,0; green,0; blue,0}, fill opacity=1] (11.25,-8.25) circle (0.125cm);
      \draw [ fill={rgb,255:red,0; green,0; blue,0}, fill opacity=1] (15,-7) circle (0.125cm);
      \draw [ fill={rgb,255:red,0; green,0; blue,0}, fill opacity=1] (17.5,-7) circle (0.125cm);
      \draw [ fill={rgb,255:red,0; green,0; blue,0}, fill opacity=1] (17.5,-8.25) circle (0.125cm);
      \draw [ fill={rgb,255:red,0; green,0; blue,0}, fill opacity=1] (15,-5.75) circle (0.125cm);
      \draw [ fill={rgb,255:red,0; green,0; blue,0}, fill opacity=1] (17.5,-5.75) circle (0.125cm);
      \draw [ fill={rgb,255:red,0; green,0; blue,0}, fill opacity=1] (15,-8.25) circle (0.125cm);
      \draw [ fill={rgb,255:red,0; green,0; blue,0}, fill opacity=1] (18.75,-7) circle (0.125cm);
      \draw [ fill={rgb,255:red,0; green,0; blue,0}, fill opacity=1] (21.25,-7) circle (0.125cm);
      \draw [ fill={rgb,255:red,0; green,0; blue,0}, fill opacity=1] (21.25,-8.25) circle (0.125cm);
      \draw [ fill={rgb,255:red,0; green,0; blue,0}, fill opacity=1] (18.75,-5.75) circle (0.125cm);
      \draw [ fill={rgb,255:red,0; green,0; blue,0}, fill opacity=1] (21.25,-5.75) circle (0.125cm);
      \draw [ fill={rgb,255:red,0; green,0; blue,0}, fill opacity=1] (18.75,-8.25) circle (0.125cm);
      \draw [ fill={rgb,255:red,0; green,0; blue,0}, fill opacity=1] (18.75,-4.5) circle (0.125cm);
      \draw [ fill={rgb,255:red,0; green,0; blue,0}, fill opacity=1] (21.25,-4.5) circle (0.125cm);
      \draw [ fill={rgb,255:red,0; green,0; blue,0}, fill opacity=1] (17.5,-4.5) circle (0.125cm);
      \draw  (7.5,-5.75) -- (10,-5.75);
      \draw  (7.5,-8.25) -- (10,-5.75);
      \draw  (7.5,-8.25) -- (10,-8.25);
      \draw  (10,-7) -- (7.5,-7);
      \draw  (7.5,-7) -- (10,-5.75);
      \draw  (10,-7) -- (7.5,-8.25);
      \draw  (11.25,-5.75) -- (13.75,-5.75);
      \draw  (13.75,-5.75) -- (11.25,-8.25);
      \draw  (11.25,-8.25) -- (13.75,-7);
      \draw  (13.75,-7) -- (11.25,-7);
      \draw  (11.25,-7) -- (13.75,-5.75);
      \draw  (13.75,-8.25) -- (11.25,-8.25);
      \draw  (15,-5.75) -- (17.5,-5.75);
      \draw  (17.5,-7) -- (15,-7);
      \draw  (15,-8.25) -- (17.5,-8.25);
      \draw  (17.5,-4.5) -- (15,-5.75);
      \draw  (17.5,-4.5) -- (15,-7);
      \draw  (17.5,-4.5) -- (15,-8.25);
      \draw  (17.5,-5.75) -- (15,-7);
      \draw  (17.5,-5.75) -- (15,-8.25);
      \draw  (17.5,-7) -- (15,-8.25);
      \draw  (18.75,-4.5) -- (21.25,-4.5);
      \draw  (21.1616,-5.6616) -- (21.1616,-5.6616);
      \draw  (21.25,-5.75) -- (18.625,-5.75);
      \draw  (18.75,-7) -- (21.25,-7);
      \draw  (21.25,-8.25) -- (18.75,-8.25);
      \draw  (18.75,-8.25) -- (21.25,-7);
      \draw  (18.75,-7) -- (21.25,-5.75);
      \draw  (18.75,-5.75) -- (21.25,-4.5);
      \draw  (21.25,-4.5) -- (18.75,-7);
      \draw  (21.25,-5.75) -- (18.75,-8.25);
      \draw  (21.25,-4.5) -- (18.75,-8.25);
      \draw [ dashed] (7.5,-5.75) circle (0.5cm);
      \draw [ dashed] (12,-4.5) circle (0cm);
      \draw [ dashed] (7.5,-5.75) circle (0cm);
      \draw [ dashed] (7.5,-7) circle (0.5cm);
      \draw [ dashed] (10,-5.75) circle (0.5cm);
      \draw [ dashed] (10.875,-6.625) rectangle  (14.125,-8.625);
      \draw [ , dashed] (7.125,-8.625) -- (10.375,-8.625) -- (10.375,-6.625) -- (9.75,-6.625) -- (7.125,-8) -- (7.125,-8.625);
      \draw [ dashed] (11.25,-5.75) circle (0.5cm);
      \draw [ dashed] (13.75,-5.75) circle (0.5cm);
      \draw [ dashed] (17.5,-4.5) circle (0.5cm);
      \draw [ dashed] (17.5,-5.75) circle (0.5cm);
      \draw [ dashed] (15,-5.75) circle (0.5cm);
      \draw [ dashed] (14.625,-6.625) rectangle  (17.875,-8.625);
      \draw [ , dashed] (21.625,-8.625) -- (18.375,-8.625) -- (18.375,-6.625) -- (21.25,-5.375) -- (21.625,-5.375) -- (21.625,-8.625);
      \draw [ dashed] (18.75,-5.75) circle (0.5cm);
      \draw [ dashed] (18.75,-4.5) circle (0.5cm);
      \draw [ dashed] (21.25,-4.5) circle (0.5cm);
    \end{tikzpicture}
  }%
  \caption{The process for building a half-graph}
  \label{fig:half_builder}
\end{figure}

\begin{example}\label{eg:half_graphs}
    For $n\in\mathbb{N}$, the \define{half graph $H_n$} is the bipartite graph with parts $u_1,\ldots,u_n$ and $v_1,\ldots,v_n$ and an edge $(u_i,v_j)$ whenever $i\leq j$. We will use the above grammar to construct all half graphs with a path co-decomposition with width at most $4$. For $n\ge 4$, let $H_n^1$ be the $(\le 4)$-coloured graph whose underlying graph is $H_n$ and whose colour classes are the singletons $A:=\{u_{n}\}$, $B:=\{u_{n-1}\}$ and $C:=\{v_n\}$ and the set of all other vertices $D:=\{u_1,\ldots,u_{n-2},v_1,\ldots,v_{n-1}\}$. Let $H_n^2$ be the $(\le 4)$-coloured graph obtained by merging the colour classes $B$ and $D$ to a new class $BD$. Let $H_n^2$ be the $(\le 4)$-coloured graph obtained by splitting the colour class $B$ to create a twin of $v_n$: label the new vertex $v_{n+1}$ and the new class $E:=\{v_{n+1}\}$. Let $H_n^3$ be the $(\le 4)$-coloured graph obtained from $H_n^2$ by merging the colour classes $BD$ and $C$ to a new class $BDC$. Let $H_n^4$ be the $(\le 4)$-coloured graph obtained by splitting $A$ to create a new vertex $u_{n+1}$ that is only incident with the vertices in $E$; that is, only with $v_{n+1}$. Observe that $H_n^4=H_{n+1}^1$. See Figure~\ref{fig:half_builder}.
\end{example}

%% file: figures/grammar.tikz
\begin{tikzpicture}
	\begin{pgfonlayer}{nodelayer}
		\node [style=none] (99) at (-1.25, 2.75) {};
		\node [style=none] (100) at (-0.75, 3.25) {};
		\node [style=none] (101) at (-1.25, 1.75) {};
		\node [style=none] (102) at (-0.75, 1.25) {};
		\node [style=none] (103) at (0.25, 1.25) {};
		\node [style=none] (104) at (0.75, 1.75) {};
		\node [style=none] (105) at (0.75, 2.75) {};
		\node [style=none] (106) at (0.25, 3.25) {};
		\node [style=none] (107) at (-0.25, 2.25) {$C_i^1$};
		\node [style=none] (108) at (4.75, 2.75) {};
		\node [style=none] (109) at (5.25, 3.25) {};
		\node [style=none] (110) at (4.75, 1.75) {};
		\node [style=none] (111) at (5.25, 1.25) {};
		\node [style=none] (112) at (6.25, 1.25) {};
		\node [style=none] (113) at (6.75, 1.75) {};
		\node [style=none] (114) at (6.75, 2.75) {};
		\node [style=none] (115) at (6.25, 3.25) {};
		\node [style=none] (116) at (5.75, 2.25) {$C_i^2$};
		\node [style=none] (117) at (1.75, -1.75) {};
		\node [style=none] (118) at (2.25, -1.25) {};
		\node [style=none] (119) at (1.75, -2.75) {};
		\node [style=none] (120) at (2.25, -3.25) {};
		\node [style=none] (121) at (3.25, -3.25) {};
		\node [style=none] (122) at (3.75, -2.75) {};
		\node [style=none] (123) at (3.75, -1.75) {};
		\node [style=none] (124) at (3.25, -1.25) {};
		\node [style=none] (125) at (2.75, -2.25) {$C_j$};
		\node [style=none] (126) at (1.5, 3) {};
		\node [style=none] (127) at (1.5, 1.5) {};
		\node [style=none] (128) at (4, 3) {};
		\node [style=none] (129) at (4, 1.5) {};
		\node [style=none] (130) at (2.75, 3.5) {$E_{i}^{12}$};
		\node [style=none] (131) at (1.25, 1) {};
		\node [style=none] (132) at (0, 0.25) {};
		\node [style=none] (133) at (1, -1.25) {};
		\node [style=none] (134) at (2.25, -0.5) {};
		\node [style=none] (135) at (4.25, 1) {};
		\node [style=none] (136) at (5.5, 0.25) {};
		\node [style=none] (137) at (4.5, -1.25) {};
		\node [style=none] (138) at (3.25, -0.5) {};
		\node [style=none] (139) at (0, -0.75) {$E_j^{1}$};
		\node [style=none] (140) at (5.5, -0.75) {$E_{j}^2$};
		\node [style=none] (141) at (-7, -2.5) {};
		\node [style=none] (142) at (-6.5, -3) {};
		\node [style=none] (143) at (-7, -1.5) {};
		\node [style=none] (144) at (-6.5, -1) {};
		\node [style=none] (145) at (-5.5, -1) {};
		\node [style=none] (146) at (-5, -1.5) {};
		\node [style=none] (147) at (-5, -2.5) {};
		\node [style=none] (148) at (-5.5, -3) {};
		\node [style=none] (149) at (-6, -2) {$C_j$};
		\node [style=none] (150) at (-7, 2) {};
		\node [style=none] (151) at (-6.5, 1.5) {};
		\node [style=none] (152) at (-7, 3) {};
		\node [style=none] (153) at (-6.5, 3.5) {};
		\node [style=none] (154) at (-5.5, 3.5) {};
		\node [style=none] (155) at (-5, 3) {};
		\node [style=none] (156) at (-5, 2) {};
		\node [style=none] (157) at (-5.5, 1.5) {};
		\node [style=none] (158) at (-6, 2.5) {$C_i$};
		\node [style=none] (159) at (-5.25, -0.25) {};
		\node [style=none] (160) at (-6.75, -0.25) {};
		\node [style=none] (161) at (-6.75, 0.75) {};
		\node [style=none] (162) at (-5.25, 0.75) {};
		\node [style=none] (163) at (-3.75, 0.25) {};
		\node [style=none] (164) at (-2.25, 0.25) {};
		\node [style=none] (165) at (-3.75, 0) {};
		\node [style=none] (166) at (-3.75, 0.5) {};
		\node [style=none] (167) at (-2.5, 0) {};
		\node [style=none] (168) at (-2.5, 0.5) {};
		\node [style=none] (169) at (-7.5, 0.25) {$E_{ij}$};
		\node [style=none] (170) at (-7, 3.75) {};
		\node [style=none] (171) at (-5, 3.75) {};
		\node [style=none] (172) at (-6.5, 4) {};
		\node [style=none] (173) at (-5.5, 4) {};
		\node [style=none] (174) at (-7.75, 4.5) {$E_{ii}$};
		\node [style=none] (175) at (-1.25, 3.5) {};
		\node [style=none] (176) at (0.75, 3.5) {};
		\node [style=none] (177) at (-0.75, 3.75) {};
		\node [style=none] (178) at (0.25, 3.75) {};
		\node [style=none] (179) at (-2, 4.25) {$E_{i}^{11}$};
		\node [style=none] (180) at (4.75, 3.5) {};
		\node [style=none] (181) at (6.75, 3.5) {};
		\node [style=none] (182) at (5.25, 3.75) {};
		\node [style=none] (183) at (6.25, 3.75) {};
		\node [style=none] (184) at (4, 4.25) {$E_{i}^{22}$};
		\node [style=none] (185) at (10.25, 2.75) {};
		\node [style=none] (186) at (10.75, 3.25) {};
		\node [style=none] (187) at (10.25, 1.75) {};
		\node [style=none] (188) at (10.75, 1.25) {};
		\node [style=none] (189) at (11.75, 1.25) {};
		\node [style=none] (190) at (12.25, 1.75) {};
		\node [style=none] (191) at (12.25, 2.75) {};
		\node [style=none] (192) at (11.75, 3.25) {};
		\node [style=none] (193) at (11.25, 2.25) {$C_i^1$};
		\node [style=none] (194) at (16.25, 2.75) {};
		\node [style=none] (195) at (16.75, 3.25) {};
		\node [style=none] (196) at (16.25, 1.75) {};
		\node [style=none] (197) at (16.75, 1.25) {};
		\node [style=none] (198) at (17.75, 1.25) {};
		\node [style=none] (199) at (18.25, 1.75) {};
		\node [style=none] (200) at (18.25, 2.75) {};
		\node [style=none] (201) at (17.75, 3.25) {};
		\node [style=none] (202) at (17.25, 2.25) {$C_i^2$};
		\node [style=none] (203) at (13.25, -1.75) {};
		\node [style=none] (204) at (13.75, -1.25) {};
		\node [style=none] (205) at (13.25, -2.75) {};
		\node [style=none] (206) at (13.75, -3.25) {};
		\node [style=none] (207) at (14.75, -3.25) {};
		\node [style=none] (208) at (15.25, -2.75) {};
		\node [style=none] (209) at (15.25, -1.75) {};
		\node [style=none] (210) at (14.75, -1.25) {};
		\node [style=none] (211) at (14.25, -2.25) {$C_j$};
		\node [style=none] (212) at (13, 3) {};
		\node [style=none] (213) at (13, 1.5) {};
		\node [style=none] (214) at (15.5, 3) {};
		\node [style=none] (215) at (15.5, 1.5) {};
		\node [style=none] (216) at (14.25, 3.5) {$E_{i}^{12}$};
		\node [style=none] (221) at (15.75, 1) {};
		\node [style=none] (222) at (17, 0.25) {};
		\node [style=none] (223) at (16, -1.25) {};
		\node [style=none] (224) at (14.75, -0.5) {};
		\node [style=none] (226) at (17, -0.75) {$E_{j}^2$};
		\node [style=none] (227) at (7.75, 0.25) {};
		\node [style=none] (228) at (9.25, 0.25) {};
		\node [style=none] (229) at (7.75, 0) {};
		\node [style=none] (230) at (7.75, 0.5) {};
		\node [style=none] (231) at (9, 0) {};
		\node [style=none] (232) at (9, 0.5) {};
		\node [style=none] (233) at (10.25, 3.5) {};
		\node [style=none] (234) at (12.25, 3.5) {};
		\node [style=none] (235) at (10.75, 3.75) {};
		\node [style=none] (236) at (11.75, 3.75) {};
		\node [style=none] (237) at (9.5, 4.25) {$E_{i}^{11}$};
	\end{pgfonlayer}
	\begin{pgfonlayer}{edgelayer}
		\draw [style=c0] (101.center) to (99.center);
		\draw [style=c0, bend left=45] (99.center) to (100.center);
		\draw [style=c0] (100.center) to (106.center);
		\draw [style=c0, bend left=45] (106.center) to (105.center);
		\draw [style=c0] (105.center) to (104.center);
		\draw [style=c0, bend left=45] (104.center) to (103.center);
		\draw [style=c0] (103.center) to (102.center);
		\draw [style=c0, bend left=45] (102.center) to (101.center);
		\draw [style=c0] (110.center) to (108.center);
		\draw [style=c0, bend left=45] (108.center) to (109.center);
		\draw [style=c0] (109.center) to (115.center);
		\draw [style=c0, bend left=45] (115.center) to (114.center);
		\draw [style=c0] (114.center) to (113.center);
		\draw [style=c0, bend left=45] (113.center) to (112.center);
		\draw [style=c0] (112.center) to (111.center);
		\draw [style=c0, bend left=45] (111.center) to (110.center);
		\draw [style=c0] (119.center) to (117.center);
		\draw [style=c0, bend left=45] (117.center) to (118.center);
		\draw [style=c0] (118.center) to (124.center);
		\draw [style=c0, bend left=45] (124.center) to (123.center);
		\draw [style=c0] (123.center) to (122.center);
		\draw [style=c0, bend left=45] (122.center) to (121.center);
		\draw [style=c0] (121.center) to (120.center);
		\draw [style=c0, bend left=45] (120.center) to (119.center);
		\draw [style=c0] (126.center) to (128.center);
		\draw [style=c0] (129.center) to (127.center);
		\draw [style=c0] (126.center) to (129.center);
		\draw [style=c0] (128.center) to (127.center);
		\draw [style=c0] (131.center) to (134.center);
		\draw [style=c0] (133.center) to (132.center);
		\draw [style=c0] (132.center) to (134.center);
		\draw [style=c0] (133.center) to (131.center);
		\draw [style=c0] (135.center) to (138.center);
		\draw [style=c0] (137.center) to (136.center);
		\draw [style=c0] (136.center) to (138.center);
		\draw [style=c0] (137.center) to (135.center);
		\draw [style=c0] (143.center) to (141.center);
		\draw [style=c0, bend right=45] (141.center) to (142.center);
		\draw [style=c0] (142.center) to (148.center);
		\draw [style=c0, bend right=45] (148.center) to (147.center);
		\draw [style=c0] (147.center) to (146.center);
		\draw [style=c0, bend right=45] (146.center) to (145.center);
		\draw [style=c0] (145.center) to (144.center);
		\draw [style=c0, bend right=45] (144.center) to (143.center);
		\draw [style=c0] (152.center) to (150.center);
		\draw [style=c0, bend right=45] (150.center) to (151.center);
		\draw [style=c0] (151.center) to (157.center);
		\draw [style=c0, bend right=45] (157.center) to (156.center);
		\draw [style=c0] (156.center) to (155.center);
		\draw [style=c0, bend right=45] (155.center) to (154.center);
		\draw [style=c0] (154.center) to (153.center);
		\draw [style=c0, bend right=45] (153.center) to (152.center);
		\draw [style=c0] (159.center) to (162.center);
		\draw [style=c0] (161.center) to (160.center);
		\draw [style=c0] (160.center) to (162.center);
		\draw [style=c0] (161.center) to (159.center);
		\draw [style=c0] (165.center) to (166.center);
		\draw [style=c0] (163.center) to (164.center);
		\draw [style=c0] (164.center) to (167.center);
		\draw [style=c0] (164.center) to (168.center);
		\draw [style=c0, bend left=105, looseness=2.75] (170.center) to (171.center);
		\draw [style=c0, bend left=75, looseness=1.75] (172.center) to (173.center);
		\draw [style=c0, bend left=90, looseness=2.25] (170.center) to (173.center);
		\draw [style=c0, bend left=90, looseness=2.25] (172.center) to (171.center);
		\draw [style=c0, bend left=105, looseness=2.75] (175.center) to (176.center);
		\draw [style=c0, bend left=75, looseness=1.75] (177.center) to (178.center);
		\draw [style=c0, bend left=90, looseness=2.25] (175.center) to (178.center);
		\draw [style=c0, bend left=90, looseness=2.25] (177.center) to (176.center);
		\draw [style=c0, bend left=105, looseness=2.75] (180.center) to (181.center);
		\draw [style=c0, bend left=75, looseness=1.75] (182.center) to (183.center);
		\draw [style=c0, bend left=90, looseness=2.25] (180.center) to (183.center);
		\draw [style=c0, bend left=90, looseness=2.25] (182.center) to (181.center);
		\draw [style=c0] (187.center) to (185.center);
		\draw [style=c0, bend left=45] (185.center) to (186.center);
		\draw [style=c0] (186.center) to (192.center);
		\draw [style=c0, bend left=45] (192.center) to (191.center);
		\draw [style=c0] (191.center) to (190.center);
		\draw [style=c0, bend left=45] (190.center) to (189.center);
		\draw [style=c0] (189.center) to (188.center);
		\draw [style=c0, bend left=45] (188.center) to (187.center);
		\draw [style=c0] (196.center) to (194.center);
		\draw [style=c0, bend left=45] (194.center) to (195.center);
		\draw [style=c0] (195.center) to (201.center);
		\draw [style=c0, bend left=45] (201.center) to (200.center);
		\draw [style=c0] (200.center) to (199.center);
		\draw [style=c0, bend left=45] (199.center) to (198.center);
		\draw [style=c0] (198.center) to (197.center);
		\draw [style=c0, bend left=45] (197.center) to (196.center);
		\draw [style=c0] (205.center) to (203.center);
		\draw [style=c0, bend left=45] (203.center) to (204.center);
		\draw [style=c0] (204.center) to (210.center);
		\draw [style=c0, bend left=45] (210.center) to (209.center);
		\draw [style=c0] (209.center) to (208.center);
		\draw [style=c0, bend left=45] (208.center) to (207.center);
		\draw [style=c0] (207.center) to (206.center);
		\draw [style=c0, bend left=45] (206.center) to (205.center);
		\draw [style=c0] (212.center) to (214.center);
		\draw [style=c0] (215.center) to (213.center);
		\draw [style=c0] (212.center) to (215.center);
		\draw [style=c0] (214.center) to (213.center);
		\draw [style=c0] (221.center) to (224.center);
		\draw [style=c0] (223.center) to (222.center);
		\draw [style=c0] (222.center) to (224.center);
		\draw [style=c0] (223.center) to (221.center);
		\draw [style=c0] (229.center) to (230.center);
		\draw [style=c0] (227.center) to (228.center);
		\draw [style=c0] (228.center) to (231.center);
		\draw [style=c0] (228.center) to (232.center);
		\draw [style=c0, bend left=105, looseness=2.75] (233.center) to (234.center);
		\draw [style=c0, bend left=75, looseness=1.75] (235.center) to (236.center);
		\draw [style=c0, bend left=90, looseness=2.25] (233.center) to (236.center);
		\draw [style=c0, bend left=90, looseness=2.25] (235.center) to (234.center);
	\end{pgfonlayer}
\end{tikzpicture}